\documentclass{article}

\usepackage{color}
\usepackage{amssymb}
\usepackage{latexsym}
\usepackage{amsmath}
\usepackage{amsfonts}
\usepackage{bbm}
\usepackage{csquotes}
\usepackage{hyperref}
\usepackage{mathtools}

\def\reff#1{{\rm(\ref{#1})}}

\def\be{\begin{eqnarray}}
\def\ee{\end{eqnarray}}
\def\b*{\begin{eqnarray*}}
\def\e*{\end{eqnarray*}}

\newtheorem{Theorem}{Theorem}[section]
\newtheorem{Definition}[Theorem]{Definition}
\newtheorem{Proposition}[Theorem]{Proposition}
\newtheorem{Assumption}[Theorem]{Assumption}
\newtheorem{Lemma}[Theorem]{Lemma}

\newtheorem{Remark}[Theorem]{Remark}
\newtheorem{Example}[Theorem]{Example}
\makeatletter \@addtoreset{equation}{section}

\makeatother \makeatletter

\def\esup{\mathop{\rm ess\;sup}}
\def\einf{\mathop{\rm ess\;inf}}

\def \be{\begin{eqnarray}}
\def \ee{\end{eqnarray}}
\def \b*{\begin{eqnarray*}}
\def \e*{\end{eqnarray*}}
\def\no{\noindent}

\def \E{\mathbb{E}}
\def \F{\mathbb{F}}

\def \L{\mathbb{L}}

\def \P{\mathbb{P}}

\def \R{\mathbb{R}}
\def \S{\mathbb{S}}

\def \[{[\,\!\![}
\def \]{]\,\!\!]}

\def \1{{\bf 1}}
\def\esup{\mathop{\rm ess\;sup}}

\def \proof{{\noindent \bf Proof}\quad}
\def \ep{\hbox{ }\hfill$\Box$}

\def\reff#1{{\rm(\ref{#1})}}

\def\Ac{{\cal A}}
\def\Bc{{\cal B}}

\def\Ec{{\cal E}}
\def\Fc{{\cal F}}

\def\Hc{{\cal H}}
\def\Jc{{\cal J}}
\def\Lc{{\cal L}}

\def\Kc{{\cal K}}

\def\Lc{{\cal L}}

\def\Pc{{\cal P}}
\def\Sc{{\cal S}}
\def\Tc{{\cal T}}

\def\Yc{{\cal Y}}

\def\cE{{\cal E}}
\def\cF{{\cal F}}
\def\cP{{\cal P}}
\def\cT{{\cal T}}
\def\cJ{{\cal J}}
\def\cL{{\cal L}}
\def\dbE{{\mathbb E}}
\def\dbF{{\mathbb F}}
\def\dbP{{\mathbb P}}
\def\dbQ{{\mathbb Q}}
\def\dbR{{\mathbb R}}
\def\dbS{{\mathbb S}}

\def\pa{{\partial}}

\newcommand{\ol}{\overline}
\newcommand{\ul}{\underline}
\newcommand{\we}{\wedge}

\def\d{{\delta}}
\def\f{{\varphi}}
\def\g{{\gamma}}
\def\o{{\omega}}
\def\O{{\Omega}}
\def\t{{\tau}}
\def\l{{\lambda}}
\def\th{{\theta}}

\def\hc{\textsc{h}}

\def\eps{\varepsilon}

\def \E{\mathbb{E}}
\def \F{\mathbb{F}}

\def \R{\mathbb{R}}

\def\P{\mathbb{P}}

\def\T{\mathbb{T}}

\def\Ac{{\cal A}}
\def\Bc{{\cal B}}

\def\Ec{{\cal E}}
\def\Fc{{\cal F}}

\def\Hc{{\cal H}}

\def\Jc{{\cal J}}
\def\Kc{{\cal K}}
\def\Lc{{\cal L}}

\def\Pc{{\cal P}}
\def\Qc{{\cal Q}}

\def\Sc{{\cal S}}
\def\Tc{{\cal T}}

\def\Yc{{\cal Y}}

\def\ep{\hbox{ }\hfill$\Box$}
\def\reff#1{{\rm(\ref{#1})}}
\def\be{\begin{eqnarray}}
\def\ee{\end{eqnarray}}
\def\beq{\begin{equation}}
\def\eeq{\end{equation}}

\def\={\;=\;}

\def\no{\noindent}
\def\.{\;.}

\def\eps{\varepsilon}

\def\1{{\bf 1}}

\def\cA{{\mathcal{A}}}

\def\cP{{\mathcal{P}}}

\def\eqref#1{\reff{#1}}

\newcommand{\bea}{\begin{eqnarray}}
\newcommand{\eea}{\end{eqnarray}}
\newcommand{\bes}{\begin{subequations}}
\newcommand{\ees}{\end{subequations}}
\newcommand{\bgt}{\begin{gather}}
\newcommand{\egt}{\end{gather}}

\newcommand{\beaa}{\begin{eqnarray*}}
\newcommand{\eeaa}{\end{eqnarray*}}

\title{\bf An overview of Viscosity Solutions of Path-Dependent PDEs
             \thanks{Ren and Touzi's research is supported by the ERC grant 321111 RoFiRM, the ANR grant ISOTACE, the Chair {\it Financial Risks} of the {\it Risk Foundation} sponsored by Soci\'et\'e G\'en\'erale, and the Chair {\it Finance and Sustainable Development} sponsored by EDF and Calyon.  }}
\author{Zhen-Jie  {\sc Ren}
             \thanks{CMAP, Ecole Polytechnique Paris, zhenjie.ren@cmap.polytechnique.fr.}   
             \and 
             Nizar {\sc Touzi}
             \thanks{CMAP, Ecole Polytechnique Paris, nizar.touzi@polytechnique.edu.}
             \and 
             Jianfeng {\sc Zhang}
             \thanks{University of Southern California, Department of Mathematics, jianfenz@usc.edu.}
             }
\date{\today}

\begin{document}
\maketitle

\begin{abstract}
This paper provides an overview of the recently developed notion of viscosity solutions of path-dependent partial differential equations. We start by a quick review of the Crandall-Ishii notion of viscosity solutions, so as to motivate the relevance of our definition in the path-dependent case. We focus on the wellposedness theory of such equations. In particular, we provide a simple presentation of the current existence and uniqueness arguments in the semilinear case. We also review the stability property of this notion of solutions, including the adaptation of the Barles-Souganidis monotonic scheme approximation method. Our results rely crucially on the theory of optimal stopping under nonlinear expectation. In the dominated case, we provide a self-contained presentation of all required results. The fully nonlinear case is more involved and is addressed in \cite{ETZ0}.
\end{abstract}

\vspace{5mm}


\vspace{5mm}


\vfill\eject

\section{Introduction}
\label{sect-Introduction}
\setcounter{equation}{0}

Let $\Omega:=\{\omega\in C^0([0,T],\R^d):\omega_0=0\}$ be the canonical space of continuous paths starting from the origin, $B$ the canonical process defined by $B_t(\o):=\omega_t$, $t\in[0,T]$, and $\F:=\{\Fc_t, t\in[0,T]\}$ the corresponding filtration. Following Dupire \cite{Dupire}, we introduce the pseudo-distance
 \bea\label{Dupire-distance}
 d\big((t,\omega),(t',\omega')\big)
 :=
 |t-t'|+\|\omega_{\wedge t}-\omega'_{\wedge t'}\|_\infty
 &\mbox{for all}&
 t,t'\in[0,T],~\omega,\omega'\in\Omega.
 \eea
Then, any process $u:[0,T]\times\Omega\longrightarrow\R$, continuous with respect to $d$, is $\F-$progressively measurable, so that $u(t,\omega)=u\big(t,(\omega_s)_{s\le t}\big)$.

The goal of this paper is to provide a wellposedness theory for the path-dependent partial differential equation (PDE):
 \bea\label{PPDE}
 -\partial_tu(t,\omega)
 - G\big(t,\omega,u(t,\omega),\partial_\omega u(t,\omega), \partial^2_{\omega\omega}u(t,\omega)
      \big)
 =0,
 &t<T,&
 \omega\in\Omega.
 \eea
with boundary condition $u(T,\omega)=\xi(\omega)$. Here, $\xi:(\Omega,\Fc_T)\longrightarrow(\R,\Bc_\R)$ is a bounded uniformly continuous function, and $G:[0,T]\times\Omega\times\R\times\R^d\times\S_d\longrightarrow \R$ is continuous in $(t,\omega)$, Lipschitz-continuous in the remaining variables $(y,z,\gamma)$, and satisfies the ellipticity condition:
 \bea\label{elliptic}
 \gamma\in\S_d \longmapsto G(t,\omega,y,z,\gamma)
 &\mbox{is non-decreasing.}&
 \eea
The unknown process $u(t,\omega)$ is required to be $\F-$progressively measurable, and the derivatives $\partial_t u,\partial_\omega u,\partial^2_{\omega\omega}u$ are $\F-$progressively measurable processes valued in $\R,\R^d,\S_d$, respectively, which will be defined later. Notice in particular that, as $\R^d-$ and $\S_d-$valued process, the derivatives $\partial_\omega u, \partial^2_{\omega\omega}u$ do not correspond to some (infinite-dimensional) gradient and Hessian with respect to the path. Consequently, the equation \eqref{PPDE} is a PDE parameterized by the path, and not a general PDE on the paths space. For this reason, the name path-dependent PDE is more relevant than PDE on the paths space.
  
There are three particular examples of such equations which can be related to the existing probability theory literature, namely  
\begin{enumerate}
\item When the nonlinearity $G$ is linear: 
  \be\label{Glin}
  G^{\mbox{\footnotesize lin}}(.,y,z,\gamma)
  :=
  \ell-ky+\frac12 {\rm Tr}(\gamma),
  \ee
for some functions $\ell,k$ defined on $[0,T]\times\Omega$ 
the natural solution of the equation \eqref{PPDE} is given by any regular version of the conditional expectation
 \be
 u^{\mbox{\footnotesize lin}}(t,\omega)
 &:=&
 \E^{\P_0}\Big[\int_t^T e^{-\int_t^s k_rdr}\ell_s ds + e^{-\int_t^T k_rdr}\xi \Big|\Fc_t\Big](\omega),
 \ee
 where $\P_0$ is the Wiener measure. Similar results hold for more general linear equations.
\item When the nonlinearity $G$ is semilinear:
 \be\label{Gs-lin}
 G^{\mbox{\footnotesize s-lin}}(.,y,z,\gamma)
 :=
 \frac12 {\rm Tr}(\gamma)
 +F(.,y,z),
 \ee
for some function $F:[0,T]\times\Omega\times\R\times\R^d\longrightarrow\R$, the natural solution of the equation \eqref{PPDE} is given by any regular version of the backward stochastic differenttial equation:
 \b*
 u^{\mbox{\footnotesize s-lin}}(t,\omega)=Y_t(\omega)
 &\mbox{where}&
 Y_s
 =
 \xi+\int_s^T F_r(Y_r,Z_r)dr-\int_s^T Z_r dB_r,~\P_0-a.s.
 \e*
 
\item The theory of second order backward stochastic differential equations introduced in \cite{CSTV} and \cite{STZ} provides a similar representation of the natural solution of the path-dependent PDE \eqref{PPDE} for  a class of fully nonlinearities $G$. 
\end{enumerate}

Another important particular example, which plays the role of a benchmark, is the so-called Markovian case when $\xi(\omega)=h(\omega_T)$, and $G(t,\omega,y,z,\gamma)=g(t,\omega_t,y,z,\gamma)$ for some functions $g$ and $h$ defined on the corresponding finite-dimensional spaces. In this context, we expect that $u(t,\omega)=v(t,\omega_t)$ for some function $v:[0,T]\times\R^d\longrightarrow\R$, and the path-dependent PDE \eqref{PPDE} reduces to the standard PDE:
 \be\label{PDE}
 -\partial_tv(t,x)
 - g\big(t,x,v(t,x),Dv(t,x), D^2v(t,x)
      \big)
 =0,
 &t<T,&
 x\in\R^d,
 \ee
where $\partial_t,D,D^2$ denotes respectively the standard time derivative, the gradient and the Hessian with respect to the space variable. In this case, it is well-known that the theory of viscosity solutions introduced by Crandall and Lions \cite{CL1,CL2} is a powerful notion of weak solution for which a solid existence and uniqueness theory has been developed, and which proved its relevance for various applications. Viscosity solutions gained importance by the contributions of Barles and Souganidis \cite{BS} to the convergence of numerical schemes, and the work of Cafarelli and Cabre \cite{CC} which makes a crucial use of viscosity solutions to obtain sharp regularity results.

Our main concern is the adaptation of the notion of viscosity solutions to the context of our path-dependent PDE \eqref{PPDE}. However, the fact that our underlying space, namely $[0,T]\times\Omega$, is not locally compact raises a major difficulty which needs to be addressed. Indeed, the stability and the uniqueness results in the theory of viscosity solutions is based on the existence of a local maximizer for an arbitrary upper semicontinuous function. 

In order to by-pass this difficulty, we introduce a convenient modification of the definition. To explain our definition, let us focus on the notion of viscosity subsolution, the case of a viscosity supersolution is symmetric. For a viscosity subsolution $u$, the standard definition considers as test functions some point $(t_0,x_0)$ all those functions $\varphi$ which are pointwisely locally tangent from above to $u$ with contact point $(t_0,x_0)$:
 \b*
 (\varphi-u)(t_0,x_0)
 =
 \min_{O_r(t_0,x_0)} (\varphi-u),
 &\mbox{for some}&
 r>0,
 \e*
where $O_r(t_0,x_0)$ denotes the open ball in $\R^{d+1}$ centered at $(t_0,x_0)$, with radius $r$. \begin{enumerate}
\item For simplicity, we first consider the case of a nonlinearity $G=G^{\mbox{\footnotesize lin}}$ as in \eqref{Glin}, or $G=G^{\mbox{\footnotesize s-lin}}$ as in \eqref{Gs-lin}, with $F(t,\omega,y,z)$ independent of the $z-$component. Our definition follows exactly the spirit of viscosity solutions, but replaces the pointwise tangency by the corresponding notion in mean:
 \b*
 (\varphi-u)(t_0,\omega_0)
 =
 \min_{\tau} \E^\P\big[(\varphi-u)_{\tau\wedge{\hc}}|\Fc_{t_0}\big](\omega_0),
 &\mbox{for some stopping time}&{\hc}>t_0,
 \e*
 where the min is over all stopping times $\t\ge t_0$.
 
\item For a more general nonlinearity $G$, our definition replaces the expectation operator $\E^\P$ by a the sublinear expectation operator $\overline{\Ec}:=\sup_{\P\in\Pc}\E^\P$ for some convenient family $\Pc$ of probability measure. We observe that $\Pc$ can be chosen to be a dominated family of measures in the semilinear case $G=G^{\mbox{\footnotesize s-lin}}$. However, in the general nonlinear case, the family $\Pc$ is not dominated. 
\end{enumerate}

The main purpose of this paper is to provide an overview of the available results on the wellposedness of the path-dependent PDE under this notion of viscosity solution. In particular, we highlight that our definition induces a richer family of test function in the Markovian case. Consequently, 

(i) the existence may be more difficult to achieve under our definition; however, we shall see that the traditional examples from the applications raise no special difficulty from the existence side; in fact, in contrast with the standard notion of viscosity solution, our definition is tight,

(ii) the uniqueness may be easier under our definition because our notion of viscosity solution is constrained by a bigger set of test functions; indeed recently comparison results were obtained in the semilinear case $G=G^{\mbox{\footnotesize s-lin}}$ with relatively simple arguments avoiding the Crandall-Ishii's lemma of the standard viscosity solution in the Markovian case; in particular, the comparison result for the linear path-dependent PDE $G=G^{\mbox{\footnotesize lin}}$ follows from the equivalence between our notion of viscosity subsolution and the (regular) submartingale property, whose proof is a simple consequence of the theory of optimal stopping.

This paper also pays a special attention to the stability of our notion of viscosity solutions, which is an essential property of standard viscosity solutions in the Markovian case, and is responsible for the denomination of this notion. We shall present the present state of stability results, together with the corresponding convergence results of numerical schemes \`a la Barles \& Souganidis \cite{BS}. 

\section{Standard viscosity solution in the Markovian case}
\label{sect-standardvisco}
\setcounter{equation}{0}

In this short section, we recall the standard definition of viscosity solutions in the Markovian case, and we review the corresponding existence and uniqueness results. In order for our notations to be consistent with the path-dependent case, our functions will be defined on $\mbox{cl}(Q)=[0,T]\times\R^d$, where $Q:=[0,T)\times\R^d$. 

\subsection{Definitions and consistency with classical solutions}

For $(t,x)\in Q$, $u\in\mbox{\rm USC}(Q)$, and $v\in\mbox{\rm LSC}(Q)$, we denote:
 \be
 \underline{A}u(t,x)
 &:=&
 \big\{\varphi\in C^{1,2}(Q): (\varphi-u)(t,x)=\min_Q(\varphi-u)
 \big\},
 \label{AsubsolMarkov}
 \\
 \overline{A}v(t,x)
 &:=&\big\{\varphi\in C^{1,2}(Q): (\varphi-v)(t,x)=\max_Q(\varphi-v)
    \big\}.
 \label{AsupersolMarkov}
 \ee

\begin{Definition}
{\rm (i)} $u\in\mbox{\rm USC}(Q)$ is a viscosity subsolution of equation \eqref{PDE} if: 
 \b*
 \big\{-\partial_t\varphi
 -g(.,u,D\varphi,D^2\varphi)\big\}(t,x)
 \le 0
 &\mbox{for all}&
 (t,x)\in Q,~\varphi\in\underline{A}u(t,x).
 \e*
{\rm (ii)} $v\in\mbox{\rm LSC}(Q)$ is a viscosity supersolution of equation \eqref{PDE} if: 
 \b*
 \big\{-\partial_t\varphi
 -g(.,u,D\varphi,D^2\varphi)\big\}(t,x)
 \ge 0
 &\mbox{for all}&
 (t,x)\in Q,~\varphi\in\overline{A}u(t,x).
 \e*
{\rm (iii)} A viscosity solution of \eqref{PDE} is a viscosity subsolution and supersolution of \eqref{PDE}.
 \end{Definition}

From the last definition, it is clear that one may add a constant to the test function $\varphi$ so that the minimum and the maximum values in \eqref{AsubsolMarkov}-\eqref{AsupersolMarkov} are zero. Then, the pictorial representation of a test function $\varphi\in\underline{A}u(t,x)$ is a smooth function tangent from above to $u$ with contact point at $(t,x)$. The symmetric pictorial representation holds for a test function $\varphi\in\overline{A}v(t,x)$. Notice that $\underline{A}v(t,x)$ may be empty, and in this case the subsolution property at $(t,x)$ holds trivially.  

We also observe that we may replace the minimum and maximum in \eqref{AsubsolMarkov}-\eqref{AsupersolMarkov} by the corresponding local notions. Moreover, by the continuity of the nonlinearity $g$, we may also assume the minimum (reps. maximum) or local minimum (resp. local maximum) to be strict, and we may restrict attention to $C^\infty(Q)$ test functions.

The following consistency property is an easy consequence of the ellipticity condition on $g$. We state it only for subsolution, but the result can be similarly stated for supersolutions.

\begin{Proposition}
Assume $g(t,x,y,z,\gamma)$ is non-decreasing in $\gamma$. Then, for a function $u\in C^{1,2}(Q)$, we have
 \b*
 u~\mbox{is a classical subsolution of \eqref{PDE}}
 &\mbox{iff}&
 u~\mbox{is a viscosity subsolution of \eqref{PDE}}.
 \e*
\end{Proposition}

\subsection{The heat equation example}
\label{subsect:heatequation}

In this subsection, we consider the equation
 \be\label{heat-Markov}
 -Lu(t,x):=
 -\partial_tu(t,x)-b(t,x)Du(t,x)-\frac12\sigma^2(t,x):D^2u(t,x) =0,
 &(t,x)\in Q.&
 \ee
where the coefficients $b:Q\longrightarrow\R^d$ and $\sigma:Q\longrightarrow\S_d$ are continuous and Lipschitz-continuous in $x$ uniformly in $t$. The purpose from studying this simple example is to gain some intuition in view of our extension to the path-dependent case. 

Under the above conditions on $b$ and $\sigma$, we may consider the unique strong solution $\{X_t,t\in[0,T]\}$ of the stochastic differential equation
 \be\label{SDE}
 X_t
 &=&
 X_0 +\int_0^t b(s,X_s)ds + \int_0^t \sigma(s,X_s)dB_s,
 ~~\P_0-a.s.
 \ee
for some given initial data $X_0$. Then, given a boundary condition $u(T,.)=\psi$ for some $\psi:\R^d\longrightarrow\R$, the natural solution of \eqref{heat-Markov} is given by:
 \b*
 u^0(t,x)
 :=
 \E^{\P_0}\big[\psi(X_T)|X_t=x\big],
 &(t,x)\in Q.
 \e*
In the remaining of this section, we verify that $u^0$ is a viscosity solution of the heat equation \eqref{heat-Markov}, and we make crucial observations which open the door for enlarging the set of test functions.
\begin{enumerate}
\item[{\bf (a)}] {\it\bf Tower property}\quad The first step is to use the Markov feature of the process $X$ in conjunction with the tower property to deduce that
 \be\label{Tower-Markov}
 u(t,x)
 =
 \E^{\P_0}\big[u(\tau,X_\tau)|X_t=x\big]
 &\mbox{for all stopping time $\tau$ with values in}&
 [t,T].
 \ee
We shall use this identity with stopping times $\tau=\tau_h:=(t+h)\wedge\inf\{s>t:|X_s-x|\ge 1\}$. For the next development, notice that $\tau_h>t$, a.s., $(s,X_s)$ is bounded on $[t,\tau_h]$, and $\tau_h\longrightarrow t$, a.s. when $h\searrow 0$.

Also, we avoid to discuss the regularity issues of the function $u^0$. For instance, if $\psi$ is Lipschitz-continuous, then $u^0$ is immediately seen to be Lipschitz-continuous with respect to the $x-$variable, uniformly in $t$, and we verify that $u^0$ is $\frac12-$H\"older-continuous with respect to the $t-$ variable, uniformly in $x$, by using the identity \eqref{Tower-Markov}. 
\item[{\bf (b)}] {\it\bf  $u^0$ is a viscosity subsolution}\quad Let $(t,x)\in Q$ and $\varphi\in\underline{A}u(t,x)$ be given. We denote by $\{X^{t,x}_s,s\in[t,T]\}$ the solution of \eqref{SDE} started from $X^{t,x}_t=x$. By definition, we have 
\begin{equation}\label{existencesharpsubsol}
 (\varphi-u)(t,x)\le (\varphi-u)
 ~\mbox{on $Q$, and then}~
 (\varphi-u)(t,x)\le \E^{\P_0}\big[(\varphi-u)(\tau_h,X^{t,x}_{\tau_h})\big],
\end{equation}
for all $h>0$. From the last inequality in mean, together with the identity \eqref{Tower-Markov}, we get
 \b*
 \varphi(t,x) 
 &\le& 
 \E^{\P_0}\big[\varphi(\tau_h,X^{t,x}_{\tau_h})\big].
 \e*
Since the test function $\varphi$ is smooth, it follows from It\^o's formula that
 \b*
 -\E^{\P_0}\Big[\int_t^{\tau_h} L\varphi(r,X^{t,x}_r)dr\Big]
 &\le&
 0.
 \e*
Dividing by $h$ and sending $h\searrow 0$, we deduce from the mean value theorem together with the dominated convergence theorem that
 \b*
 -L\varphi(t,x) &\le& 0.
 \e*
\item[{\bf (c)}] {\it\bf  $u^0$ is a viscosity supersolution}\quad For $(t,x)\in Q$ and $\varphi\in\overline{A}u(t,x)$, notice that we have the analogue of \eqref{existencesharpsubsol}:
\begin{equation}\label{existencesharpsupersol}
 (\varphi-u)(t,x)\ge (\varphi-u)
 ~\mbox{on $Q$, and then}~
 (\varphi-u)(t,x)\ge \E^{\P_0}\big[(\varphi-u)(\tau_h,X^{t,x}_{\tau_h})\big].
 \end{equation}
Following the same line of argument as in (b), it follows that $L\varphi(t,x)\ge 0$, as required.
\end{enumerate}

\no {\bf Crucial observation}\quad Notice that only the right-hand sides of \eqref{existencesharpsubsol} and \eqref{existencesharpsupersol} have been useful to prove that $u^0$ is a subsolution and supersolution, respectively, of the heat equation \eqref{heat-Markov}. The right-hand sides of \eqref{existencesharpsubsol} and \eqref{existencesharpsupersol} express that the test function $\varphi$ is tangent to $u$ in mean, locally along the trajectory of the underlying process $(s,X^{t,x}_s)$. Of course, the set of smooth functions which are tangent from above (reps. from below) in mean is larger than $\underline{A}u$ (resp. $\overline{A}u$). Consequently, we may have used an alternative definition of viscosity solution with a richer family of test functions (defined by the right-hand sides of \eqref{existencesharpsubsol} and \eqref{existencesharpsupersol}), and still get the same existence result. The benefit from such a stronger definition may be that the uniqueness theory can be simplified by suitable use of the additional test functions. 

\begin{Remark}\label{rem:smoothprocesses}{\rm 
The $C^{1,2}$ smoothness of the test function $\varphi$ is only needed in order to apply It\^o's formula
 \b*
 \varphi\big(\tau_h,X^{t,x}_{\tau_h}\big)-\varphi(t,x)
 &=&
 \int_t^{\tau_h}L\varphi(r,X^{t,x}_r)dr
 +\int_t^{\tau_h}D\varphi(r,X^{t,x}_r)\sigma(r,X^{t,x}_r)dB_r,
 ~~\P_0-\mbox{a.s.}
 \e*
Motivated by this observation, we shall take It\^o's formula as a starting point for the definition of smooth processes in the path-dependent case.
}
\end{Remark}

\subsection{Existence for HJB equations}
\label{subsect:HJBequation}

In this subsection, we show that the crucial observation from the previous subsection holds in the context of the fully nonlinear Hamilton-Jacobi-Bellman (HJB) equation:
 \be\label{HJB-Markov}
 -\partial_tu
 -\sup_{k\in K}\Big\{ b(.,k)Du+\frac12\sigma^2(.,k):D^2u\Big\}
 =0,
 &(t,x)\in Q.&
 \ee
Here, for simplicity, we consider the case of a bounded set of controls $K$. The controlled coefficients $b:Q\times K\longrightarrow\R^d$ and $\sigma:Q\times K\longrightarrow\S_d$ are continuous in $t$, Lipschitz-continuous in $x$ uniformly in $(t,\kappa)$. The controls set is denoted by $\Kc$, and consists of all $\F-$progressively measurable process with values in $K$.
For all control process $\kappa\in\Kc$, we introduce the controlled process $X^{\kappa}$ as the unique strong solution of the SDE
 \b*
 X^\kappa_t=X_0+\int_0^t b(s,X^\kappa_s,\kappa_s)ds
                           +\int_0^t \sigma(s,X^\kappa_s,\kappa_s)dB_s,
 &\P_0-\mbox{a.s.}&
 \e*
and we denote by $X^{\kappa,t,x}$ the solution corresponding to the initial data $X^{\kappa,t,x}_t=x$. The Dynkin operator associated to $X^\kappa$ is denoted:
 \b*
 L^k
 &:=&
 \partial_t+b(.,k)D+\frac12\sigma^2(.,k):D^2.
 \e*
Given a boundary condition $u(T,.)=\psi$ for some $\psi:\R^d\longrightarrow\R$, the natural solution of \eqref{HJB-Markov} is given by:
 \b*
 u^1(t,x)
 &=&
 \sup_{\kappa\in\Kc}
 \E^{\P_0}\big[ \psi(X^{\kappa,t,x}_T)\big],
 ~~(t,x)\in Q.
 \e*
In the remaining of this section, we verify that $u^1$ is a viscosity supersolution of the HJB equation \eqref{HJB-Markov}, and we focus on the crucial observation that only the tangency condition in mean is used for this purpose. The subsolution property can be obtained by similar standard arguments, and the reader can verify that only tangency in mean is needed, again.
\begin{enumerate}
\item[{\bf (a)}] {\it\bf Dynamic programming principle}\quad In the present nonlinear case, the substitute for the tower property identity \eqref{Tower-Markov} is the following dynamic programming identity:
 \begin{equation}\label{DPP-Markov}
 u(t,x)
 =
 \sup_{\kappa\in\Kc}
 \E^{\P_0}\big[u(\tau^\kappa,X^{\kappa,t,x}_{\tau^\kappa})\big]
 ~\mbox{for all stopping times $\tau^\kappa$ with values in}~
 [t,T].
 \end{equation}
This identity will be used with stopping times $\tau^\kappa=\tau^\kappa_h:=(t+h)\wedge\inf\{s>t:|X^{\kappa,t,x}_s-x|\ge 1\}$. For the next development, notice that $\tau_h>t$, a.s., $(s,X^\kappa_s)$ is bounded on $[t,\tau^\kappa_h]$, and $\tau^\kappa_h\longrightarrow t$, a.s. when $h\searrow 0$.

The proof of \eqref{DPP-Markov} is a difficult task relying on involved measurable selections techniques, see \cite{STZ} for the regular case (which does not require measurable selection arguments), \cite{ET1, ET2} for the general irregular case, and \cite{BT} for a weak dynamic programming principle which is sufficient for the task of deriving the viscosity property, while by-passing the measurable selection arguments.

We also avoid here to discuss the regularity issues of the function $u^1$. For instance, if $\psi$ is Lipschitz-continuous, then $u^1$ is immediately seen to be Lipschitz-continuous with respect to the $x-$variable, uniformly in $t$, and we verify that $u^1$ is $\frac12-$H\"older-continuous with respect to the $t-$ variable, uniformly in $x$, by using the identity \eqref{DPP-Markov}. 

\item[{\bf (b)}] {\it\bf  $u^1$ is a viscosity supersolution}\quad Let $(t,x)\in Q$ and $\varphi\in\overline{A}u(t,x)$ be given. Fix an arbitrary control process $\kappa\in\Kc$. For the purpose of the present argument, we may take this control porcess to be constant $\kappa_s=k$ for all $s\in[t,T]$. By definition, we have 
\begin{equation}\label{existencesharpsupersol-HJB}
 (\varphi-u)(t,x)\ge (\varphi-u)
 ~\mbox{on $Q$, and then}~
 (\varphi-u)(t,x)\ge \E^{\P_0}\big[(\varphi-u)(\tau^\kappa_h,X^{\kappa,t,x}_{\tau^\kappa_h})\big],
 \end{equation}
for all $h>0$. From the last inequality in mean, together with the identity \eqref{DPP-Markov}, we get
 \b*
 \varphi(t,x) 
 &\ge& 
 \E^{\P_0}\big[\varphi(\tau_h,X^{\kappa,t,x}_{\tau_h})\big].
 \e*
Since the test function $\varphi$ is smooth, it follows from It\^o's formula that
 \b*
 -\E^{\P_0}\Big[\int_t^{\tau_h} L\varphi(r,X^{\kappa,t,x}_r)\Big]
 &\ge&
 0.
 \e*
Dividing by $h$ and sending $h\searrow 0$, we deduce from the mean value theorem together with the dominated convergence theorem that
 \b*
 -L^k\varphi(t,x) &\ge& 0.
 \e*
By the arbitrariness of $k\in K$, this proves the required supersolution property.
\end{enumerate}

\no {\bf Crucial observation}\quad Here again, only the right-hand side of \eqref{existencesharpsupersol-HJB} has been useful to prove that $u^1$ is a supersolution of the HJB equation \eqref{HJB-Markov}. The right-hand side of \eqref{existencesharpsupersol-HJB} expresses that the test function $\varphi$ is tangent to $u$ in mean, locally along the trajectory of the underlying process $(s,X^{\kappa,t,x}_s)$, for all possible control process $\kappa\in\Kc$. 
The latter is a new feature which appears in the present nonlinear case: while the linear case involves the tangency condition under the expectation operator $\E^{\P_0}$, the present nonlinear case requires the use of a sub linear expectation defined by an additional maximization with respect to all possible choices of control process $\kappa\in\Kc$. 

This additional feature however does not alter the observation that the set of smooth functions which are tangent from below in (sublinear) mean is larger than $\overline{A}u$. Consequently, we may have used an alternative definition of viscosity solution with a richer family of test functions (defined by the right-hand side of \eqref{existencesharpsupersol-HJB}), and still get the same existence result. Similar to the case of the linear heat equation, the benefit from such a stronger definition may be that the uniqueness theory can be simplified by suitable use of the additional test functions. 

\subsection{Comparison of viscosity solutions}

The uniqueness result of viscosity solution of second order fully nonlinear elliptic PDEs is usually obtained as a consequence of the comparison result, which corresponds to the maximum principle. 

\begin{Definition}
We say that the equation \eqref{PDE} satisfies comparison of bounded solutions if for all bounded viscosity subsolution $u$, and bounded viscosity supersolution $v$, we have 
 \b*
 (u-v)(T,.) \le 0~~\mbox{on}~~\R^d
 &\mbox{implies}&
 u-v\le 0~~\mbox{ on cl}(Q).
 \e*
\end{Definition}

Comparison results for viscosity solution are available for a wide class of equations. The most accessible results are for the case of first order equations where the beautiful trick of doubling variables is remarkably efficient. 

For second order equations, comparison results are more difficult and require to introduce a convenient regularization, typically by inf-convolution. The most general approach which covers possibly degenerate equations relies crucially on the Crandall-Ishii Lemma which provides the substitute of first and second order conditions at a local maximum point when the objective function is only upper semicontinuous. 

In the context of uniformly elliptic equations, the argument of Caffarelli \& Cabre \cite{CC} avoids the technique of doubling variables, but still relies crucially on the inf-convolution regularization. We refer to Wang \cite{Wang} for the extension to the uniformly parabolic case which requires a more involved regularization technique. 

All available comparison results for second order elliptic and parabolic equations use the restriction of test functions to paraboloids. This leads to the notion of superjets and subjets. For notations consistency, we continue our discussion with the parabolic case.

For $q\in\R$, $p\in\R^d$, and $\gamma\in\S_d$, we introduce the paraboloid function:
 \b*
 \phi^{q,\beta,\gamma}(t,x)
 &:=&
 q t+p\cdot x
 +\frac12 \gamma x\cdot x,
 ~~ (t,x)\in Q.
 \e* 
For $u\in\mbox{USC}(Q)$, let $(t_0,x_0)\in Q$, $\varphi\in\underline{A}u(t_0,x_0)$, define $q:=\partial_t\varphi(t_0,x_0)$, $p:=D\varphi(t_0,x_0)$, and $\gamma:=D^2\varphi(t_0,x_0)$. Then, it follows from a Taylor expansion that:
 \b*
 u(t,x)
 \le
 u(t_0,x_0)+\phi^{q,p,\gamma}(t-t_0,x-x_0)
 +\circ\big(|t-t_0|+|x-x_0|^2\big).
 \e*
Motivated by this observation, we introduce the {\it superjet} $J^+u(t_0,x_0)$ by
 \be\label{superjet}
 J^+u(t_0,x_0)
 &:=&
 \big\{(q,p,\gamma)\in\R\times\R^d\times\S_d:~\mbox{for all}~(t,x)\in Q\\
  &&\hspace{5mm}
         u(t,x)\le u(t_0,x_0)+\phi^{q,p,\gamma}(t-t_0,x-x_0)
                     +\circ\big(|t-t_0|+|x-x_0|^2\big)
 \big\}.
 \nonumber
 \ee
Then, it can be proved that a function $u\in \mbox{USC}(Q)$ is a viscosity subsolution of the equation \eqref{PDE} if and only if 
 \b*
 F(t,x,u(t,x),p,\gamma) \le 0
 &\mbox{for all}&
 (q,p,\gamma)\in J^+u(t,x).
 \e*
The nontrivial implication of the previous statement requires to construct, for every $(q,p,A)\in  J^+u(t,x)$, a smooth test function $\varphi$ such that the difference  $(\varphi-u)$ has a local minimum at $(t,x)$. 

Similarly, we define the  {\it subjet} $J^-v(t_0,x_0)$ of a function $v\in\mbox{LSC}(Q)$ at the point $(t_0,x_0)\in Q$ by
 \be\label{subjet}
 J^-v(t_0,x_0)
 &:=&
 \Big\{(q,p,\gamma)\in\R^d\times\S_d: ~\mbox{for all}~(t,x)\in Q\\
 &&\hspace{5mm}
         v(x)\ge v(t_0,x_0)+\phi^{q,p,\gamma}(t-t_0,x-x_0)
                      +\circ\big(|t-t_0|+|x-x_0|^2\big)
 \Big\},
 \nonumber
 \ee
and $v\in \mbox{LSC}(Q)$ is a viscosity supersolution of the equation \eqref{PDE} if and only if 
 \b*
 F(t,x,v(t,x),p,\gamma) \ge 0
 &\mbox{for all}&
 (q,p,\gamma)\in J^-u(t,x).
 \e*
By continuity considerations, we can even enlarge the semijets $J^\pm$ to the following closure
\b*
 \tilde{J}^\pm w(t,x)
 &:=&
 \Big\{(q,p,\gamma)\in\R^d\times\S_d: 
        (t_n,x_n, w(t_n,x_n),q_n,p_n,\gamma_n)\longrightarrow (t,x,w(t,x),q,p,\gamma) 
 \\ & & \hspace{20mm}
        \mbox{for some sequence}
        ~(t_n,x_n,q_n,p_n,\gamma_n)_n\subset\mbox{\mbox{Graph}}(J^\pm w)
 \Big\},
 \e*
where $(t_n,x_n,q_n,p_n,\gamma_n)\in\mbox{\rm Graph}(J^\pm w)$ means that $(q_n,p_n,\gamma_n)\in J^\pm w(t_n,x_n)$. The following result is obvious, and provides an equivalent definition of viscosity solutions.

\begin{Proposition}\label{propsemijet}
Let $u\in\mbox{\rm USC}(Q)$, and $v\in\mbox{\rm LSC}(Q)$.
\\
{\rm (i)}\quad Assume that $g$ is lower-semicontinuous. Then, $u$ is a viscosity subsolution of \eqref{PDE} iff:
 \b*
 -q-g(t,x,u(t,x),p,\gamma) \le 0
 &\mbox{for all}&
 (t,x)\in Q~\mbox{and}~(q,p,\gamma)\in\tilde{J}^+u(t,x).
 \e*
{\rm (ii)}\quad Assume that $g$ is upper-semicontinuous. Then, $v$ is a viscosity supersolution of \eqref{PDE} iff:
 \b*
 -q-g(t,x,v(t,x),p,\gamma) \ge 0
 &\mbox{for all}&
 (t,x)\in Q~\mbox{and}~(q,p,\gamma)\in\tilde{J}^-v(t,x).
 \e*
\end{Proposition}

\subsection{Stability of viscosity solutions}

We conclude this section by reviewing the stability property of viscosity solutions. The following result is expressed in the context of our parabolic fully-nonlinear equation. However, the reader can see from its proof that it holds for general degenerate second order elliptic equations. We consider a family of equations parameterized by $\eps>0$:
 \bea\label{geps}
 -\partial_tu-g^\eps(x,u,Du,D^2u)=0
 &\mbox{on}&
 Q,
 \eea
and we consider the convergence problem of a corresponding family of subsolutions $(u^\eps)_{\eps>0}$. The main ingredient for the stability result is the notion of relaxed semi limits introduced by Barles and Perthame \cite{BarlesPerthame}:
 \beaa
 \overline{u}(t,x):=\limsup_{(\eps,t',x')\to(0,t,x)} u^\eps(t',x')
 &\mbox{and}&
 \overline{g}(\zeta):=\limsup_{(\eps,\zeta')\to(0,\zeta)} g^\eps(\zeta'),
 \eeaa
where $\zeta=(t,x,y,z,\gamma)$. Notice that the semilimits here are taken both in the variables and the small parameter $\eps$, and are finite whenever the functions of interest are locally bounded in the corresponding variables and the small parameter $\eps$.

\begin{Theorem}\label{thmstability}
Let $u^\eps\in\mbox{USC}(Q)$ be a viscosity subsolution of \eqref{geps} for all $\eps>0$. Suppose that the maps $(\eps,x)\longmapsto u_\eps(x)$ and $(\eps,\zeta)\longmapsto g^\eps(\zeta)$
are locally bounded. Then, $\overline{u}\in\mbox{USC}(Q)$ is a viscosity subsolution of
the equation
 \bea\label{gbar}
 -\partial_t\overline{u}-\overline{g}(x,\overline{u},D\overline{u},D^2\overline{u}) = 0
 &\mbox{on}&
 Q,
 \eea
A similar statement holds for supersolutions.
\end{Theorem}

\proof The fact that $\overline{u}$ is upper semicontinuous is an easy exercise. Let $\varphi\in\underline{A}\overline{u}(t,x)$. Without loss of generality, we may assume that the test function $\varphi$ is strictly tangent from above to $\overline{u}$ at the point $(t,x)$, i.e.
 \bea
 (\varphi-\overline{u})(t,x) < (\varphi-\overline{u})(t',x')
 &\mbox{for all}&
 (t',x')\in Q,~(t',x')\neq (t,x).
 \eea
By definition of $\overline{u}$, there is a sequence $(\eps_n,x_n)$ $\in$ $(0,1]\times\R^d$ such that
 \beaa
 (\eps_n,t_n,x_n) \;\longrightarrow\; (0,t,x)
 &\mbox{and}&
 u^{\eps_n}(t_n,x_n) \;\longrightarrow\; \overline{u}(t,x).
 \eeaa
Let $O$ be an open subset of $Q$ containing $(t,x)$ and $(t_n,x_n)_n$. Let $(\bar t_n,\bar x_n)$ be a
minimizer of $\varphi-u^{\eps_n}$ on cl$(O)$. We claim that
 \bea\label{claimstability}
 (\bar t_n,\bar x_n) \longrightarrow (t,x)
 ~\mbox{and}~
 u^{\eps_n}(\bar t_n,\bar x_n) \longrightarrow \overline{u}(t,x)
 &\mbox{as}&
 n\to\infty.
 \eea
Before verifying this, let us complete the proof. We first deduce that $(\bar t_n,\bar x_n)$ is an interior minimizer of the difference $(\varphi-u^{\eps_n})$. Then, it follows from the viscosity subsolution property of $u^{\eps_n}$ that:
 \beaa
 0
 &\ge&
 \big\{-\partial_t\varphi
         -g^{\eps_n}\left(.,u^{\eps_n},D\varphi,D^2\varphi\right)\big\}(\bar t_n,\bar x_n).
 \eeaa
Then, taking limits on both sides, we see that
 $$
 0
 \ge
 -\partial_t\varphi(t,x)
 -\limsup_{n\to\infty} 
   g^{\eps_n}\left(.,u^{\eps_n},D\varphi,D^2\varphi\right)(\bar t_n,\bar x_n)
 \ge
 \big\{-\partial_t\varphi
 -\overline{g}\left(.,\overline{u},D\varphi,D^2\varphi\right)\big\}(t,x),
 $$
by \eqref{claimstability} and the definition of $\overline{g}$.

It remains to prove Claim \reff{claimstability}. Recall that $(\bar t_n,\bar x_n)_n$ is valued in the compact set cl$(O)$. Then, there is a subsequence, still named $(\bar t_n,\bar x_n)_n$, converging to some $(\bar t,\bar x)\in\mbox{cl}(O)$. We now prove that $(\bar t,\bar x)=(t,x)$ and obtain the second claim in \reff{claimstability} as a by-product. By the fact that $(\bar t_n,\bar x_n)$ is a minimizer of $(\varphi-u^{\eps_n})$ on cl$(O)$, together with the definition of $\overline{u}$, we see that
 \beaa
 0 \;=\; (\varphi-\overline{u})(t,x)
   &=& \lim_{n\to\infty}\; \left(\varphi-u^{\eps_n}\right)(t_n,x_n) 
   \\
   &\ge&
   \limsup_{n\to\infty}\; \left(\varphi-u^{\eps_n}\right)(\bar t_n,\bar x_n) 
   \\
   &\ge& 
   \liminf_{n\to\infty}\; \left(\varphi-u^{\eps_n}\right)(\bar t_n,\bar x_n) \\
   &\ge& (\varphi-\overline{u})(\bar t,\bar x).
 \eeaa
We now obtain \reff{claimstability} from the fact that $(t,x)$ is a strict minimizer of the difference $(\varphi-\overline{u})$.
 \ep

\section{Viscosity solution of path-dependent PDEs}
\label{sect-visco}
\setcounter{equation}{0}

We now turn to the main purpose of this paper, namely the theory of viscosity solutions for path-dependent PDEs  \reff{PPDE}:
 \beaa
 -\partial_tu(t,\omega)
 - G\big(t,\omega,u(t,\omega),\partial_\omega u(t,\omega), \partial^2_{\omega\omega}u(t,\omega)
      \big)
 =0,
 &t<T,&
 \omega\in\Omega.
 \eeaa
  where the generator $G:[0,T]\times\Omega\times\R\times\R^d\times\S_d\longrightarrow \R$ is a continuous  map satisfying the ellipticity condition  \reff{elliptic}. We recall that  $\Omega:=\{\omega\in C^0([0,T],\R^d):\omega_0=0\}$ is the underlying canonical space,  $B_t(\omega):=\omega_t$, $t\in[0,T]$,  is the canonical process, $\dbP_0$ is the Wiener measure, $\F:=\{\Fc_t, t\in[0,T]\}$ with $\Fc_t=\sigma(B_s,s\le t)$   is the natural filtration equipped with the pseudo-distance $d$ defined in \reff{Dupire-distance}. Moreover, denote 
  \beaa
  \Theta := [0, T)\times \O,\quad \overline \Theta := [0, T]\times \O,
  \eeaa
 and   $C^0(\ol \Theta)$ is the set of continuous  processes on $\ol\Theta$. We note that  any  $u \in C^0(\ol \Theta)$ is $\F-$progressively measurable, namely  $u(t,\omega)=u\big(t,(\omega_s)_{s\le t}\big)$. 


 
 \subsection{Differentiability}
 
Before introducing the notion of viscosity solutions for this path-dependent PDE, we first need to specify the meaning of the time derivatives $\partial_tu(t,\omega)$ and the spatial derivatives $\partial_\omega u(t,\omega)$ and $\partial^2_{\omega\omega} u(t,\omega)$. Once these derivatives are clearly defined, we would have, on one hand a natural definition of classical solutions for the path-dependent PDE, and on the other hand a natural set of smooth functions to serve as test functions for our notion of viscosity solutions.



These path derivatives were first introduced by Dupire \cite{Dupire}. In particular,  
\cite{Dupire} defines the vertical derivatives (our spatial derivatives) by bumping the path at time $t$.
While such a definition is natural in the larger space of discontinuous paths, our paths space $\Omega$ would require an extension of the map $u$ to the set of discontinuous paths. We refer to Cont \& Fourni\'e \cite{CF} for this approach, where it is proved in particular that such a vertical derivative, if exists, does not depend on the choice of the extension of $u$ to the set of discontinuous paths. Motivated by Remark \ref{rem:smoothprocesses}, we adopt the following notion of smoothness.

\begin{Definition}[Smooth processes]\label{def:C12}
Let $\Pc$ be a set of probability measures on $\Omega$ with $B$ a $\P-$semimartingale for all $\P\in\Pc$. We say that $u\in C^{1,2}_\Pc(\Theta)$ if
 $u\in C^0(\ol\Theta)$ and there exist processes $\alpha, Z,\Gamma\in C^0(\Theta)$ valued in $\R$, $\R^d$ and $\S_d$, respectively, such that:
\b*
 du_t
 &=&
 \alpha_tdt + \frac12 \Gamma_t:d\langle B\rangle_t + Z_t dB_t,
 ~~\P-\mbox{a.s. for all}~\P\in\Pc.
 \e*
 The processes $\alpha$, $Z$ and $\Gamma$ are called the time derivative, spacial gradient and spatial Hessian, respectively, and we denote $\pa_t u := \alpha$, $\partial_\omega u_t:=Z_t$, $\partial^2_{\omega\omega} u_t:=\Gamma_t$.

\end{Definition}

We observe that any $C^{1,2}$ process in the Dupire sense is in $C^{1,2}_\Pc(\Theta)$. This is an immediate consequence of the functional It\^o formula proved in Dupire \cite{Dupire} and \cite{CF}. In particular, our notion of smooth processes is weaker than the corresponding one in \cite{Dupire}. We also note that, when $\cP$ is rich enough, our path derivatives are unique.

\begin{Remark}{\rm
The previous definition does not require that $\partial^2_{\omega\omega}u_t$ be the derivative (in some sense) of $\partial_\omega u_t$. This is very well illustrated by the following example communicated by Mete Soner. Let $d=2$, and $u_t:=\int_0^tB^1_sdB^2_s$ which is defined pathwise due to the results of Karandikar \cite{Karandikar}. 
\\
$\bullet$ Clearly $\partial_tu=0$. Since $du_t=B^1_t dB^2_t$, under any semimartingale measure, we also deduce that $\partial_\omega u_t=(0,B^1_t)^{\rm T}$, and $\partial^2_{\omega\omega} u_t=0_{\Sc_2}$. Hence $u\in C^{1,2}_\Pc(\Theta)$ for any subset $\Pc$ of the collection of all semimartingale measures for $B$. 
\\
$\bullet$ Let $\partial_\omega^{\rm\footnotesize D}u_t$ and $\partial_{\omega\omega}^{{\rm\footnotesize D}^2}u_t$ denote the vertical first and second derivatives in the Dupire sense. Direct calculation reveals that $\partial_\omega^{\rm\footnotesize D}u_t=(0,B^1_t)^{\rm T}=\partial_\omega u_t$. However,
 $$
 \partial_{\omega\omega}^{{\rm\footnotesize D}^2}u_t
 =
 \left(\begin{array}{cc}
        0 & 0
        \\
        1 & 0
        \end{array}
 \right), 
 $$
which is not symmetric~!
\\
$\bullet$ However, we need to point out that in this example $u$ does not belong to $C^0(\overline\Theta)$.
\\
$\bullet$ We complement this example by the following observation from a private communication with Bruno Dupire. By considering the Dupire vertical derivative as originally defined on the set of discontinuous paths, we see by direct calculation that $\partial_{\omega\omega}^{{\rm\footnotesize D}^2}u_t=0_{\Sc_2}=\partial^2_{\omega\omega}u_t$. 
}
\end{Remark}

\begin{Definition}[Classical solution]
Let $\Pc$ be a set of probability measures on $\Omega$ with $B$ a $\P-$semimartingale for all $\P\in\Pc$.
\\
{\rm (i)} $u\in C^{1,2}_\Pc(\Theta)$ is a $\Pc-$classical subsolution of the path-dependent PDE \eqref{PPDE} if 
 \b*
 -\partial_tu-G\big(.,u,\partial_\omega u,\partial^2_{\omega\omega}u\big)
 \le
 0
 &\mbox{on}&
 \Theta.
 \e*
{\rm (ii)} $v\in C^{1,2}_\Pc(\Theta)$ is a $\Pc-$classical supersolution of the path-dependent PDE \eqref{PPDE} if 
 \b*
 -\partial_tv-G\big(.,v,\partial_\omega v,\partial^2_{\omega\omega}v\big)
 \ge
 0
 &\mbox{on}&
 \Theta.
 \e*
{\rm (iii)} A $\Pc-$classical solution of \eqref{PPDE} is both classical subsolution and supersolution.
\end{Definition}

\begin{Example}\label{P-Heatequation}
Let $u(t,\omega):=\E^{\dbP_0}[\xi|\Fc_t]$ for some $\xi\in\L^1(\P_0,\Fc_T)$, and assume $u\in C^{1,2}_{\P_0}(\Theta)$. By definition, this implies that
 \b*
 du_t
 &=&
 \big(\partial_tu_t+\frac12\partial^2_{\omega\omega}u_t\big) dt
 +\partial_\omega u_t dB_t,
 ~~\P_0-\mbox{a.s.}
 \e*
Since the process $u$ is a martingale, it follows that:
 \b*
 \partial_tu_t+\frac12\partial^2_{\omega\omega}u_t
 =0,
 &(t,\omega)\in\Theta.&
 \e*
In other words, $u$ is a $\P_0-$ classical solution of the path-dependent heat equation.
\end{Example}

\begin{Example}
For $\xi\in\L^2(\P_0,\Fc_T)$, consider the backward stochastic differential equation:
 \b*
 du_t
 =
 -F_t(\omega,u_t,Z_t)dt + Z_t dB_t,
 &u_T=\xi,&
 \P_0-\mbox{a.s.}
 \e*
where $F:[0,T]\times\Omega\times\R\times\R^d\longrightarrow\R$ is continuous, uniformly Lipschitz in $(y,z)$, with $F(0,0)$ a square integrable process. Assume $u\in C^{1,2}_{\P_0}(\Theta)$. By definition, this implies that
 \b*
 du_t
 =
 \big(\partial_tu_t+\frac12\partial^2_{\omega\omega}u_t\big) dt
 +\partial_\omega u_t dB_t
 =
 -F_t(\omega,u_t,Z_t)dt + Z_t dB_t,
 ~~\P_0-\mbox{a.s.}
 \e*
Identifying the martingale terms, we see that $\partial_\omega u_t = Z_t$. Next, identifying the drift term, it follows that $u$ is a $\P_0-$classical solution of the path-dependent semilinear PDE:
 \b*
 -\partial_tu_t-\frac12\partial^2_{\omega\omega}u_t
 -F_t(\omega,u_t,\partial_\omega u_t)
 =0,
 &(t,\omega)\in\Theta.&
 \e*
\end{Example}

\begin{Remark}{\rm
(i) In the Markovian case, strong regularity results are induced by the ellipticity of the underlying diffusion coefficient. The simplest example is when the diffusion is the identity matrix. Let $u(t,x):=\E^{\P_0}[h(B_T)|B_t=x]$. Then $u\in C^\infty([0,T)\times\R^d)$.
\\
(ii) The path-dependency induces specific non-smoothness as outlined by the following example. Let $u(t,\omega):=\E^{\P_0}[B_{\frac{T}{2}}|\Fc_t]=\omega_{t\wedge\frac{T}{2}}$ for all $t\in[0,T]$. Clearly, $\partial_tu_t=0$, and $du_t=\1_{t\le \frac{T}{2}}dB_t$ implying that $\partial_\omega u_t$ is not continuous. Hence $u\not\in C^{1,2}(\Theta)$.
}
\end{Remark}

\subsection{Viscosity solutions of path-dependent PDEs}

\subsubsection{Notations}
First recall our canonical setting $(\O, B, \dbF, \dbP_0)$.
We denote by $\Tc$ the set of all $\F$-stopping times, $\Tc^+\subset\Tc$ the collection of all strictly positive stopping times, and $\Tc^t\subset\Tc$ the subset of the $\F$-stopping times larger than $t$. 

For $\o,\o'\in\O$ and $t\in [0,T]$, we define
 $$
 (\o\otimes_t\o')_s 
 := 
 \o_s\1_{\{s< t\}}
 +(\o_t+\o'_{s-t})\1_{\{s\ge t\}}.
 $$ 
Let $\xi:\O\rightarrow \R$ be $\Fc_T$-measurable random variable. For any $(t,\o)\in\Theta$, define
 \beaa
 \xi^{t,\o}(\o')
 :=
 \xi\big(\o\otimes_t\o'\big)
 &\mbox{for all}&
 \o'\in\O.
 \eeaa
Clearly, $\xi^{t,\o}$ is $\Fc_{T-t}$-measurable, and thus $\Fc_T$-measurable. Similarly, given a process $X$ defined on $\O$, we denote:
$$X^{t,\o}_s(\o^{'}):=X_{t+s}(\o\otimes_t\o^{'}),\ \text{for }s\in[0,T-t].$$
Clearly, if $X$ is $\F$-adapted, then so is $X^{t,\o}$.

Let $\Pc$ be a family of probability measures on $\O$. We also introduce the sublinear and superlinear expectation operators associated to $\Pc$:
 \b*
 \overline{\Ec}^\Pc
 :=\sup_{\P\in\Pc} \E^\P
 &\mbox{and}&
 \underline{\Ec}^\Pc
 :=\inf_{\P\in\Pc} \E^\P.
 \e*
\subsubsection{Definition of Viscosity Solutions}

We recall that the nonlinearity $G$ satisfies the ellipticity condition in \eqref{elliptic}. We assume in addition that it is $L_0-$Lipschitz with respect to the arguments $(y,z,\gamma)$, uniformly in $(t,\omega)$:
 \be\label{L0-Lip}
 \big|G(t,\omega,y,z,\gamma)-G(t,\omega,y',z',\gamma')\big|
 &\le&
 L_0\big(|y-y'|+|z-z'|+|\gamma-\gamma'|\big)
 \ee
for all $y,y'\in\R$, $z,z'\in\R^d$, $\gamma,\gamma'\in\S_d$, $(t,\omega)\in\Theta$.

Ou Definition \ref{def:C12} of smooth processes involves a family of probability measures that we intentionally did not discuss so far. We now introduce a specific family of semimartingale measures which will be needed for our notion of viscosity solutions.  

\begin{Definition}\label{def:PL} By $\Pc_L$ we denote the collection of all continuous semimartingale measures $\dbP$ on $\O$ whose drift and diffusion characteristics are bounded by $L$ and $\sqrt{2L}$, respectively.  


\end{Definition}

 We refer to  \cite{ETZ0} for properties of this class.
In our subsequent analysis, the family of probability measures $\Pc$ is a subset of $\Pc_L$ for some $L>0$. 

Motivated by the crucial observations of Subsections \ref{subsect:heatequation} and \ref{subsect:HJBequation}, we introduce the sets of test processes:
 \b*
 \underline{\Ac}^\Pc u_t(\omega)
 &:=&
 \Big\{\varphi\in C^{1,2}_\Pc(\Theta):
          (\varphi-u^{t,\o})_0
          =
          \min_{\tau\in\Tc} \underline{\Ec}^\Pc\big[(\varphi-u^{t,\o})_{\tau\wedge \hc}\big]
          ~\mbox{for some}~\hc\in\Hc^\Pc_+
 \Big\},
 \\
 \overline{\Ac}^\Pc v_t(\omega)
 &:=&
 \Big\{\varphi\in C^{1,2}_\Pc(\Theta):
          (\varphi-v^{t,\o})_0
          =
          \max_{\tau\in\Tc} \overline{\Ec}^\Pc\big[(\varphi-v^{t,\o})_{\tau\wedge \hc}\big]
          ~\mbox{for some}~\hc\in\Hc^\Pc_+
 \Big\},
 \e*
where $\Hc^\Pc_+\subset\Tc^+$ satisfies the following properties, for all $\hc,\hc'\in\Hc^\Pc_+$:
\begin{equation}\label{propH}
\begin{array}{l}
\mbox{{\bf H1} (stability by minimization)}~~\hc\wedge\hc'\in\Hc^\Pc_+,
\\
\mbox{{\bf H2} (stability by localization)}~~\hc\wedge \hc_\eps\in\Hc^\Pc_+,~\mbox{where}~ \hc_\eps:=\eps\wedge\inf\big\{t>0:|B|_t\ge\eps\big\}.
\end{array}
\end{equation}
Later, we will call $\hc$ the localizing stopping time (or the localization) of the corresponding test process $\f$.

\begin{Definition}[Viscosity solution of path-dependent PDE]\label{def:visco} Let $u,v\in C^0(\ol\Theta)$.
\\
{\rm (i)} $u$ is a $\Pc-$viscosity subsolution of \eqref{PPDE} if:
 \b*
 \big\{-\partial_t\varphi
         -G\big(.,u,\partial_\omega\varphi,\partial^2_{\omega\omega}\varphi\big)
 \big\}(t,\omega)
 \le
 0
 &\mbox{for all}&
  (t,\omega)\in\Theta,
 ~\varphi\in\underline{\Ac}^\Pc u_t(\omega).
 \e*
{\rm (ii)} $v$ is a $\Pc-$viscosity supersolution of \eqref{PPDE} if:
 \b*
 \big\{-\partial_t\varphi
         -G\big(.,v,\partial_\omega\varphi,\partial^2_{\omega\omega}\varphi\big)
 \big\}(t,\omega)
 \ge
 0
 &\mbox{for all}&
 (t,\omega)\in\Theta,
 ~\varphi\in\overline{\Ac}^\Pc v_t(\omega).
 \e*
{\rm (iii)} A $\Pc-$viscosity solution of \eqref{PPDE} is both a $\Pc-$subsolution and a $\Pc-$supersolution.
\end{Definition}

\begin{Remark}{\rm
in the Markovian case, we may as well use the last definition as an alternative to the standard notion of viscosity solutions. Compared to the standard notion reviewed in Section \ref{sect-standardvisco}, we see that any $\phi\in\underline{A}u(t,x)$ induces a process $\varphi(t,\omega):=\phi(t,\omega_t)$ which obviously lies in $\underline{\Ac}^\Pc u_t(\omega)$. However, even in the Markovian case $u_t(\omega)=u(t,\omega_t)$, a test process in $\underline{\Ac}^\Pc u_t(\omega)$ does not necessarily induce a test function in $\underline{A} u(t,\omega_t)$. 
Thus, our notion of viscosity solution involves more test functions than the standard notion. A viscosity subsolution/supersolution in sense of Definition \ref{def:visco} is restricted by a richer family of test functions. Consequently:
\\
$\bullet$ under our definition, we may hope to take advantage of the richer family of test functions in order to obtain an easier uniqueness proof,
\\
$\bullet$ under our definition, the existence problem is more restricted than under the standard theory of viscosity solutions.
}
\end{Remark}

\begin{Remark}
\label{rem-cA}
{\rm Due to the stability property of the set $\Hc^\Pc_+$ by localization, the viscosity property introduced in Definition \ref{def:visco} is a local property. Indeed, in order to check the viscosity property of $u$ at $(t,\omega)$, it suffices to know the value of $u^{t,\o}$ on $[0, \hc_\eps]$ for an arbitrarily small $\eps>0$. In particular, since $u$ and $\varphi$ are locally bounded, there is no integrability issue in the definition of the set of test functions $\underline{\Ac}^\Pc$ and $\overline{\Ac}^\Pc$. 
}
\end{Remark}

\subsection{Semijets definition and punctual differentiability}

Similar to the standard notion of viscosity solutions in finite-dimensional spaces, we will now prove that we may reduce our Definition \ref{def:visco} to paraboloids:
 \b*
 \phi^{q,p,\gamma}_s(\omega)
 :=
 qs+p\cdot\omega_s+\frac12\gamma:\omega_s\o_s^{\rm T},
 &s\in[0,T-t],&
 \omega\in\Omega,
 \e*
for some $(q,p,\gamma)\in\R\times\R^d\times\S_d$. We then introduce the corresponding subjet and superjet:
 \b*
 \underline{\Jc}^\Pc u_t(\omega)
 &:=&
 \big\{ (q,p,\gamma)\in\R\times\R^d\times\S_d: 
           \phi^{q,p,\gamma}\in\underline{\Ac}^\Pc u_t(\omega)
 \big\},
 \\
 \overline{\Jc}^\Pc v_t(\omega)
 &:=&
 \big\{ (q,p,\gamma)\in\R\times\R^d\times\S_d: 
           \phi^{q,p,\gamma}\in\overline{\Ac}^\Pc v_t(\omega)
 \big\}.
 \e*

\begin{Proposition}
Let $\Pc\subset\Pc_L$ for some $L>0$. A process $u\in C^0(\ol\Theta)$ is a $\Pc-$viscosity subsolution of \eqref{PPDE} if and only if:
 \be\label{def:viscojets}
 -q-G(t,\omega,u_t(\omega),p,\gamma)\le 0
 &\mbox{for all}&
 (t,\omega)\in \Theta,~(q,p,\gamma)\in\underline{\Jc}^\Pc u_t(\omega).
 \ee
The corresponding statement holds for supersolutions.
\end{Proposition}

\proof We focus on the nontrivial direction, assuming that \eqref{def:viscojets} holds. For $(t,\omega)\in\Theta$ and $\varphi\in\underline{\Ac}^\Pc u_t(\omega)$, we have to prove that $-q-G(t,\omega,u_t(\omega),p,\gamma)\le 0$, where 
 \b*
 q:=\partial_t\varphi(t,\omega),
 &p:=\partial_\omega\varphi(t,\omega),&
 \gamma:=\partial^2_{\omega\omega}\varphi(t,\omega).
 \e*
Without loss of generality, we take $(t,\omega)=(0,0)$. For $\eps>0$, we denote $q_\eps:=q+\eps(1+2L)$, and $\phi:=\phi^{q_\eps,p,\gamma}$.  By the smoothness of $\varphi$, we may find $\delta_\eps>0$, such that
 \b*
 |\partial_t\varphi-q|\le\eps,
 &|\partial_\omega\varphi-p-\gamma\omega_t|\le\eps,&
 \mbox{and}~~|\partial^2_{\omega\omega}\varphi-\gamma|\le\eps
 ~~
 \mbox{on}~~\Qc_\eps:=\{(t,\omega):t\le\delta_\eps,|\omega|_t\le\delta_\eps\}.
 \e*
Let $\hc$ be the stopping time corresponding to $\varphi$, and set $\hc_\eps:=\hc\wedge\inf\{t>0:(t,\omega)\not\in\Qc_\eps\}$. Then, for all stopping time $\tau\in\Tc_0$:
 \b*
 (\phi-u)_0-\underline{\Ec}^\Pc\big[(\phi-u)_{\tau\wedge\hc_\eps}\big]
 &\!\!\!\!\!\le&\!\!\!\!\!
 (\varphi-u)_0-\underline{\Ec}^\Pc\big[(\varphi-u)_{\tau\wedge\hc_\eps}\big]
 +\overline{\Ec}^\Pc\big[(\varphi-\varphi_0-\phi)_{\tau\wedge\hc_\eps}\big]
 \\
 \!\!\!\!\!
 &\!\!\!\!\!\le&\!\!\!\!\!
 \overline{\Ec}^\Pc\Big[\int_0^{\tau\wedge\hc_\eps}\!\!\!\!
                                  (\partial_t\varphi_s\!-\!q_\eps)ds
                                  +(\partial_\omega\varphi_s\!-\!p\!-\!\gamma B_s)dB_s
                                  +(\partial^2_{\omega\omega}\varphi_s\!-\!\gamma)d\langle B\rangle_s
                          \Big].
 \e*
Since $\Pc\subset\Pc_L$, it follows that the integral term inside the nonlinear expectation $\overline{\Ec}^\Pc$ is non-positive, implying that $(\phi-u)_0-\underline{\Ec}^\Pc\big[(\phi-u)_{\tau\wedge\hc_\eps}\big]\le 0$. Consequently $(q_\eps,p,\gamma)\in\underline{\Jc}^\Pc u_0$ and therefore $-q_\eps-G(t,\omega,u_t(\omega),p,\gamma)\le 0$ by \eqref{def:viscojets}. The required result follows by sending $\eps\searrow 0$.
\ep

\vspace{5mm}

\begin{Proposition}\label{prop:jetadditive}
For $u^i,v^i\in C^0(\bar\Theta)$, $i=0,1$, we have
 \b*
 \underline{\Jc}^\Pc u^0_t(\omega)
 +\underline{\Jc}^\Pc u^1_t(\omega)
 \subset
 \underline{\Jc}^\Pc(u^0+u^1)_t(\omega)
 &\mbox{and}&
 \overline{\Jc}^\Pc v^0_t(\omega)
 +\overline{\Jc}^\Pc v^1_t(\omega)
 \subset
 \overline{\Jc}^\Pc(v^0+v^1)_t(\omega)
 \e*
\end{Proposition}

\proof
We only report the argument for the subjets. Let $\theta^i=(q^i,p^i,\gamma^i)\in\underline{\Jc}^\Pc u^i_t(\omega)$, $i=0,1$. By definition, this means that the corresponding paraboloids $\phi^{\theta^i}\in\underline{\Ac}^\Pc u^i_t(\omega)$, i.e. there is $\hc^i\in\Hc^\Pc_+$ such that
 \b*
 -u^i_t
 \le
 \E^\P\big[(\phi^{\theta^i}-(u^i)^{t,\omega})_{\tau\wedge\hc^i}
         \big]
 &\mbox{for all}&
 \tau\in\Tc~\mbox{and}~\P\in\Pc.
 \e*
With $\hc:=\hc^0\wedge\hc^1\in\Hc^\P$, this implies that
 \b*
 -(u^0+u^1)_t
 \le
 \E^\P\big[\big(\phi^{\theta^0}+\phi^{\theta^1}-(u^0+u^1)^{t,\omega}\big)_{\tau\wedge\hc}
         \big]
 &\mbox{for all}&
 \tau\in\Tc~\mbox{and}~\P\in\Pc.
 \e*
Since $\phi^{\theta^0}+\phi^{\theta^1}=\phi^{\theta^0+\theta^1}$, this shows that $\theta^0+\theta^1\in\underline{\Jc}^\Pc(u^0+u^1)_t(\omega)$.
\ep

\subsection{Punctual differentiability}

The following notion is adapted from Caffarelli and Cabre \cite{CC}.

\begin{Definition}
A process $u$ is $\Pc-$punctually $C^{1,2}$ at some point $(t,\omega)\in\Theta$ if
 \b*
 \Jc^\Pc u_t(\omega)
 &:=&
 \mbox{cl}\big(\underline{\Jc}^\Pc u_t(\omega)\big)
 \cap
 \mbox{cl}\big(\overline{\Jc}^\Pc u_t(\omega)\big)
 \;\neq\;
 \emptyset.
 \e*
\end{Definition}

The next (immediate) result states that the viscosity property reduces to a pointwise property at points of punctual differentiability.

\begin{Proposition}\label{prop:punctualpointwise}
Assume that the nonlinearity $G$ is continuous in $(z,\gamma)$, and
let $u\in C^0(\ol\Theta)$ be a $\Pc-$viscosity solution of \eqref{PPDE}. Then, if $u$ is $\Pc-$punctually $C^{1,2}$ at some point $(t,\omega)\in\Theta$, we have
 \b*
 -q-G(t,\omega,u(t,\omega),p,\gamma)=0
 &\mbox{for all}&
 (q,p,\gamma)\in\Jc^\Pc u_t(\omega).
 \e*
\end{Proposition} 

For our subsequent analysis, we need the following additivity property of punctual differentiability, which is a direct consequence of Proposition \ref{prop:jetadditive}.

\begin{Proposition}\label{prop:punctualadditive}
Let $u,v$ be two processes which are $\Pc-$punctually $C^{1,2}$ at some point $(t,\omega)\in\Theta$. Then, $u+v$ is $\Pc-$punctually $C^{1,2}$ at $(t,\omega)$, and 
 \b*
 \Jc^\Pc u_t(\omega)
 +\Jc^\Pc u_t(\omega)
 &\subset&
 \Jc^\Pc (u+v)_t(\omega).
 \e*
\end{Proposition}

\subsection{Consistency of path-dependent viscosity solutions}

We conclude this definition subsection by proving consistency of our notion of viscosity solution with classical solutions.

\begin{Proposition}
Let $G$ be continuous, elliptic and uniformly $L_0-$Lipschitz-continuous in $(y,z,\gamma)$. Let $\cP_{L_0}\subset\cP\subset\cP_L$ for some $L\ge L_0$. Then, for $u\in C^{1,2}_\Pc(\Theta)$, the following are equivalent:
\\
{\rm (i)} $u$ is a $\Pc-$classical subsolution (reps. supersolution) for some $L>0$,
\\
{\rm (ii)} $u$ is a $\Pc-$viscosity subsolution (reps. supersolution).
\end{Proposition}

\proof We only report the proof of the subsolution property. The supersolution property follows by the same line of argument. If $u$ is a 
$\Pc-$viscosity subsolution and $u\in C^{1,2}_\Pc(\Theta)$, then it is clear that $u^{t,\o}\in\underline{\cA}^\Pc u_t(\o)$ for all $t<T$ and $\o\in\O$, and therefore $u$ is a $\Pc-$classical subsolution.

We next assume that $u$ is a classical subsolution, and we assume to the contrary that $5c:=-\partial_t\varphi-G(.,u,\partial_\o\varphi,\partial^2_{\o\o}\varphi)>0$ for some $t<T$, $\o\in\O$, and $\varphi\in\underline{\Ac}^\Pc u_t(\o)$. Without loss of generality, we may assume $(t,\o)=(0,0)$. Let $\bar\alpha\in\dbR^d$, $\bar\beta\in\dbS_d$ be arbitrary constants with $|\bar\alpha|\le L_0$ and $\frac12 \mbox{Tr}[\bar\beta^2]\le L_0$, to be fixed later, and denote by $\bar\dbP:=\dbP^{\bar\alpha,\bar\beta}$ the corresponding probability measure in $\cP$, and $\bar\Lc:=\bar\alpha\cdot\partial_\o+\frac12\bar\beta^2:\partial^2_{\o\o}$. By the continuity of $G$, and the fact that $u,\varphi\in C^{1,2}$,
 \beaa
 &-\partial_t\varphi-G_0(u_0,\partial_\o\varphi_0,\partial^2_{\o\o}\varphi_0)\ge 4c
 ,~
 |\bar\Lc\varphi-\bar\Lc\varphi_0|\le c,&
 \\
 &\mbox{and}~
 \big|G(u,\partial_\o u,\partial^2_{\o\o}u)
        -G_0(u_0,\partial_\o u_0,\partial^2_{\o\o}u_0)\big|\le c,
 ~|\bar\Lc u-\bar\Lc u_0|\le c,
 &~~\mbox{on}~
 [0,\hc_\eps],
 \eeaa
for $\eps>0$ sufficiently small, where $\hc_\eps:=\eps\wedge\inf\{s>0:|\o_s|\ge\eps\}$. Since $u$ is a $\Pc-$classical subsolution, we compute for every $\tau\in\Tc$ that
 \beaa
 (\varphi-u)_0
 -\E^{\bar\P}\big[(\varphi-u)_{\tau\wedge\hc_\eps}\big]
 &=&
 \E^{\bar\P}\Big[\int_0^{\tau\wedge\hc_\eps}d(u-\varphi)_s\Big]
 \\
 &=&
 \E^{\bar\P}\Big[\int_0^{\tau\wedge\hc_\eps}
                          \big\{\partial_t(u-\varphi)_s
                                   +\bar\Lc(u-\varphi)_s
                          \big\}ds
                  \Big]
 \\
 &\ge&
 \E^{\bar\P}\Big[\int_0^{\tau\wedge\hc_\eps}\!\!\!\!
                          \big\{G_0(u_0,\partial_\o\varphi_0,\partial^2_{\o\o}\varphi_0)
  \\ && \hspace{14mm}
                                   -G_0(u_0,\partial_\o u_0,\partial^2_{\o\o}u_0)
                                  +\bar\Lc(u-\varphi)_0
                          \big\}ds
                  \Big]
 \!+\!c\bar\P[\tau\wedge\hc_\eps].
 \eeaa
By the definition of $\cP$, we may find $\bar\alpha$ so that $G_0(u_0,\partial_\o\varphi_0,\partial^2_{\o\o}\varphi_0)-G_0(u_0,\partial_\o u_0,\partial^2_{\o\o}u_0)+\bar\Lc(u-\varphi)_0=0$. Then,whenever $\tau>0$, $\bar\P-$a.s., we have $(\varphi-u)_0>\E^{\bar\P}\big[(\varphi-u)_{\tau\wedge\hc_\eps}\big]$, contradicting the fact that $\varphi\in\underline{\Ac}^\Pc u_0$.
\ep

\section{Wellposedness of the path-dependent heat equation}

In this section, we consider the heat equation
 \be\label{p-heat}
 -\partial_tu-\frac12\mbox{Tr}\big[\partial^2_{\omega\omega}u]
 &=&
 0
 \ee
where, for simplicity, the diffusion matrix is taken to be the identity matrix. We recall that $\P_0$ denotes the Wiener measure. In addition to the previous notations, we denote $\F^*$ as the filtration augmented by all $\P_0$-null sets. Also, denote $\cT_*$ (resp. $\cT_*^t$) as the set of all $\F^*$-stopping times taking values in $[0,T]$ (resp. $[t,T]$).  In this section, we take 
 \b*
 \Pc:=\{\P_0\}
 &\mbox{and}&
 \Hc^\Pc_+:=\Tc^+.
 \e*
In this section about the heat equation, the relevant space for our comparison result is
 \beaa
 C^0_{2,\P_0}(\Theta,\dbR)
 &:=&
 \Big\{ u\in C^0(\bar\Theta,\dbR):\E^{\P_0}\Big[\sup_{t+s\le T}\big|u^{t,\o}_s\big|^2\Big]<\infty
                                              ~\mbox{for all}~(t,\o)\in\Theta
 \Big\}.
 \eeaa

\subsection{Facts from optimal stopping theory}

Let $X\in C^0_{2,\P_0}(\Theta,\dbR)$. Our main result uses the Snell envelope characterization of the optimal stopping stopping problem:
 \b*
 V_0
 := 
 \sup_{\tau\in\Tc_*} \E^{\P_0}[X_\tau],
 \e* 
The standard characterization of this problem uses the dynamic formulation of this problem:
 \b*
 Y^0_t
 &:=&
 \esup_{\tau\in\Tc_*^t} \E^{\P_0}\big[X_{\tau\wedge T}\big|\Fc_t\big],
 ~~0\le t\le T,
 \e*
so that $Y^0_0=V_0$ by the Blumenthal zero-one law. In this context, an optimal stopping rule is well-known to be defined by the first hitting time
 \b*
 \tau^* 
 &:=&
 \inf\big\{t\ge 0: Y^0_t=X_t\big\}.
 \e*
In addition to the standard result, we need an additional refinement by introducing the variable:
 \b*
 \Yc_\t(\omega)
 :=
 \sup_{\th\in\Tc_*}
 \E^{\P_0}\big[X^{\t(\o),\omega}_\th\big],
 &\mbox{for all}&
 \t\in\Tc_*,\o\in\Omega.
 \e*

\begin{Theorem}\label{thm:optimalstop0}
Let $X\in C^0_{2,\P_0}(\Theta,\dbR)$. Then, there exists an $\F-$adapted version $Y$ of $Y^0$ satisfying:
 \b*
 Y_{\tau\wedge T}
 =
 \Yc_{\tau\wedge T},
 &\P_0-\mbox{a.s. for all}&
 \tau\in\Tc_*.
 \e* 
Moreover, $Y$ is a pathwise continuous $\P_0-$supermartingale, $Y_{\wedge\tau^*}$ is a $\P_0-$martingale, and $\tau^*$ is an optimal stopping rule.
\end{Theorem}

This result follows from the more general Theorem \ref{thm:optimalstop} below.

\subsection{Existence, comparison, and uniqueness}

\begin{Definition}\label{def:regularsubmart}
An $\F-$progressively measurable process $m$ is a regular $\P_0-$submartingale (resp. supermartingale) if, for any $(t,\omega)\in\Theta$, we have 
 \b*
 m_t(\omega) - \E^{\P_0}\big[m^{t,\omega}_\tau\big] \le 0 
 &\mbox{(resp. $\ge 0$)}&
 \mbox{for all}~~\tau\in\Tc.
 \e*
\end{Definition}

\begin{Lemma}\label{lem:regsubmart}
Let $u\in C^0_{2,\P_0}(\Theta,\dbR)$, $(t,\omega)\in\Theta$, and $\hc\in\Tc^+$, be such that $u_t(\omega) > \E^\P\big[u^{t,\omega}_\hc\big]$. Then, 
 \b*
 0\in\underline{\Ac}^{\P_0}u_{t+t^*}(\o\otimes_t\o^*)
 &\mbox{for some}&
 (t^*,\o^*)~~\mbox{with the localization}~~\hc^*:=\hc^{t^*,\o^*}-t^*\in\Tc^+. 
 \e*
\end{Lemma}

\proof
Without loss of generality, we may assume that $(t,\o)=(0,0)$. Consider the optimal stopping problem $V_0:=\sup_{\tau\in\Tc_*}\E^{\P_0}\big[u_{\tau\wedge\hc}\big]$. Set $X_s:=u_{s\wedge\hc}$ and let $Y$ be the $\F-$adapted Snell envelope as introduced in Theorem \ref{thm:optimalstop0}, $\tau^*$ the corresponding optimal stopping rule. From the strict inequality $u_0 > \E^{\P_0}\big[u_\hc\big]$, it follows that $\P_0[\tau^*<\hc]>0$. By Theorem \ref{thm:optimalstop0}, we also have $Y_{\tau^*}=\Yc_{\tau^*}$, $\P_0-$a.s. We may then find $\omega^*$ such that $t^*:=\tau^*(\omega^*)<\hc(\omega^*)$, and:
 \b*
 u_{t^*}(\omega^*)
 &=&
 Y_{t^*}(\o^*)
 \;=\;
 \max_{\tau\in\Tc} 
 \E^{\P_0}\big[(u_{\hc\wedge\cdot})^{t^*,\omega^*}_\t
                                    \big],
 \e*
By definition of $\underline{\Ac}^{\P_0}u$, this is exactly the required result.
\ep

\vspace{5mm}

The main result of this section is the following.

\begin{Theorem}\label{thm:subsol-submart}
For a process $u\in C^0_{2,\P_0}(\Theta,\dbR)$, the following are equivalent:
\\
{\rm (i)} $u$ is a regular $\P_0-$submartingale (resp. supermartingale),
\\
{\rm (ii)} $u$ is $\P_0-$viscosity subsolution (resp. supersolution) of the path-dependent heat equation \eqref{p-heat}.
\end{Theorem}

\proof
(i) $\Longrightarrow$ (ii): For arbitrary $(t,\omega)\in\Theta$ and $\varphi\in\underline{\Ac}^{\P_0}u_t(\omega)$, we have for some $\hc\in\Tc_+$:
 \b*
 \varphi_0-u_t(\omega) \le \E^{\P_0}\big[\varphi_{\t\wedge\hc}-u^{t,\omega}_{\tau\wedge\hc}\big]
 &\mbox{for all}&
 \tau\in\Tc.
 \e*
For all $\eps>0$, define $\hc_\eps(\o^{'}):=\hc(\o^{'})\wedge\inf\{s\geq 0: |\o^{'}_s|\ge\eps\}$. Then, since $u$ is a regular $\P_0-$submartingale, it follows that
 \b*
 0 
 \;\ge\;
 u_t(\omega) - \E^{\P_0}\big[u^{t,\omega}_{\hc_\eps}\big]
 &\ge&
 \varphi_0 - \E^{\P_0}\big[\varphi_{\hc_\eps}\big]
 \;=\;
 \E^{\P_0}\Big[\int_0^{\hc_\eps}
                  (-\partial_t\varphi-\frac12\!\!:\!\partial^2_{\omega\omega}\varphi)_sds
           \Big]
 \e*
by the smoothness of $\varphi$. Sending $\eps\searrow 0$, we see that $(-\partial_t\varphi-\frac12\sigma^2\!\!:\!\partial^2_{\omega\omega}\varphi)_0\le 0$, as required.

(ii) $\Longrightarrow$ (i): Clearly, it is sufficient to prove that the process $\bar u:=u^\eps_t:=u_t+\eps t$ is a regular $\P_0-$submartingale for all $\eps>0$, as the required claim will follow by sending $\eps$ to zero. By (ii), we deduce immediately that $\bar u$ is a $\P_0-$viscosity subsolution of the equation $\eps-\partial_t\bar u-\frac12{\rm Tr}[\partial^2_{\omega\omega}\bar u] \le 0$ on $\Theta$. In particular, this implies that
 \be\label{0notinA}
 0\not\in\underline{\Ac}^{\P_0}\bar u_t(\omega)
 &\mbox{for all}&
 (t,\omega)\in\Theta.
 \ee
Suppose to the contrary that $\bar u$ is not a regular $\P_0-$submartingale, i.e. $\bar u_t(\omega)>\E^{\P_0}\big[\bar u^{t,\omega}_\hc\big]$ for some $(t,\omega)\in\Theta$ and $\hc\in\Tc_+$. Then, Lemma \ref{lem:regsubmart} induces a contradiction of \eqref{0notinA}. 
\ep

\vspace{5mm}

As an immediate consequence of Theorem \ref{thm:subsol-submart}, we obtain the wellposedness of the path-dependent heat equation.

\begin{Theorem}[Comparison and existence for the heat equation]~\\
{\rm (i)} Let $u,v\in C^0_{2,\P_0}(\Theta,\dbR)$ be $\P_0-$viscosity subsolution and supersolution, respectively, of the path-dependent heat equation \eqref{p-heat}, with $u_T\le v_T$ on $\Omega$. Then $u\le v$ on $[0,T]\times\Omega$.
\\
{\rm (ii)} For an $\Fc_T$ r.v. $\xi$ such that $u_t(\o):=\E^{\P_0}[\xi^{t,\o}]\in C^0_{2,\P_0}(\Theta,\dbR)$, the process $u$ is the unique $\P_0-$viscosity solution of the path-dependent heat equation \eqref{p-heat} with boundary condition $u_T=\xi$ on $\Omega$.
\end{Theorem}

\proof
(i) By Theorem \ref{thm:subsol-submart}, we have $u_t(\omega)\le \E^{\P_0}[(u_T)^{t,\omega}]$ and $\E^{\P_0}[(v_T)^{t,\omega}]\ge v_t(\omega)$ for all $(t,\omega)\in\Theta$. Then $u_T\le v_T$ on $\Omega$ implies that $u\le v$ on $[0,T]\times\Omega$.

(ii) Uniqueness is a direct consequence of the comparison result of (i). Clearly the process $u_t(\omega):=\E^{\P_0}\big[\xi^{t,\omega}\big]$ is uniformly continuous on $[0,T]\times\Omega$. Since $u$ is a $\P_0-$martingale, it follows from Theorem \ref{thm:subsol-submart} that it is both a viscosity subsolution and supersolution.
\ep

\vspace{5mm}

%
%

\section{Wellposedness of semilinear path-dependent PDEs}
\label{sect-comparisonslin}
\setcounter{equation}{0}

In this section, we consider the equation
 \be\label{p-semilinear}
 - \partial_tu
 - \frac12\mbox{Tr}\big[\partial^2_{\omega\omega}u\big]
 - F(.,u,\partial_\omega u)
 =
 0
 &\mbox{on}&
 \Theta.
 \ee
The nonlinearity $F:\Theta \times\R\times\R^d\longrightarrow\R$ is assumed to satisfy the following assumptions which consists with the general assumption as (\ref{L0-Lip}).

\begin{Assumption}\label{assum:gen}
The nonlinearity $F: (t,\o,y,z)\in\Theta\times\mathbb{R}\times\mathbb{R}^d \longmapsto F(t,\o,y,z)\in\mathbb{R}$ satisfies the following conditions:
\\
{\rm (i)}\quad $F$ is uniformly continuous in $(t,\o)$,
\\
{\rm (ii)}\quad $F$ is uniformly $L_0-$Lipschitz continuous in $(y,z)$, for some $L_0\ge 0$, i.e. 
 \beaa 
 |F(\cdot,y,z)-F(\cdot,y',z')|
 \le 
 L_0\left(|y-y'|+|z-z'|\right)
 &\mbox{for all}&
 y,y'\in\mathbb{R},~z,z'\in\mathbb{R}^d. 
 \eeaa
{\rm (iii)}\quad The process $F(\cdot,0,0)$ is bounded.
\end{Assumption}

For all bounded $\F-$progressively measurable process $\lambda$, we denote:
 \b*
 d\P_\lambda := Z^\lambda_T\cdot d\P_0
 &\mbox{on $\Fc_T$, where}&
 Z^\lambda_T:=e^{\int_0^T\lambda_t\cdot dB_t-\frac12\int_0^T |\lambda_t|^2dt}.
 \e*
In this section, we take 
 \be\label{Psemilinear}
 \Pc
 :=
 \Big\{\P_\lambda:~\lambda~\mbox{bounded by}~L
 \Big\},
 \ee
where $L\ge L_0$ is arbitrary. Notice that $\P_0$ is a dominating measure for the family $\Pc$. For simplicity, we say $\P_\l\in\Pc$ by implying that $\l$ is the corresponding bounded process. Similar to the section of the heat equation,  we denote $\F^*$ as the filtration augmented by all $\P_0$-null sets. Also, we consider the set of localizing stopping times as:
 \b*
 \Hc^\Pc_+
 &:=&
 \Tc^+.
 \e*
In this section about the semilinear equation, the relevant space for our comparison result is
 \beaa
 C^0_{2,\cP}(\Theta,\dbR)
 &:=&
 \Big\{ u\in C^0(\Theta,\dbR):\ol\cE^{\cP}\Big[\sup_{t+s\le T}\big|u^{t,\o}_s\big|^2\Big]<\infty
                                              ~\mbox{for all}~(t,\o)\in\Theta
 \Big\}.
 \eeaa

\subsection{Optimal stopping under dominated nonlinear expectation}

For $X\in C^0_{2,\cP}(\Theta,\dbR)$, we consider the optimal stopping stopping problem under dominated nonlinear expectation:
 \b*
 V_0
 := 
 \sup_{\tau\in\Tc_*} \overline{\Ec}^\Pc[X_\tau].
 \e*
The corresponding dynamic formulation is defined by:
 \b*
 Y^0_t
 &:=&
 \esup_{\tau\in\Tc_*^t} \overline{\Ec}^\Pc\big[X_{\tau\wedge T}\big|\Fc_t\big],
 ~~0\le t\le T.
 \e*
with first hitting time:
 \b*
 \tau^* 
 &:=&
 \inf\big\{t\ge 0: Y^0_t=X_t\big\}.
 \e*
Since the dominating measure $\P$ satisfies the Blumenthal zero-one law, it follows that $Y^0_0=V_0$. We also introduce the pointwise optimal stopping problem:
 \b*
 \Yc_t(\omega)
 :=
 \sup_{\tau\in\Tc_*}
 \overline{\Ec}^{\Pc}\big[X^{t,\omega}_{\tau\wedge (T-t)}\big],
 &\mbox{for all}&
 (t,\omega)\in\bar\Theta.
 \e*

\begin{Theorem}\label{thm:optimalstop}
Let $X\in C^0_{2,\cP}(\Theta,\dbR)$. Then, there exists an $\F-$adapted version $Y$ of $Y^0$ satisfying:
\\
{\rm (i)} for all $\tau\in\Tc$, we have $Y_{\tau\wedge T}=\Yc_{\tau\wedge T}$, $\P_0-$a.s.
\\
{\rm (ii)} $Y$ is a pathwise continuous $\P-$supermartingale for all $\P\in\Pc$, and $\tau^*$ is an optimal stopping rule,
\\
{\rm (iii)} $Y_t=\esup_{\tau\in\Tc_*^t}\E^{\P^*}\big[X_{\tau\wedge T}|\Fc_t\big]$ for all $t\in[0,T]$, $\P_0-$a.s. for some $\P^*\in\Pc$, and
 \b*
 &Y=Y_0+M^*-K^*
 ~~\mbox{with}~~M^*_0=K^*_0=0,~~\mbox{and}~~
 \int(Y-X)dK^*=0,~~\P_0-\mbox{a.s.}&
 \e*
for some pathwise continuous martingale $M^*$ and predictable nondecreasing process $K^*$.
\end{Theorem}

This result can be proved by referring to the corresponding literature in the theory of reflected backward stochastic differential equations, see Remark 7.3 in \cite{ETZ1}. For the convenience of those readers who are not familiar with this literature, we report in Section \ref{sect:optimalstop} a proof purely based on arguments from optimal stopping theory.

\subsection{Punctual smoothness of submartingales}

In this subsection, we prove that a process $u\in C^0_{2,\cP}(\Theta,\dbR)$ which is $\P-$submartingale for some $\P\in\Pc$ is punctually $C^{1,2}_\Pc-$ Leb$\otimes\P-$a.e. This is our natural extension of the well-known result that any non-decreasing function is differentiable a.e. and our proof builds on the corresponding standard results in analysis that we quickly review. For a function $f:[0,T]\longrightarrow\R$ with finite variation, we use the following notations for the left-semigradients:
 \b*
 \dot{\underline{f}}^\ell(t):=\liminf_{\eps\uparrow 0}\frac{f(t+\eps)-f(t)}{\eps}
 &\mbox{and}&
 \dot{\overline{f}}^\ell (t):=\limsup_{\eps\uparrow 0}\frac{f(t+\eps)-f(t)}{\eps}.
 \e*
The right-semigradients $\dot{\underline{f}}^r$ and $\dot{\overline{f}}^r$ are defined similarly by sending $\eps\downarrow 0$. The function $f$ is differentiable at a point $t$ if 
 \b*
 \dot{f}(t)=\lim_{\eps\to 0}\frac{f(t+\eps)-f(t)}{\eps}
 &\mbox{exists, and therefore}&
 \dot{f}(t)
 =
 \dot{\underline{f}}^\ell(t)
 =
 \dot{\overline{f}}^\ell(t)
 =
 \dot{\underline{f}}^r(t)
 =
 \dot{\overline{f}}^r(t).
 \e*
Our smoothness results uses crucially the two following properties:
\begin{enumerate}
\item[$\mathbf{\mbox{\bf FV}_1}$] The set of points of differentiability of $f$ has full Lebesgue measure. 
\item[$\mathbf{\mbox{\bf FV}_2}$] If $f$ is absolutely continuous, then $\lim_{\eps\to 0}\frac{1}{\eps}\int_t^{t+\eps}|\dot{f}(s)-\dot{f}(t)|ds=0$, Leb-a.e. on $[0,T]$.
\end{enumerate}
For a subset $\Theta_0\subset \bar\Theta$, we denote $\T^{\Theta_0}:=\{t:(t,\o)\in\Theta_0~\mbox{for some}~\o\in\O\}$ and $\O^{\Theta_0}_t:=\big\{\o:(t,\o)\in\T^{\Theta_0}\big\}$.

\begin{Theorem}\label{Thm:pathwisesmooth}
Let $\P_\theta\in\Pc$ and $u\in C^0_{2,\cP}(\Theta,\dbR)$ be $\P_
\theta-$submartingale. Then $u$ is $\Pc-$punctually $C^{1,2}$ on $\Theta_0$, for some $\Theta_0$ with 
 \bea\label{section-ae}
 \mbox{Leb}\big[\T^{\Theta_0}\big]=T
 &\mbox{and}&
 \P_0\big[\O^{\Theta_0}_t\big]=1
 ~~\mbox{for all}~~t\in\T^{\Theta_0}.
 \eea
\end{Theorem}

\no {\bf Sketch of the proof.}\quad For a proof in more details, we refer to \cite{RTZ}.
We proceed in two steps.
\\
\underline{\it Step 1}: By the Doob-Meyer decomposition, we have $u=u_0+M+A$, $\P_0-$a.s. for some $\P_\theta-$martingale $M$ and nondecreasing predictable process $A$, with $M_0=A_0=0$. Then, process
$M^0:=M-\int_0^\cdot \theta_sd\langle M,B\rangle_s$
 defines a $\P_0$-martingale.

Since all $\P_0$-martingale have the martingale representation, it follows that $t\longmapsto H_t:=\langle M,B\rangle_t=\langle M^0,B\rangle_t$ is absolutely continuous on $[0,T]$, $\P_0-\mbox{a.s.}$, i.e.
 $$
 h_t:=\dot{\overline{H}}^\ell_t=\dot{H}_t
 ~\mbox{for a.e.}~t\in[0,T],
 ~\P_0-\mbox{a.s.}
 $$
By the above property {\bf FV}$\mathbf{_2}$ together with the Fubini theorem, we see that
 \be
 &\mbox{Leb}\otimes\P_0[\Theta_1]=T
 ~~\mbox{where}~~
 \Theta_1:=\big\{(t,\omega): \lim_{\eps\to 0}\frac{1}{\eps}\int_t^{t+\eps}|h_s-h_t|ds=0
                      \big\}.&
 \ee

Further, applying property {\bf FV}$\mathbf{_1}$ to the finite variation process $A^\theta:=A+\int_0^.\theta_s dH_s$, and using again the Fubini theorem, we see that:
 \be
 \mbox{Leb}\otimes\P_0[\Theta_2]=T
 &\mbox{where}&
 \Theta_2:=\big\{(t,\omega): a_t(\o):=\dot{A}^\theta_t(\omega)~\mbox{exists}
                      \big\}.
 \ee
\underline{\it Step 2}: In this step, we prove that for $(t,\omega)\in\Theta_0:=\Theta_1\cap\Theta_2$.
 \be\label{pC12-i}
 (q^\eps,p,0)\in\underline{\Jc}^\Pc u_t(\omega),
 &\mbox{where}&
 q^\eps:=a_t(\omega)-\eps(1+L),~~p:=h_t(\omega).
 \ee
We define
 \b*
 &\hc(\o^{'})
 :=
 \inf\big\{s>0: (A^\theta)^{t,\o}_s(\o^{'})-A^\theta_t(\o)\le (a_t(\o)-\eps)s
                     ~\mbox{or}~
                     \int_0^s |h^{t,\o}_r(\o^{'})-h_t(\o)|dr\ge\eps s
      \big\}.&
 \e* 
 Since $\o\in\Theta_0$, we have $\hc\in\Tc^+$. Also, note that
$M^{\l,t}_s:=M^{t,\o}_s-M_t(\o)-\int_0^s(\theta-\l)_rh_rdr$ defines a $\P_\l$-martingale. Further, rewriting the Doob-Meyer decomposition, we have 
$$u^{t,\o}_s=u_t(\o)+(A^\theta)^{t,\o}_s-A^\theta_t(\o)+M^{\l,t}_s-\int_0^s\l_rh_rdr,\ \P_\l\mbox{-a.s.}$$
So, for all $\tau\in\Tc$, $\P_\lambda\in\Pc$:
 \b*
 \E^{\P_\l}\big[\big(\phi^{q^\eps,p,0}-u^{t,\omega}\big)_{(\tau\wedge\hc)^{t,\omega}}                                         
                              \big]
 &\!\!\!=&\!\!\!
 -u_t(\omega)
 +
 \E^{\P_\l}\Big[(a_t(\omega)-\eps)(\tau\wedge\hc)
                                     -A^\theta_{\tau\wedge\hc}+A^\theta_t(\o)
                                     \\ &&\hspace{25mm} 
         -\eps L (\tau\wedge\hc)+\int_0^{\tau\wedge\hc} (h_s-h_t)\lambda_sds                         
                              \Big]
 \le
 -u_t(\omega),
 \e*
by the definition of $\hc$. Then \eqref{pC12-i} holds.
\\
\underline{\it Step 3}: From the previous step, it follows that $\big(a_t(\omega),h_t(\omega),0\big)\in\mbox{cl}\big(\underline{\Jc}^\Pc u_t(\omega)\big)$. By a similar argument, we may show that $\big(a_t(\omega),h_t(\omega),0\big)\in\mbox{cl}\big(\overline{\Jc}^\Pc u_t(\omega)\big)$. Consequently, $\big(a_t(\omega),h_t(\omega),0\big)\in\Jc^\Pc u_t(\omega)$, and $u$ is punctually $C^{1,2}_\Pc$.  
\ep

\subsection{Comparison}
In this subsection, we are going to show the comparison principle for the semilinear path-dependent equation.

\begin{Theorem}\label{thm:comparison}
Let Assumption \ref{assum:gen} hold true. Let $u,v\in C^0_{2,\cP}(\Theta,\dbR)$ be $\cP$-viscosity subsolution and  supersolution, respectively, of the equation (\ref{p-semilinear}). Assume further that $u_{T}\leq v_{T}$ on $\O$. Then $u\le v$ on $\bar\Theta$.
\end{Theorem}

To show Theorem \ref{thm:comparison}, we need some preparation. The following lemma is the analog of Lemma \ref{lem:regsubmart} in the context of the semilinear path-dependent PDEs. We omit the proof, since it is similar to that of Lemma \ref{lem:regsubmart}.

\begin{Lemma}\label{lem:maximumpoint}
Let $u\in  C^0_{2,\cP}(\Theta,\dbR)$, $(t,\omega)\in\Theta$, and $\hc\in\Tc^+$, be such that $u_t(\omega) > \ol\Ec^{\Pc}\big[u^{t,\omega}_\hc\big]$. Then, 
 \b*
 0\in\underline{\Ac}^{\Pc}u_{t+t^*}(\o\otimes_t\o^*)
 &\mbox{for some}&
 (t^*,\o^*)~~\mbox{with the localization}~~\hc^*:=\hc^{t^*,\o^*}-t^*\in\cT^+. 
 \e*
\end{Lemma}

The next main ingredient is the partial comparison result.

\begin{Proposition}\label{prop:partialcomparison}
In the setting of Theorem \ref{thm:comparison}, assume in addition that $v\in C^{1,2}_{\Pc}(\Theta)$. Then $u\le v$ on $\Theta$.
\end{Proposition}

\proof
First, by possibly transforming the problem to the comparison of $\tilde u_t:=e^{\lambda t}u_t$ and $\tilde v_t:=e^{\lambda t}v_t$, it follows from the Lipschitz property of the nonlinearity $F$ in $y$ that we may assume without loss of generality that $F$ is decreasing in $y$.

Suppose to the contrary that $c:=(u-v)_t(\o)>0$ at some point $(t,\o)\in\Theta$. Let $c_0:=\frac{c}{2T}$, and define $f_s:=(u-v)^+_s+c_0(s-t)$, $s\in[t,T]$. Since $(u-v)_T\le 0$, it follows that $f_t(\o)>\ol\cE^{\cP}[f^{t,\o}_{T-t}]$. By Lemma \ref{lem:maximumpoint}, we may find a point $(t^*,\o^*)$ such that $t^*\in[t,T)$ and $0\in\ul\cA^{\cP}f_{t^*}(\o^*)$. In particular, this implies that 
 \beaa
 -(u-v)^+_{t^*}(\o^*)-c_0(t^*-t)
 &\le&
 \ul{\cE}^{\cP}\big[-((u-v)^+_T)^{t^*,\o^*}-c_0(T-t)\big]
 \;=\;
 -c_0(T-t),
 \eeaa 
so that $(u-v)^+_{t^*}(\o^*)\ge c_0(T-t^*)>0$. Then, since $(u-v)^+\ge u-v$, we deduce from $0\in\ul\cA^{\cP}f_{t^*}(\o^*)$ that
 \beaa
 (\varphi-u)_{t^*}(\o^*)
 \le
 \ul{\cE}^{\cP}\big[(\varphi-u)^{t^*,\o^*}_{\tau\wedge T}\big]
 &\mbox{for all}~\tau\in\cT,~\mbox{where}&
 \varphi_s(\o):= v_s(\o)-c_0(s-t).
 \eeaa
Since $v\in C^{1,2}_{\Pc}(\Theta)$, this means that $\f^{t^*,\o^*}\in\ul\cA^{\cP}u_{t^*}(\o^*)$. Then, by the viscosity subsolution property of $u$, and the classical supersolution property of $v$, we deduce that
 \beaa
 0 
 &\ge&
 \big\{- \partial_t\f - \frac12\mbox{Tr}\big[\partial^2_{\omega\omega}\f\big]-F(.,u,\partial_\o \varphi)\big\}(t^*,\o^*)
 \\
 &=&
 c_0+\big\{ - \partial_tv - \frac12\mbox{Tr}\big[\partial^2_{\omega\omega}v\big]-F(.,u,\partial_\o v)\big\}(t^*,\o^*)
 \\
 &\ge&
 c_0+\{F(.,v,\partial_\o v)
           -F(.,u,\partial_\o v)\}(t^*,\o^*)
 \;\ge\; 
 c_0,
 \eeaa
where the last inequality follows from the non-increase of  $F$ in $y$ and the fact that $u_{t^*}(\o^*)\ge v_{t^*}(\o^*)$. Since $c_0>0$, this is the required contradiction.
\ep

\begin{Lemma}\label{v a.s. PC2}
Under Assumption (\ref{assum:gen}), there is a constant $C$ such that
\\
{\rm (i)} the process $\big\{u_t+\int_0^t|u_s|ds+Ct,t\in[0,T]\big\}$ is a $\dbP^u$-regular submartingale, for some $\dbP^u\in\cP$,
\\
{\rm (ii)} the process $\big\{v_t-\int_0^t|v_s|ds-Ct,t\in[0,T]\big\}$ is a $\dbP^v$-regular supermartingale, for some $\dbP^v\in\cP$,
\\
{\rm (iii)} $u$ and $v$ are $\cP-$punctually $C^{1,2}$ on $\Theta^u$ and $\Theta^v$, respectively, where $\Theta^u$ and $\Theta^v$ satisfy \eqref{section-ae}.
\end{Lemma}

\proof
Assertion (iii) is a direct consequence of (i) and (ii) together with Theorem \ref{Thm:pathwisesmooth}. By Assumption \ref{assum:gen}, we may find a constant $C$ such that:
 \beaa
 |F(t,\o,y,z)|
 &\le& 
 C-1+L_0(|y|+|z|)
 \eeaa
Then, it is easy to verify that $\bar{u}_t:=u_t+Ct$ and $\bar{v}_t:=v_t-Ct$ are $\cP$-viscosity subsolution and supersolution, respectively of:
 \begin{equation}\label{eq f(0)=0}
 -\cL \bar{u}-L_0(|\bar u-Ct|+|\pa_\o \bar{u}|) +1 \le 0
 ~\mbox{and}~
 -\cL \bar{v}+L_0(|\bar v+Ct|+|\pa_\o \bar{v}|) -1 \ge 0
 ~~\mbox{on}~~[0,T)\times\O.
 \end{equation}
In the rest of this proof, we shall show that $\bar{u}$ and $\bar{v}$ are $\ol{\cE}^{\cP}$-regular submartingale and $\ul{\cE}^{\cP}$-regular supermartingale, respectively. In addition, we prove in Appendix (Proposition \ref{cE submartingale}) that a continuous $\ol{\cE}^{\cP}$-regular submartingale is a $\P$-submartingale for some $\P\in\cP$. This leads to the desired result.

We only prove that $\bar{u}$ is $\ol{\cE}^{\cP}$-regular submartingale, as the corresponding statement for $\bar{v}$ follows from the same line of argument. 

Suppose to the contrary that $\bar{u}_t(\o)>\ol{\cE}^{\cP}[\bar{u}^{t,\o}_\hc]$ for some $(t,\o)\in[0,T)\times\O$ and some stopping time $\hc\in\cT^+$. Then, it follows from Lemma \ref{lem:maximumpoint} that there exist $t^*$ and $\o^*$ such that $0\in\ul{\cA}^{\cP}\bar{u}_{t^*}(\o^*)$, i.e. there exists $\hc'\in\cT^+$ such that
$$-\bar u_{t^*}(\o^*)\geq \ol\cE^{\cP}\Big[-\bar u^{t^*,\o^*}_{\t\we\hc'}-L_0\int_0^{\t\we\hc'} |u^{t^*,\o^*}_s|ds\Big].$$
As a result, function $\f_t:=-L_0\int_0^t |u^{t^*,\o^*}_s|ds$ is in $\ul\cA^{\cP} u_{t^*}(\o^*)$. Since $\bar{u}$ is a $\cP$-viscosity subsolution of the left equation of  \eqref{eq f(0)=0}, this leads to
 $$
 L_0|u_{t^*}(\o^*)|-L_0|u_{t^*}(\o^*)|+1\leq 0,
 $$
which is the required contradiction, thus completing the proof of (i).
\ep

\vspace{5mm}

We are now ready for the key-result for the proof of the comparison result. We observe that this statement is an adaptation of the approach of Caffarelli and Cabre \cite{CC} to the comparison in the context of the standard theory of viscosity solutions in finite dimensional spaces. 

\begin{Proposition}\label{prop:wsubsol}
Let Assumption \ref{assum:gen} hold, and consider the $L$ in the definition of $\cP$ (recall that $L\ge L_0$).
Let $u,v\in C^0_{2,\cP}(\Theta,\dbR)$ be $\cP$-viscosity subsolution and supersolution, respectively, of the path-dependent PDE \eqref{p-semilinear}. Then, $w:=u-v$ is a $\cP$-viscosity subsolution of
\begin{equation}\label{equation w}
-\cL w(t,\o)-L|w_t(\o)|-L|\pa_\o w_t(\o)|\leq 0.
\end{equation}
\end{Proposition}

\no {\bf Sketch of Proof}\quad
Without loss of generality, we only check the viscosity property at
$(t,\o)=(0,0)$. For an arbitrary $(a,\beta,0)\in\ul{\cJ}^{\cP}(u-v)_0$, we have to prove that
 \bea\label{requiredu-vsubsol}
 -a-L\left|(u-v)_{0}\right|-L|\beta|
 &\le&
 0.
 \eea
{\bf 1.} Denote as usual by $\phi^{a,\beta}=\phi^{a,\beta
,0}$ the corresponding paraboloid process. By definition, there exists $\hc\in\cT^+$ such that
 \beaa
 c_0
 \;:=\;
 -(u-v)_{0}
 &=&
 \min_{\t\in\cT}\ul{\cE}^{\cP}\big[(\phi^{a,\beta}-u+v)_{\t\we\hc}\big].
 \eeaa
For $\delta>0$, $r>0$, and $\hc_r:=\hc\we\inf\{t:|\o_t|\geq r\}$, define the Snell envelop:
\beaa
 \hat{m}_{t}
 :=
 \einf_{\t\in\cT_*^t,\dbP\in\cP}\dbE^\dbP\left[m_{\t\we\hc_{r}}|\cF_t\right],
 &t\in[0,T],~\mbox{where}&
 m
:=
\phi^{a+\delta,\beta}-u+v.
 \eeaa
Clearly, 
 \bea\label{mmhat}
 m_{0}=c_0,~~\ul{\cE}^{\cP}\left[m_{\hc_r}\right]>c_0,
 &\hat{m}_{0}\leq m_{0},&
 \mbox{and}~~
 \hat{m}_{\hc_{r}}=m_{\hc_{r}},~~\dbP_0-\mbox{a.s.}
 \eea
Further, from Theorem \ref{thm:optimalstop}, we have that:
 \bea\label{lambda*}
 \hat{m}_{t}=\einf_{\t\in \cT_*^t}\dbE^{\dbP_{\l^*}}\left[m_{\t\we\hc_{r}}|\cF_t\right],
 &\dbP_0-\mbox{a.s. for some}&
 \|\lambda^*\|\le L.
 \eea
{\bf 2.} By classical optimal stopping theory, $\hat m$ is a $\dbP_{\l^*}-$submartingale with Doob-Meyer decomposition
 \beaa
 \hat{m}=\hat{m}_{0}+\hat{A}+\hat{M},
 &\mbox{with}&
 \hat{A}
 =
 \int_0^. \1_{\{m=\hat{m}\}}(s)d\hat{A}_{s},
 ~~\dbP_{\l^*}\mbox{-a.s.}
 \eeaa
for some $\dbP_{\l^*}-$martingale $\hat{M}$, and some nondecreasing process $\hat A$. In addition, we may prove that $\hat A$ is absolutely continuous $\P_0$-a.s. (see Step 4 in the proof of Proposition 7.3 in \cite{RTZ}).  Then, it follows from \eqref{mmhat} that:
 \beaa
 0
 < 
 \ul{\cE}^{\cP}\left[m_{\hc_{r}}-m_{0}\right]
 \leq 
 \ul{\cE}^{\cP}\left[\hat{m}_{\hc_{r}}-\hat{m}_{0}\right]
 &\leq& 
 \dbE^{\dbP_{\l^*}}\left[\int_{0}^{\hc_{r}}\1_{\{m=\hat{m}\}}(t)d\hat{A}_{t}\right]
 \\
&=&
 \lim_{M\to\infty}
 \dbE^{\dbP_{\l^*}}\left[\int_{0}^{\hc_{r}}\1_{\{m=\hat{m}\}}(t)
                                                              \1_{\{|\dot{\hat{A}}_t|\le M\}}\dot{\hat{A}}_t\;dt
                             \right]
 \\
 &\le& 
 \lim_{M\to\infty} M\;\dbE^{\dbP_{\l^*}}\left[\int_{0}^{\hc_{r}}\1_{\{m=\hat{m}\}}(t)dt\right].
\eeaa
This implies that $\mbox{Leb}\otimes\dbP_0\big[ t< \hc_{r},m=\hat{m}\big]>0$, 
so that, with the subsets $\Theta^u,\Theta^v$ from Proposition \ref{v a.s. PC2}, we have:
 \beaa
 \mbox{Leb}\otimes\dbP_0\big[ \{t\in[0,\hc_r),m=\hat{m}\}\cap \Theta^u\cap\Theta^v \big]
 &>&
 0.
 \eeaa
 Further, by taking in account (i) of Theorem \ref{thm:optimalstop}, we may find a point $(t^{*},\omega^{*})$ such that
 $$
 \begin{array}{c}
 \hc^{t^*,\o^*}_{r}-t^*\in\cT^+,~
 m_{t^{*}}(\o^{*})
 =\hat{m}_{t^{*}}(\o^{*})
 =\inf_{\t\in\cT_*}\ul{\cE}^{\cP}[m^{t^*,\o^*}_{\t\we (\hc_r^{t^*,\o^*}-t^*)}],
 \\
 \mbox{and}~u,v~\mbox{are}~\cP-\mbox{punctually}~C^{1,2}~\mbox{at}~(t^*,\o^*).
 \end{array}
 $$
{\bf 3.} By Proposition \ref{prop:punctualadditive}, it follows that $m$ is $\cP$-punctually $C^{1,2}$ at $(t^{*},\omega^{*})$, and $(a^m,\beta^m):=(a+\delta-a^u+a^v,\beta-\beta^u+\beta^v)\in\cJ^{\cP} m(t^*,\o^*)$ for any $(a^u,\beta^u)\in\cJ^{\cP} u(t^*,\o^*)$ and $(a^v,\beta^v)\in\cJ^{\cP} v(t^*,\o^*)$. 
Then, by using the viscosity subsolution property of $u$ together with Proposition \ref{prop:punctualpointwise} and the Lipschitz property of $F$ from Assumption \ref{assum:gen}, we see that:
 \beaa
 0 
 &\ge& 
 -a^u-F(t^{*},\o^{*},u_{t^{*}}(\o^{*}),\beta^u)
 \\
 &=&  (-a^v+a^m)-a-\d-F(t^{*},\o^{*},(u-v+v)_{t^{*}}(\o^{*}),\beta+\beta^v-\beta^m)
 \\
 &\ge& 
 a^m-L|\beta^m|-a-\delta-L\left|(u-v)_{t^{*}}(\o^{*})\right|-L|\beta|-a^v-F(t^{*},\o^{*},v_{t^{*}}(\o^{*}),\beta^v)
 \eeaa
We shall prove in Step 5 below that
 \bea\label{Claim2}
 a^m-L|\beta^m|
 &\ge& 0.
 \eea
Together with the viscosity supersolution property of $v$, this provides:
 \beaa
 0 
 &\ge&
 -a-\delta-L\left|(u-v)_{t^{*}}(\o^{*})\right|-L|\beta|.
 \eeaa
Since $t^*\rightarrow 0$ as $r\rightarrow0$, and $u,v\in C^0$, this provides $-a-\d-L\left|(u-v)_{0}|-L|\beta\right|\leq0$, which implies \eqref{requiredu-vsubsol} by sending $\delta\rightarrow0$.\\

\no {\bf 4.}\quad It remains to prove \eqref{Claim2}. For the sake of simplicity,
we set $t^{*}=0$. Recall that $(a^m,\beta^m)\in\cJ^{\cP} m_0$ and $m_{0}=\hat{m}_{0}=\inf_{\t\in\cT}\ul{\cE}^{\cP}[m_{\t\we H_r}]$. Suppose to the contrary that $a^m-L|\beta^m|<0$. Then, there exists $(\hat a,\hat\beta)\in\ul{\cJ}^{\cP} m_0$ such that $\hat a-L|\hat\beta|<0$. By definition of $\ul{\cJ}^{\cP} m_0$, we have
 \beaa
 m_{0}
 \;=\;
 \sup_{\t}\ol{\cE}^{\cP}\left[m_{\t\we\hat\hc}-\phi^{\hat a,\hat\beta}_{\t\we\hat\hc}\right]
 &\mbox{for some}&
 \hat\hc\in\cT^+~\mbox{with}~\hat\hc\le\hc_r.
 \eeaa
Then, considering the process $\l:=-L\text{sgn}(\hat\beta)$, we see that:
 \beaa
 \hat m_0
 = m_0
 \ge \dbE^{\dbP_\l}\big[m_{\hat\hc}-\phi^{\hat a,\hat\beta}_{\t\we\hat\hc}\big]
 = \dbE^{\dbP_\l}\big[m_{\hat\hc}\big]
    -(\hat a-L\hat\beta) \dbE^{\dbP_\l}[\hat\hc]
 > \dbE^{\dbP_{\l}}\big[m_{\hat\hc}\big].
 \eeaa
Since $\hat\hc\le\hc_r$ and $\dbP_\l\in\cP$, this is in contradiction with the definition of $\hat m_0$.
\ep

\vspace{5mm}

The previous proposition, together with the partial comparison result of Proposition \ref{prop:partialcomparison}, lead directly to the comparison result.

\vspace{5mm}

\no {\bf Proof of Theorem \ref{thm:comparison}}\quad
By Proposition \ref{prop:wsubsol}, $u-v\in C^0_{2,\cP}(\Theta,\dbR)$ is a $\cP-$viscosity subsolution of the path-dependent equation (\ref{equation w}). Clearly, $0$ is a classical supersolution of the same equation. Since $(u-v)_T\le 0$, we conclude from the partial comparison result of Proposition \ref{prop:partialcomparison} that $u-v\le 0$ on $\Theta$.
\ep

\subsection{Existence}

To establish an existence result of $\cP-$viscosity solutions of the equation \eqref{p-semilinear} under the above assumption \ref{assum:gen}, we consider a terminal condition defined by an $\Fc_T-$measurable r.v. $\xi$. Then, the PPDE \reff{p-semilinear} with terminal condition $u(T,\o) = \xi(\o)$ is closely related to the following backward stochastic differential equation (BSDE):
\bea
\label{BSDE}
Y^0_t = \xi + \int_t^T F(s, B, Y^0_s, Z^0_s) ds - \int_t^T Z^0_s dB_s,
~~ 0\le t\le T,~~\dbP_0\mbox{-a.s.}
\eea
We refer to the seminal paper by Pardoux and Peng \cite{PP} for the wellposedness of such BSDEs. On the other hand, for any $(t,\o)\in[0,T]\times\Omega$, by \cite{PP} the following BSDE on $[t,T]$ has a unique solution:
\bea
\label{Y0t}
Y^{0, t,\o}_s = \xi^{t,\o} + \int_s^T  F^{t,\o}(r,B^t, Y^{0,t, \o}_r, Z^{0,t,\o}_r) dr - \int_s^T Z^{0,t,\o}_rdB^t_r,~~\dbP_0^t\mbox{-a.s.}
\eea
By the Blumenthal 0-1 law, $Y^{0,t,\o}_t$ is a constant and we thus define
\bea
\label{u=Y}
u^0(t,\o) := Y^{0,t,\o}_t.
\eea

\begin{Theorem}
\label{thm-existence}
Let $\xi\in\mbox{UCB}(\Omega)$ be an $\Fc_T-$measurable r.v. Then, under Assumption \ref{assum:gen}, $u^0$ is a viscosity solution of PPDE \reff{PPDE} with terminal condition $u^0_T=\xi$.
\end{Theorem}

\proof  Under our assumptions on the nonlinearity $F$, it follows from the boundedness and uniform continuity of $\xi$ that $u^0$ is uniformly continuous on $[0,T]\times\Omega$, see \cite{EKTZ}. We show that $u^0$ is a $\cP-$viscosity subsolution, the same line of argument allows to prove that $u^0$ is a $\cP-$viscosity subsolution. We proceed by contradiction, assuming that $u^0$ is not a viscosity subsolution. Then, there exist $(t,\o)\in [0,T)\times\O$ and $\f\in \underline\cA^{\Pc}u^0_t(\o)$ such that:
 \beaa
 2c \;:=\; 
  - \partial_t\f_0
 - \frac12\mbox{Tr}\big[\partial^2_{\omega\omega}\f_0\big]-F_t(\o,u^0_t(\o),\partial_\o\f_0)
 &>&
 0.
 \eeaa
Without loss of generality, we assume $u_t(\o)=\f_0$, and we set $(t,\o)=(0,0)$. Denote:
 \beaa
 &\phi_s:=  \partial_t\f_s + \frac12\mbox{Tr}\big[\partial^2_{\omega\omega}\f_s\big]+F_s(\f_s,\partial\f_s)~~\mbox{so that}~~\phi_0=-2c,&
 \\
 &\mbox{and}~~\tilde Y_s:=  \f_s,
 ~~\tilde Z_s:= \pa_\o\f_s, 
 ~~\d Y_s:=\tilde Y_s - Y_s,
 ~~\d Z_s :=\tilde Z_s- Z_s,
 ~~, s\in [0,T].&
 \eeaa
Applying It\^{o}'s formula, we have
 \beaa
 d (\d Y_s) 
 &=& 
 \big(\partial_t\f_s + \frac12\mbox{Tr}\big[\partial^2_{\omega\omega}\f_s\big]\big)ds+\tilde Z_s\cdot dB_s 
 +F_s(Y_s,Z_s)ds - Z_s\cdot dB_s
 \\
 &=&
 \big[\phi_s+F_s(Y_s,Z_s)-F_s(\tilde Y_s,\tilde Z_s)\big]ds
 +\d Z_s\cdot dB_s
 ,~~\P_0-\mbox{a.s.}
 \eeaa
Since $\d Y_0=0$, it follows from the $L_0-$Lipschitz property of $F$ that for all stopping time $\tau\in\Tc$:
 \beaa
 0
 &\ge&
 (\f-u)_\tau
 -\int_0^\tau\big(\phi_s-L_0|\delta Y_s|\big)ds
 +\int_0^\tau\big(\d Z_s \cdot dB_s+L_0|\delta Z_s|ds\big)
 ,~~\P_0-\mbox{a.s.}
 \eeaa
Define $\hc_\eps:=\eps\wedge\inf\{s>0:|B_s|\ge\eps\}\wedge\inf\{s>0:\phi_s-L_0|\d Y_s|\ge -c\}$, and notice that $\hc_\eps>0$, $\P_0-$a.s. since $\d Y_0=0$. Then, 
 \beaa
 0
 &\ge&
 (\f-u)_{\hc_\eps}
 +c\;\hc_\eps
 +\int_0^{\hc_\eps}\big(\d Z_s \cdot dB_s+L_0|\delta Z_s|ds\big)
 ,~~\P_0-\mbox{a.s.}
 \eeaa
By the Girsanov theorem, we may find a probability measure $\bar\P\in\Pc^{L_0}\subset\Pc$ such that $B+L_0\int_0^.\mbox{sgn}(Z_s)ds$ is a $\bar\P-$Brownian motion. 
 Then, it follows from the previous inequality that
$\E^{\bar\P}\big[(\f-u)_\hc\big]\le-c\;\E^{\bar\P}[\hc]<0$, contradicting the fact that $\f\in \underline\cA^{\Pc}u^0_0$.
\ep

\section{Wellposedness of fully nonlinear path-dependent PDEs}
\label{sect-comparisonnonlin}
\setcounter{equation}{0}

In this section, we outline the main results established in \cite{ETZ2} in the context of the fully nonlinear path-dependent PDE:
 \bea\label{PPDE-nonlin}
\cL u := -\partial_tu
 -G\big(.,u,\partial_\o u,\partial^2_{\o\o}u\big)
 =0
& \mbox{on}&
 [0,T)\times\O.
 \eea

\begin{Assumption}\label{assum-G}
The nonlinearity $G$ satisfies:
\\
{\rm (i)} The process $G(., y,z,\gamma)$ is  continuous, and $G(., 0, {\bf 0}, {\bf 0})$ is bounded. 
\\
{\rm (ii)} $G$ is elliptic, i.e. nondecreasing in $\gamma$. 
\\
{\rm (iii)} $G$ is  $L_0-$Lipschitz in $(y,z,\g)$, uniformly in $(t,\o)$.
\end{Assumption}

In the present fully nonlinear context, we shall consider Definition \ref{def:visco} of viscosity solutions with the sets of test processes $\underline{\Ac}$ and $\overline{\Ac}$ defined by means of
 \beaa
 \Pc:=\Pc_L~\mbox{for some}~L\ge L_0,
 &\mbox{and}&
 \Hc:=\{\hc=t\wedge\hc_O: t\in[0,T],0\in O\subset\R^d,~\mbox{bounded convex}\},
 \eeaa
where $\hc_O:=\inf\{t>0:B_t\not\in O\}$. Observe that, unlike the semilinear case, the set $\Pc_L$ of Definition \ref{def:PL} is a non-dominated family of probability measures.

Following the same line of argument as in the semilinear case, it is shown in \cite{ETZ1} that the following partial comparison results hold true.

\begin{Theorem}
Let $u,v\in\mbox{UCB}(\Omega)$ be viscosity subsolutions and supersolution, respectively of the equation \eqref{PPDE-nonlin}, with $u_T\le v_T$ on $\Omega$. Assume further that either one of them is in $C^{1,2}_\Pc(\Theta)$. Then, under Assumption \ref{assum-G}, $u\le v$ on $\Theta$.
\end{Theorem}

We next report the wellposedness result from \cite{ETZ2} which requires further conditions on the path-frozen PDE:
 \beaa
 {\rm (E)}^{t,\o}_{\eps}
 &g_{t,\o}\big(s,v(s,x),Dv(s,x),D^2v(s,x)\big)=0,&
 (s,x)\in Q_t^\eps:=[t,T]\times B_{\R^d}(\eps),
 \eeaa
where $B_{\R^d}(\eps)$ is the centered open ball of $\R^d$ with radius $\eps$. We denote the parabolic boundary of the domain $Q_t^\eps$ by $\partial Q_t^\eps:=[t,T)\times B_{\R^d}(\eps) \cup \{T\}\times\mbox{cl}\big[B_{\R^d}(\eps)\big]$.

\begin{Assumption}\label{assum:nonlinearPDE}
{\rm (i)} The process $G(\cdot, y,z,\g)$ is uniformly continuous, uniformly in $(y,z,\g)$;

{\rm (ii)} For all $\eps>0$, $(t,\o)\in\Theta$, and $h\in C^0\big(\partial Q_t^\eps\big)$, we have $\overline{v}=\underline{v}$, where:
 \bea\label{barv}
 \left.\begin{array}{lll}
 \overline v(s,x) 
 &:=& 
 \inf\Big\{w(s,x): w~\mbox{classical supersolution of }
                          {\rm (E)}^{t,\o}_{\eps}
                     ~\mbox{and}
                     ~w \ge h 
                     ~\mbox{on}~\partial Q^{\eps}_t\Big\},
 \\
 \underline v(s,x) 
 &:=& 
 \sup\Big\{w(s,x): w~\mbox{classical subsolution of }
                          {\rm (E)}^{t,\o}_{\eps}
                    ~\mbox{and}
                    ~w \le h 
                    ~\mbox{on}~\partial Q^{\eps}_t\Big\}.
 \end{array}\right.
 \eea 
\end{Assumption}

\begin{Remark}{\rm The following sufficient condition for the nonlinearity $g:=g_{t,\omega}$ to satisfy Assumption \ref{assum:nonlinearPDE} (ii) is reported from Proposition 8.2 of \cite{ETZ2}:
\\
(i) The nonlinearity $g\big(s,y,z,\gamma\big)$ is continuous in $s$, uniformly Lipschitz in $(y,z,\gamma)$, and non-decreasing in $\gamma$,
\\
(ii) The PDE ${\rm (E)}^{t,\o}_{\eps}$ satisfies existence and comparison in the sense of viscosity solutions within the class of bounded functions,
\\
(iii) Either one of the following conditions holds:
\\
$~~~~~$(iii-1) $g$ is convex in $(y, z, \g)$, $ g_\delta(.,\gamma):=\inf_{A\in\dbS_d, A\ge 0}\big\{ g(.,\gamma+A)- \mbox{Tr}[A]\big\}>-\infty$ for $0\le\d\le c_0$, for some $c_0>0$, and $g_\d\longrightarrow g$ as $\d\searrow 0$,
\\
$~~~~~${\rm (iii-2)} $g$ is convex in $\gamma$ and  uniformly elliptic:  for some constant $c_0>0$,
\beaa
g(., \gamma) - g(., \gamma') \ge c_0 \mbox{Tr}[\gamma-\gamma'] 
&\mbox{for any}& 
\gamma\ge \gamma'.
\eeaa
$~~~~~${\rm (iii-3)} $g$ is uniformly elliptic and $d\le 2$.
}
\end{Remark}

We finally formulate a technical condition on the final condition $\xi$. We shall denote $\overline{\omega}:=\max_{s\le t}\omega_s$, $\underline{\omega}:=\min_{s\le t}\omega_s$, and $\omega^t_s:=\omega_s-\omega_t$ for all $0\le t\le s\le T$.

\begin{Assumption}\label{assum:xi}
$\xi=g\big(\big(\omega_{t_i},\overline{\omega}_{t_i},\underline{\omega}_{t_i}\big)_{1\le i\le n},\omega\big)$ for some $0=t_0<\ldots<t_n=T$ and some function $g\in\mbox{UCB}(\R^{3dn}\times\Omega)$ satisfying for all $\theta\in\R^{3dn}$, $i<n$, and $\omega,\omega'\in\Omega$, there exist some$p>0$ and continuity modulus $\rho$ such that:
 \beaa
 |g(\theta,\omega)-g(\theta,\omega')| 
 \le
 \rho\big(\big\|(\omega-\omega')\big\|^p_{\L^p([t_i,t_{i+1}])}\big)
 &\mbox{whenever}&
  \omega_{\wedge t_i}=\omega'_{\wedge t_i}
 ~\mbox{and}~
 \omega^{t_{i+1}}={\omega'}^{t_{i+1}}.
 \eeaa
\end{Assumption}

We are now able for the wellposedness result proved in \cite{ETZ2}.

\begin{Theorem}\label{thm-wellposedness}
Let Assumptions \ref{assum-G}, \ref{assum:nonlinearPDE}, \ref{assum:xi} hold true. 
\\
{\rm (i)}  Let $u,v\in \mbox{UCB}(\Omega)$ be $\Pc-$viscosity subsolution and supersolution, respectively, of PPDE (\ref{PPDE-nonlin}) with $u_T \le \xi\le v_T$. Then  $u\le v$ on $\Theta$.
\\
{\rm (ii)} The PPDE \reff{PPDE-nonlin} with terminal condition $\xi$ has a unique viscosity solution $u\in \mbox{UCB}(\Theta)$.
\end{Theorem}

\section{Stability of viscosity solutions of path-dependent PDEs}
\label{sect-stability}
\setcounter{equation}{0}

\subsection{Stability}
We shall establish the stability in the context of fully nonlinear PPDE, and thus we use the setting in Section \ref{sect-comparisonnonlin}. We first report the fully nonlinear analogue of Lemmas \ref{lem:regsubmart} and \ref{lem:maximumpoint}.
\begin{Lemma}\label{lem:maximumpoint-nonlinear}
Let $u\in  \mbox{UCB}(\Theta)$, $(t,\omega)\in\Theta$, and $\hc\in\Hc$, be such that $u_t(\omega) > \ol\Ec^{\Pc}\big[u^{t,\omega}_\hc\big]$. Then, 
 \b*
 0\in\underline{\Ac}^{\Pc}u_{t+t^*}(\o\otimes_t\o^*)
 &\mbox{for some}&
 (t^*,\o^*)~~\mbox{with the localization}~~\hc^*:=\hc^{t^*,\o^*}-t^*\in\Hc. 
 \e*
\end{Lemma} 
We remark that in this case $\cP$ has no dominating measure, and consequently the Dominated Convergence Theorem fails under $\overline\cE^\cP$. The proof of Lemma \ref{lem:maximumpoint-nonlinear} relies on the theory of optimal stopping under nondominated nonlinear expectation. The Snell envelop approach in this context is rather technical, and makes crucially use of the regularity of $X$ and the special structure of $\hc$, see \cite{ETZ0}. 




We now present the stability result. Fix $\cP=\cP_L$ and simplify the notations: $\overline\cE := \overline\cE^\cP$, $\underline\cE := \underline \cE^\cP$.
\begin{Theorem}
\label{thm-stability}  Let  $G$, $G^\eps$ satisfy Assumption \ref{assum-G} with a common $L_0 \le L$, and $u, u^\eps\in \mbox{UCB}(\Omega)$, for each $\eps>0$. Assume 

{\rm (i)} for each $\eps>0$, $u^\eps$ is a viscosity subsolution (resp. supersolution) of PPDE \reff{PPDE-nonlin} with generator $G^\eps$;

{\rm (ii)} as $\eps\to 0$, $(G^\eps, u^\eps)$ converge to $(G, u)$  locally uniformly.

\no Then $u$ is a viscosity subsolution (resp. supersolution) of PPDE \reff{PPDE-nonlin} with generator $G$.

\end{Theorem}
\proof Without loss of generality we shall only prove the viscosity subsolution property at $(0, {\bf 0})$. Let $\f\in {\ul\cA}^{\cP}\!u(0,{\bf 0})$  with corresponding $\hc \in \Hc$, $\d_0>0$ be a constant  such that $\hc_{\d_0} \le \hc$ and $\lim_{\eps\to 0} \rho(\eps,\d_0)=0$, where $\rho(\eps,\d_0)$ is the bound of $|G^\eps-G|+|u^\eps-u|$ on the $\d_0$-neighborhood of $(0,0,  y_0, z_0, \g_0):=(0,0,  u_0, \pa_\o \f_0, \pa_{\o\o}^2\f_0)$.

Now for $0 <\d \le \d_0$, denote $\f_\d(t,\o) := \f(t,\o) + \d t$. One can easily show that  $ \underline \cE_0[\hc_\d]>0$, see  \cite{ETZ1}. Then we have
\beaa
(\f_\d-u)_0&=& (\f-u)_0 \le \underline \cE\Big[(\f-u)_{\hc_\d}\Big]  = \underline \cE\Big[(\f_\d -u)_{\hc_\d} - \d \hc_\d\Big] \\
&\le&\underline \cE\Big[(\f_\d -u)_{\hc_\d} \Big] - \d  \underline \cE[\hc_\d] <  \underline \cE\Big[(\f_\d -u)_{\hc_\d}\Big].
\eeaa 
By the local uniform convergence of $G^\eps$ and $u^\eps$, there exists $\eps_\d>0$ small enough such that
\bea
\label{stabilityest1}
(\f_\d-u^\eps)_0  <  \underline \cE\Big[(\f_\d -u^\eps)_{\hc_\d}\Big],\quad \forall \eps\le \eps_\d.
\eea
By Lemma \ref{lem:maximumpoint-nonlinear}, we may find a point $(t^*,\o^*)$ such that
$$0\in\ul\cA^{\cP}(u^\epsilon-\f_\d)_{t^*}(\o^*)~\mbox{with the localization}~\hc^*:=\hc_\d^{t^*,\o^*}-t^*\in\Hc.$$
It is straightforward to check that $ \f_\d^{t^*, \o^*} \in {\underline\cA}^{\cP} u^\eps(t^*, \o^*)$. Since $u^\eps$ is a viscosity subsolution of PPDE \reff{PPDE-nonlin} with generator $G^\eps$, we have 
\bea
\label{stabilityest2}
0 \ge \Big[-\pa_t \f_\d - G^\eps(\cdot, u^\eps, \pa_\o \f_\d, \pa_{\o\o}^2 \f_\d)\Big](t^*,  \o^*)
=  \Big[-\pa_t \f - \d - G^\eps(\cdot, u^\eps, \pa_\o \f, \pa_{\o\o}^2 \f)\Big](t^*, \o^*).
\eea

Note that $t^* < \hc_\d(\o^*)$, then $|u^\eps- u|(t^*, \o^*)  \le \rho(\eps, \d)\le \rho(\eps, \d_0)$.  By local uniform convergence, we may set $\d$ small enough and then $\eps$ small enough so that $(\cdot, u^\eps, \pa_\o \f, \pa_{\o\o}^2 \f)(t^*, \o^*)$ is in the   $\d_0$-neighborhood of $(0,0,  y_0, z_0, \g_0)$. Thus, \reff{stabilityest2}  and Assumption \ref{assum-G} lead to
\beaa
0 &\ge&  \Big[-\pa_t \f  - G(\cdot, u^\eps, \pa_\o \f, \pa_{\o\o}^2 \f)\Big](t^*, \o^*) - \d - \rho(\eps, \d_0)\\
&\ge&  \Big[-\pa_t \f  - G(\cdot, u, \pa_\o \f, \pa_{\o\o}^2 \f)\Big](t^*, \o^*) - \d - \rho(\eps, \d_0) - C \rho(\eps, \d)\\
&\ge& \cL\f_0 - \sup_{(t, \o): t < \hc_\d(\o)} \Big|G(\cdot, u, \pa_\o \f, \pa_{\o\o}^2 \f)(t, \o) - G(\cdot, u, \pa_\o \f, \pa_{\o\o}^2 \f)(0,0)\Big|   -  \d -C \rho(\eps, \d_0).
\eeaa
Now by first sending $\eps\to 0$ and then $\d\to 0$ we obtain $\cL\f_0  \le 0$. Since  $\f\in {\ul\cA}^\cP \!u(0,{\bf 0})$ is arbitrary, we see that $u$ is a viscosity subsolution of PPDE \reff{PPDE-nonlin} with generator $G$ at $(0, {\bf 0})$ and thus complete the proof.
 \ep
 
 \subsection{Monotone scheme for PPDEs}
As an important application of the above stability result (in spirit),  in this subsection we  study discretization schemes for PPDEs. For any $(t,\o)\in \Theta$ and $h \in (0, T-t)$, we denote $\cF^{t,\o}_{t+h}:=\cF_{t+h}\cap\{B_{t\wedge\cdot}=\o_{t\wedge\cdot}\}$. Let $\T^{t,\o}_h$ be an operator on $\L^0(\cF^{t,\o}_{t+h})$.  For $n\ge 1$, denote $h:= {T\over n}$, $t_i := i h$, $i=0,1,\cdots, n$, and define:
\bea
\label{PPDE-uh}
u^h(t_n, \o) := \xi(\o),~~ u^h(t, \o) := \T^{t,\o}_{t_i - t}\big[u^h(t_{i},\cdot)\big],~ t \in [t_{i-1}, t_i),~ i=n,\cdots, 1.
\eea
where we abuse the notation that: 
\beaa
\T^{t,\o}_h[\f_s] := \T^{t,\o}_h[\f^{t,\o}_{s-t}]~\mbox{for process} ~\f.
\eeaa

 \begin{Assumption}
\label{assum-num1} Assumption \ref{assum-G} holds, and

{\rm (i)}   $\xi:\O\to\dbR$ is bounded and uniformly continuous.

{\rm (ii)}  Comparison principle for PPDE \reff{PPDE}  holds in the class of bounded viscosity solutions.
\end{Assumption}

\begin{Assumption}
\label{assum-num2}
The descritization operator  $\T^{t,\o}_h$ satisfies the following  conditions:
 
{\rm (i)} Consistency:   for any $(t,\o)\in \Theta$ and $\f\in C^{1,2}(\Theta)$,
\beaa
\lim_{(t', \o', h, c) \to (t,{\bf 0},0,0)}{[c+\f](t',\o')- \T^{t', \o\otimes_t \o'}_h\big[[c+\f](t'+h, \cdot)\big]\over h} = \cL\f(t,\o).
\eeaa
where $(t', \o')\in \Theta$, $h\in (0, T-t)$, $c\in \dbR$.

{\rm (ii)} Monotonicity: for some constant $L\ge L_0$ and any $\f, \psi \in \mbox{UCB}(\cF^t_{t+h})$, 
\beaa
\ol \cE^{\cP}[ (\f -\psi)^{t,\o}] \le 0\quad\mbox{implies}\quad \T^{t,\o}_h[\f] \le \T^{t,\o}_h[\psi].
\eeaa

{\rm (iii)} Stability:  $u^h$ is uniformly bounded and uniformly continuous in $\o$, uniformly on $h$. Moreover, there exists a modulus of continuity function $\rho$, independent of $h$,  such that
\beaa
|u^h(t,\o) - u^h(t', \o_{\cdot\wedge t})|\le \rho\Big((t'-t)\vee h\Big),\quad\mbox{for any $t<t'$ and any $\o\in \O$}.
\eeaa
\end{Assumption}

We now report the result from \cite{ZZ}, which extends the seminal work Barles and Souganidis \cite{BS} to our path dependent case.

\begin{Theorem}
\label{thm-monotone}
Let Assumptions \ref{assum-num1} and \ref{assum-num2} hold.  Then PPDE \reff{PPDE-nonlin} with terminal condition $u(T,\cdot) = \xi$ has a unique bounded viscosity solution $u$, and $u_h$ converges to  $u$ locally uniformly as $h\to 0$.
 \end{Theorem}
 
\proof By the stability, $u^h$ is bounded. Define
\bea
\label{ulu}
\ul u(t,\o) :=\liminf_{h\rightarrow 0} u^h(t,\o),\quad \ol u(t,\o) :=\limsup_{h\rightarrow 0} u^h(t,\o).
\eea
Clearly $\ul u (T,\o) = \xi(\o) = \ol u(T,\o)$, $\ul u\le \ol u$, and $\ul u, \ol u$ are bounded and uniformly continuous. We shall show that $\ul u$ (resp. $\ol u$) is a viscosity supersolution (resp. subsolution) of PPDE \reff{PPDE-nonlin}. Then by the comparison principle we see that $\ol u \le \ul u$ and thus $u:=\ol u = \ul u$ is the unique viscosity solution of PPDE \reff{PPDE-nonlin}.   The convergence of $u^h$ is obvious now, which, together with the uniform regularity of $u^h$ and $u$, implies further the locally uniform convergence.

Without loss of generality, we shall only prove  by contradiction that $\ul u$  satisfies the viscosity supersolution  property  at $(0, {\bf 0})$. Assume not, then there exists $\f^0 \in \ol \cA^{\cP} \ul u(0,{\bf 0})$ with corresponding $\hc \in \Hc$ such that
$
-c_0:=\cL \f^0(0,0)<0.
$
Denote 
\bea
\label{f0}
\f(t,\o) := \f^0(t,\o) - {c_0\over 2}t.
\eea
Then
\bea
\label{cLf0}
\cL \f(0,0) = -{c_0\over 2}<0.
\eea

Denote $X^0:=\f-\ul u$, $X^h := \f-u^h$, $\hc_\eps := \hc^0_\eps\wedge \eps^5:= \inf\{t: |B_t|\ge \eps\}  \wedge \eps^5$, and $c_\eps := {1\over 3}c_0\eps^5$. Note that $\hc_\eps  \le \hc$  for $\eps$ small enough, and by \cite{ETZ1} (2.8),
\bea
\label{cheeta}
\sup_{\dbP\in \cP_{L}} \dbP(\hc_\eps \neq \eps^5) = \sup_{\dbP\in \cP_{L}}  \dbP(\hc_\eps^0 < \eps^5) \le CL^4 \eps^{-4} \eps^{10} \le C\eps c_\eps .
\eea
Then
\beaa
\ol\cE[\eps^5 - \hc_\eps] \le \ol\cE\big[\eps^5 \1_{\{ \hc_\eps\neq \eps^5\}}\big] \le C\eps c_\eps.
\eeaa
  Thus, for $\eps$ small, it follows from $\f^0\in \ol\cA^{L} \ul u(0,{\bf 0})$ that
\bea
\label{ce}
X^0_0 - \ol \cE[X^0_{\hc_\eps}] &=& [\f^0-\ul u]_0 -  \ol \cE\Big[(\f^0-\ul u)_{\hc_\eps} - {c_0\over 2} {\hc_\eps}\Big] \nonumber\\
&\ge&  \ol \cE\Big[(\f^0-\ul u)_{\hc_\eps} \Big] - \ol \cE\Big[(\f^0-\ul u)_{\hc_\eps} - {c_0\over 2} {\hc_\eps}\Big] \\
&\ge&  \ul \cE\Big[{c_0\over 2} {\hc_\eps}\Big]  = {c_0\eps^5\over 2} - {c_0\over 2} \ol\cE[\eps^5-{\hc_\eps}] \ge {3c_\eps\over 2} - C\eps c_\eps \ge c_\eps >0.\nonumber
\eea
Let $h_k\downarrow 0$ be a sequence such that
\bea
\label{hk}
\lim_{k\to\infty} u^{h_k}_0 = \ul u_0,
\eea
and simplify the notations: $u^k := u^{h_k}$, $X^k := X^{h_k}$.  Then \reff{ce} leads to
\beaa
c_\eps \le [\f_0-\liminf_{h\to 0} u^h_0] -  \ol \cE\Big[\f_{\hc_\eps} - \liminf_{h\to 0} u^h_{\hc_\eps} \Big]  \le [\f_0-\lim_{k\to \infty} u^{k}_0] -  \ol \cE\Big[\f_{\hc_\eps} - \liminf_{k\to\infty} u^{k}_{\hc_\eps} \Big].
\eeaa
Note that $X^k$ is uniformly bounded. Then by \reff{cheeta} we have
\beaa
\ol\cE\Big[|X^k_{\hc_\eps}-X^k_{\eps^5}|\Big] \le C\eps c_\eps.
\eeaa
Since $u^h$ is uniformly continuous, applying the monotone convergence theorem under nonlinear expectation $\ol\cE$, see e.g. \cite{ETZ0} Proposition 2.5, we have
\beaa
c_\eps  & \le& \lim_{k\to \infty} [\f_0- u^{k}_0] -  \ol \cE\Big[ \limsup_{k\to\infty} [\f_{\hc_\eps} - u^{k}_{\hc_\eps}] \Big] \\
&\le& \lim_{k\to \infty} X^k_0 -  \ol \cE\Big[ \limsup_{k\to\infty} X^k_{\eps^5} \Big] + C\eps c_\eps =  \lim_{k\to \infty} X^k_0 -  \ol \cE\Big[ \lim_{m\to\infty}\sup_{k\ge m} X^k_{\eps^5} \Big] + C\eps c_\eps\\
&=& \lim_{k\to \infty} X^k_0 -   \lim_{m\to\infty}\ol \cE\Big[\sup_{k\ge m} X^{k}_{\eps^5}\Big]+ C\eps c_\eps \le  \lim_{k\to \infty} X^k_0 -   \limsup_{k\to\infty}\ol \cE\Big[X^{k}_{\eps^5}\Big]+ C\eps c_\eps \\
&\le& \lim_{k\to \infty} X^{k}_0 -  \limsup_{k\to\infty} \ol \cE\Big[X^{k}_{\hc_\eps} \Big]+ C\eps c_\eps=  \liminf_{k\to\infty}\Big[ X^{k}_0-\ol \cE\big[X^{k}_{\hc_\eps} \big]\Big] + C\eps c_\eps.
\eeaa
Then, for all $\eps$ small enough and  $k$ large enough,
\bea
\label{Xk0}
X^{k}_0-\ol \cE\big[X^{k}_{\hc_\eps} \big] \ge {c_\eps\over 2}.
\eea

Now, applying Lemma \ref{lem:maximumpoint-nonlinear}, we obtain that
\b*
 0\in\underline{\Ac}^{\Pc}X^k_{t^k_*}(\o^*)
 &\mbox{for some}&
 (t^k_*,\o^*)~~\mbox{with the localization}~~\hc^k_\eps:=\hc_\eps^{t_*^k, \o^*}-t^k_*. 
\e*
Moreover, in this case, we may prove that 
\be\label{tk*}
\sup_{\dbP\in\Pc}\dbP\big[\hc^k_\epsilon\leq \d\big]\leq C \delta^2
\ee 
for all $\delta\leq h_k$ (see \cite{ZZ}).
Let $\{t^k_i, i=0,\cdots, n_k\}$ denote the time partition corresponding to $h_k$, and assume $t^k_{i-1} \le t^k_* < t_i^k$. Note that 
\beaa
X^k_{t^k_*} (\o^k) = Y^k_{t^k_*} (\o^k)  \ge \overline\cE\Big[(X^k)^{t^k_*,\o^k}_{\t\wedge \hc_\eps^k}\Big],\quad\forall \t \in \cT.
\eeaa
Set $\d_k := t^k_i - t^k_* \le h_k$ and $\t:= \d_k$. Combine the above inequality and \reff{tk*} we have
\beaa
[\f - u^k](t^k_*,\o^k) \ge \overline\cE\Big[(\f-u^k)^{t^k_*,\o^k}_{\d_k\wedge \hc_\eps^k}\Big] \ge \overline\cE\Big[(\f-u^k)^{t^k_*,\o^k}_{\d_k}\Big] - C \d_k^2.
\eeaa
This implies
\beaa
 \overline\cE\Big[\Big(\f^{t^k_*,\o^k}_{\d_k} -  [\f - u^k](t^k_*,\o^k) - C \d_k^2\Big) - (u^k)^{t^k_*,\o^k}_{\d_k}\Big] \le 0.
 \eeaa
 By the monotonicity condition (Assumption \ref{assum-num2} (ii)) we have
 \bea
 \label{uktk}
 u^k(t^k_*,\o^k) = \T^{t^k_*,\o^k}_{\d_k}[u^k_{t^k_i}]  \le  \T^{t^k_*,\o^k}_{\d_k}\Big[\f_{t^k_i} -  [\f - u^k](t^k_*,\o^k) - C \d_k^2\Big].
 \eea
 
We next use the consistency condition (Assumption \ref{assum-num2} (i)). For $(t,\o) = (0, {\bf 0})$, set
\beaa
t':= t^k_*,\quad \o' := \o^k,\quad h:= \d_k,\quad c := -  [\f - u^k](t^k_*,\o^k) - C \d_k^2.
\eeaa
By first sending $k\to \infty$ and then $\eps\to 0$, we see that
\beaa
d((t^k_*, \o^k), (0,{\bf 0})) \le \hc_\eps + \sup_{0\le t\le \hc_\eps} |\o^k_t|\le 2\eps \to 0,\quad h \le h_k\to 0,
\eeaa
which, together with \reff{f0}, \reff{hk}, and the uniform continuity of $\f$ and $u^k$, implies
\beaa
|c| \le \Big| [\f - u^k](t^k_*,\o^k) -  [\f - u^k](0, {\bf 0})\Big| + | u^k_0 - \ul u_0| + C\d_k^2 \to 0.
\eeaa
Then, by the consistency condition, we obtain from \reff{uktk} that
\beaa
0 &\le& {u^k(t^k_*,\o^k)  - \T^{t^k_*,\o^k}_{\d_k}\Big[\f_{t^k_i} -  [\f - u^k](t^k_*,\o^k) - C \d_k^2\Big]\over \d_k} \\
&=&{ [c + \f](t^k_*, \o^k) - \T^{t^k_*,\o^k}_{\d_k}\Big[ [c+\f]_{t^k_i}\Big]\over \d_k} + C\d_k\to  \cL \f(0, {\bf 0}).
\eeaa
This contradicts with \reff{cLf0}. 

\section{Optimal stopping under dominated nonlinear expectation} 
\label{sect:optimalstop}

The objective of this section is to provide a self-contained proof of Theorem \ref{thm:optimalstop}. We follow the setting in Section \ref{sect-comparisonslin}. In particular, the family $\Pc$ of equivalent probability measures is defined as in (\ref{Psemilinear}).

We emphasize that the main results of this section are available in the literature in reflected backward stochastic differential equations, see \cite{ElkarouiKPPQ,Hamadene,PengXu}. We collect them here for the convenience of the readers who might not be familiar with this literature. Our presentation in Subsections \ref{subsect:DMSnell} and \ref{subsect:F*Snell} is inspired from El Karoui \cite{ElKaroui} and Appendix D of Karatzas and Shreve \cite{KaratzasShreve}, which are focused on the standard optimal stopping under linear expectation. 

\subsection{Preliminaries}

For ease of notation, we simply write $\ol\cE:=\ol{\cE}^\cP$. We start by the dominated convergence theorem under $\ol{\cE}$ which holds by the fact that $\cP$ is dominated by $\dbP_0$.

\begin{Lemma}\label{ConvergenceThm}
Let $X_n$ be a sequence of random variables. Assume that $X_n^{1+\alpha}$ are uniformly integrable under probability $\P_0$ and $X_n\rightarrow X$ $\dbP_0$-a.s. Then, we have $\ol\cE[|X_n-X|]\rightarrow 0$.
\end{Lemma}
\proof
For any $\dbP_\l\in\cP$, we have
\beaa
\dbE^{\dbP_\l}[|X_n-X|] &=& \dbE^{\dbP_0}\Big[e^{\int_0^T \l_s dB_s-\frac12\int_0^T |\l_s|^2ds}|X_n-X|\Big]\\
&\leq & \big(\dbE^{\dbP_0}[e^{\int_0^T q\l_s dB_s-\int_0^T \frac{q}{2}|\l_s|^2ds}]\big)^\frac{1}{q}\big(\dbE^{\dbP_0}[|X_n-X|^p] \big)^\frac{1}{p},
\eeaa
where $p=1+\alpha$ and $\frac{1}{p}+\frac{1}{q}=1$. Since $\l$ is bounded, we know that $\dbE^{\dbP_0}\big[e^{\int_0^T q\l_s dB_s-\int_0^T \frac{q}{2}|\l_s|^2ds}\big]$ is bounded. Then, by the convergence theorem, we obtain that $\dbE^{\dbP_0}[|X_n-X|^{1+\alpha}]\rightarrow 0$. The proof is complete.
\ep

\begin{Lemma}
If $X\geq 0$ $\dbP_0$-a.s. and $\ul\cE[X]=0$, then $X=0$ $\dbP_0$-a.s..
\end{Lemma}
\proof
Since $\ul\cE[X]=0$, for any $\epsilon>0$ there exists $\dbP^\epsilon\in\cP$ such that $\dbE^{\dbP^\epsilon}[X]<\epsilon$. Also, by Cauchy-Schwarz inequality, we have the estimate:
$$\dbE^{\dbP_0}[X^\frac12] = \dbE^{\dbP^\epsilon}\Big[e^{-\int_0^T \l^\epsilon dB_s-\frac12\int_0^T |\l^\epsilon|^2ds}X^\frac12\Big]\leq C\dbE^{\dbP^\epsilon}[X]^\frac12<C\epsilon^\frac12.$$
Since $\epsilon$ is arbitrary, we get $\dbE^{\dbP_0}[X^\frac12]=0$. So, we conclude that $X=0$ $\dbP_0$-a.s..
\ep

\vspace{5mm}

Finally, we state the following lemma, which is a direct consequence of Proposition 3.1. in El Karoui, Peng and Quenez \cite{ElkarouiPengQuenez}.

\begin{Lemma}\label{lem:optimalprob}
Let $\xi\in\L^2(\P_0)$, and $v_t:=\esup_{\P\in\cP}\E^\P[\xi|\cF_t]$. Then, $v_t=\E^{\bar\P}[\xi|\cF_t]$ $\P_0$-a.s. for all $t\in [0,T]$ for some $\bar\P\in\cP$.
\end{Lemma}

\subsection{RCLL version of the $\dbF^*-$Snell envelop}
\label{subsect:F*Snell}

Throughout this section, we consider a process $X:[0,T]\times\Omega\longrightarrow\dbR$ satisfying the following condition.

\begin{Assumption}\label{assum X}
The process $X$ is piecewise pathwise continuous $\dbF$-adapted on $[0,T]$, and $\sup_{t\in[0,T]}|X_t|\in\L^2(\cP)$, i.e. 
 \b*
 \E^\P\Big[\sup_{t\in[0,T]}|X_t|^2\big]
 \;<\;\infty,
 &\mbox{for all}&
 \P\in\cP.
 \e* 
\end{Assumption}

Our starting point is the classical Snell envelop process:
 \beaa
 Y_t
 &:=&
 \esup_{\t\in\cT_*^t,\dbP\in\cP}\;\dbE^{\dbP}[X_\t|\cF_t]
 ,~~t\in[0,T].
 \eeaa
Clearly, $Y_t$ is $\cF^*_t$-measurable for all $t\in[0,T]$.

\begin{Lemma}\label{lattice}
For any $t\in[0,T)$, $\big\{\dbE^{\dbP}[X_\t|\cF_t];(\t,\dbP)\in\cT_*^t\times\cP\big\}$ satisfies the lattice property.
\end{Lemma}

\proof
Let $\t_1,\ \t_2\in\cT_*^t$ and $\dbP_1,\ \dbP_2\in\cP$. Let $A:=\{\dbE^{\dbP_1}[X_{\t_1}|\cF_t]\geq \dbE^{\dbP_2}[X_{\t_2}|\cF_t]\}$, and define
 \beaa
 \bar{\t}:=\t_1 1_A+\t_2 1_{A^c}
 &\mbox{and}&
 \bar\dbP(D)
 :=
 \dbE^{\dbP_1}\Big[\dbE^{\dbP_1}[1_{A\cap D}|\cF_t]+\dbE^{\dbP_2}[1_{A^c\cap D}|\cF_t]\Big],
 ~~D\in\Fc_T.
 \eeaa
Clearly, $\bar\t\in\cT_*^t$, $\bar\dbP\in\cP$, and we immediately verify that
$$\dbE^{\bar\dbP}[X_{\bar\t}|\cF_t]\geq\max\{\dbE^{\dbP_1}[X_{\t_1}|\cF_t],\dbE^{\dbP_2}[X_{\t_2}|\cF_t]\},\ \dbP_0\text{-a.s.}$$
\ep

\vspace{5mm}

We next introduce the concatenation $\dbP_1\otimes_t\dbP_2$ of two probability measures $\dbP_1,\dbP_2\in\cP$ by:
 \beaa
 (\dbP_1\otimes_t\dbP_2)(D)
 :=
 \dbE^{\dbP_1}\Big[\dbE^{\dbP_2}[1_D|\cF_t]\Big]
 &\mbox{for all}&
 D\in \cF_T,
 \eeaa
and we observe that $\dbP_1\otimes_t\dbP_2\in\cP$.

\begin{Lemma}\label{lem:Y supmartingale}
$Y$ is an $\ol\cE$-supermartingale with $\sup_{t\in[0,T]}\E^{\P_0}\big[|Y_t|^2\big]<\infty$ and $\ol{\cE}[Y_t]=\sup_{\t\in\cT_*^t}\ol\cE[X_\t]$ for all $t\in[0,T]$.
\end{Lemma}

\proof
Denote $|X|^*_T:=\sup_{t\in[0,T]}|X_t|$. By the definition of $Y$, we have
 \beaa
 \sup_{t\in[0,T]}\E^{\P_0}\big[|Y_t|^2\big]
 \leq 
 \E^{\P_0}\big[\esup_{\P\in\cP}\E^\P\big[(|X|^*_T)^2 \big|\cF_t\big]\big]
 \leq 
 \sup_{\P\in\cP}\E^\P\big[(|X|^*_T)^2\big]<\infty.
 \eeaa
For arbitrary $\dbP\in\cP$ and $s\leq t$, it follows from Lemma \ref{lattice} and the property of the $\esup$ that:
$$
\dbE^\dbP[Y_t|\cF_s] 
=\!\!
\esup_{\t\in\cT_*^t;\dbP'\in\cP}\dbE^{\dbP\otimes\dbP'}[X_\t|\cF_s] 
\le\!\! 
\esup_{\t\in\cT_*^t;\dbP'\in\cP}\dbE^{\dbP'}[X_\t|\cF_s]
\le\!\!
\esup_{\t\in\cT_*^s;\dbP'\in\cP}\dbE^{\dbP'}[X_\t|\cF_s]
=
Y_s,
~\dbP_0-\mbox{a.s.}
$$
which proves that $Y$ is an $\ol\cE-$supermartingale. 

We finally prove the last claim. For all $\t\in\cT_*^t$ and $\dbP\in\cP$, we have $Y_t\geq\dbE^\dbP[X_\t|\cF_t],\ \dbP_0\text{-a.s.}$ Hence, we obtain for any $\t\in\cT_*^t$ and  $\dbP,\ \dbP^{'}\in\cP$ that
$\ol\cE[Y_t]\geq\dbE^{\dbP^{'}}[Y_t]\geq\dbE^{\dbP^{'}\otimes\dbP}[X_\t]$, and therefore $\ol\cE[Y_t]\geq\sup_{\t\in\cT_*^t}\ol\cE[X_\t]$. On the other hand, it follows from Lemma \ref{lattice} that:
 \beaa
 \dbE^{\dbP}[Y_t]
 =
 \sup_{\t\in\cT_*^t,\dbP^{'}\in\cP}\dbE^{\dbP\otimes\dbP^{'}}[X_\t]\leq \sup_{\t\in\cT_*^t}\ol\cE[X_\t]
 &\mbox{for all}&
 \dbP\in\cP.
 \eeaa
\ep

\begin{Proposition}\label{DPP}{\rm (Dynamic programming principle)} For all $t\in[0,T)$ and $\th\in\cT_*^t$:
 \b*
 Y_t
 \;=\;
 \esup_{\t\in\cT_*^t,~\dbP\in\cP}
 \dbE^\dbP\Big[X_\t 1_{\{\t<\th\}}+Y_\th 1_{\{\t\geq \th\}}\big|\cF_t\Big],
 &\dbP_0\text{-a.s}.&
 \e*
\end{Proposition}

\proof
Since $X\le Y$, we have for all $\th\in\cT_*^t$
$$
Y_t 
\le 
\esup_{\t\in\cT_*^t,\dbP\in\cP}\dbE^\dbP[X_\t 1_{\{\t<\th\}}+Y_\t 1_{\{\t\geq \th\}}|\cF_t]
\le 
\esup_{\t\in\cT_*^t,\dbP\in\cP}\dbE^\dbP[X_\t 1_{\{\t<\th\}}+Y_\th 1_{\{\t\geq \th\}}|\cF_t], 
~\dbP_0-\text{a.s.}
$$
where the last inequality is due to the $\ol\cE$-supermartingale property of $Y$ of Lemma \ref{lem:Y supmartingale}. On the other hand, since $Y$ is $\ol\cE$-supermartingale, we have for all $\t\in\cT_*^t$ and $\dbP\in\cP$:
 $$
 Y_t 
 \ge
 \dbE^{\dbP}[Y_{\th\we \t}|\cF_t]
 =
 \dbE^\dbP[Y_{\th}1_{\{\t\geq\th\}}+Y_\t 1_{\{\t<\th\}}|\cF_t]
 \ge
 \dbE^\dbP[Y_{\th}1_{\{\t\geq\th\}}+X_\t 1_{\{\t<\th\}}|\cF_t],
 ~\dbP_0\text{-a.s}.
 $$
The proof is completed by taking $\esup$ over $\t\in\cT_*^t$ and $\dbP\in\cP$.
\ep

\begin{Lemma}\label{lem: RCLL version}
$Y$ has a $\dbP_0$-a.s. RCLL $\dbF^*-$adapted version. Moreover, there exists $\bar\P\in\cP$ such that $\E^{\bar \P}[\sup_{t\in[0,T]}|Y_t|^2]<\infty$.
\end{Lemma}

\proof
{\it Step 1.} Let $\{t_n\}\subset[0,T]$ be such that $t_n\searrow t$. By Lemma \ref{lem:Y supmartingale}, we know that $\ol\cE[Y_{t_n}]=\sup_{\t\in\cT_*^{t_n}}\ol\cE[X_\t]\leq \sup_{\t\in\cT_*^{t}}\ol\cE[X_\t] \leq \ol\cE[Y_t]$. On the other hand, for any $\t\in\cT_*^t$, denoting $\t_n:=\t\vee t_n$, it follows from the continuity of $X$ and the $\P_0-$uniform integrability of $\{X_{\t_n}^2,n\ge 1\}$ that $\ol\cE[X_\t]=\lim_{n\rightarrow\infty}\ol\cE[X_{\t_n}]\leq\liminf_{n\rightarrow\infty}\ol\cE[Y_{t_n}].$
Using again Lemma \ref{lem:Y supmartingale}, we obtain that $\ol\cE[Y_t]\leq \liminf_{n\rightarrow\infty}\ol\cE[Y_{t_n}]$. Hence,
 \beaa
 \ol\cE[Y_t] 
 &=& 
 \lim_{s\downarrow t}\ol\cE[Y_{s}].
 \eeaa
{\it Step 2.} It follows from Lemma \ref{lem:Y supmartingale} that $Y$ is a $\dbP_0$-supermartingale in the right-continuous filtration $\dbF^*$. By classical martingale theory, we know that for any $t\in[0,T)$,
$$Y_{t+}:=\lim_{s\downdownarrows t,s\in\dbQ}Y_s\text{ exists }\dbP_0\text{-a.s.}$$
Note that $Y_{t+}$ is $\cF^*_t$-measurable. Also, we have the properties that $\{Y_{t+}\}_t$ is RCLL and $Y_{t+}=\dbE[Y_{t+}|\cF^*_t]\leq Y_t$, $\dbP_0$-a.s. 

We now show that $Y_{t+}=Y_t$, $\dbP_0$-a.s. Suppose to the contrary that $\dbP_0[Y_{t+}<Y_t]>0$. Then, we have $\dbE^{\dbP_0}\big[\sqrt{Y_t-Y_{t+}}\big]>0$, implying that $\ul\cE[Y_t-Y_{t+}]>0$. Then, there exists $\d>0$ such that:  
 \bea\label{estm1 Yt-Yt+}
 \dbE^\dbP[Y_t-Y_{t+}]
 \;\ge\;
 \d \;>\; 0
 &\mbox{for all}&
 \dbP\in\cP.
 \eea
By the definition of $Y_{t+}$ and the fact that $\{Y_t^\frac{3}{2}\}$ are uniformly integrable (by Lemma \ref{lem:Y supmartingale}), we obtain by Lemma \ref{ConvergenceThm} that $\ol\cE[Y_t] = \lim_{s\downarrow t}\ol\cE[Y_{s}]=\ol\cE[Y_{t+}]$. This means that for all $\dbP\in\cP$ and $\epsilon>0$, there exists $\dbP^{'}\in\cP$  such that $\dbE^\dbP[Y_t]-\epsilon\leq \dbE^{\dbP^{'}}[Y_{t+}]$. Together with (\ref{estm1 Yt-Yt+}), this implies that $\dbE^\dbP[Y_t]-\epsilon\leq \dbE^{\dbP^{'}}[Y_t]-\d\leq\ol\cE[Y_t]-\d$, and therefore $\ol\cE[Y_t]-\epsilon\leq\ol\cE[Y_t]-\d$. By arbitrariness of $\epsilon>0$, this provides that $\ol\cE[Y_t]\leq \ol\cE[Y_t]-\d$, which is the required contradiction. So, we have proved that $Y_{t+}$ is an $\dbF^*-$adapted RCLL version of $Y_t$.
\\
{\it Step 3.} With $|X|^*_T:=\sup_{t\in[0,T]}|X_t|$, we have:
 \beaa
 \sup_{t\in[0,T]}|Y_t|^2\leq \sup_{t\in[0,T]}\esup_{\P\in\cP}\E^\P\big[\big(|X|_T^*\big)^2|\cF_t\big].
 \eeaa
By Lemma \ref{lem:optimalprob}, there exists $\bar\P\in\cP$ such that $\E^{\bar\P}[X^*|\cF_t]=\esup_{\P\in\cP}\E^\P[X^*|\cF_t]$ for all $t$, $\P_0$-a.s. Then, it follows from the Doob inequality that:
 \beaa
 \E^{\bar\P}\Big[\sup_{t\in[0,T]}|Y_t|^2\Big]
 &\le&
 \E^{\bar\P}\Big[\sup_{t\in[0,T]}\E^{\bar\P}\big[(|X|^*_T)^2\big|\cF_t\big]\Big]
 \;\le\; 
 4\E^{\bar\P}\big[(|X|^*_T)^2\big],
 \eeaa
which provides the desired result  by Assumption \ref{assum X}.
\ep

\subsection{Doob-Meyer decomposition of the RCLL $\dbF^*-$Snell envelop}
\label{subsect:DMSnell}

From now on, we consider $Y$ in its $\dbF^*-$adapted RCLL version of Lemma \ref{lem: RCLL version}. For a vector $x=(x_1,\ldots,x_d)\in\dbR^d$, we denote $|x|_1:=\sum_{i=1}^d|x_i|$.

\begin{Proposition}\label{prop DM}
There exist $H\in\mathbb{H}_{loc}$ and a non-decreasing previsible process $K$ such that
 \b*
 Y_t
 \;=\;
 Y_0+(H\cdot B)_t-L\int_0^t|H_s|_1ds-K_t,
 &t\in[0,T],&
 \dbP_0-\mbox{a.s.},
 \e*
with $\E^{\P_0}\big[\sup_{t\in[0,T]}|(H\cdot B)_t|\big]<\infty$.
\end{Proposition}

\proof
\no {\bf 1.}\quad By Lemma \ref{lem:Y supmartingale}, $Y$ is a $\dbP-$supermartingale, with Doob-Meyer decomposition, 
 \bea\label{P doob-meyer}
 Y
 \;=\;
 Y_0+M^\dbP-A^\dbP,
 &\dbP_0-\mbox{a.s.  for all}&
 \dbP\in\cP,
 \eea
for some $\dbP$-martingale $M^\dbP$ and some non-decreasing previsible process $A^\dbP$. By the martingale representation property, $M^{\dbP_0}=(H\cdot B)$, $\dbP_0$-a.s. for some $H\in\mathbb{H}_{loc}$. By the Girsanov theorem, the process $\tilde{M}^{\dbP_\l}:=M^{\dbP_0}-\int_0^.\l^T_sH_sds$ defines a $\dbP_\l$-local martingale. Then, it follows from the uniqueness of the Doob-Meyer decomposition that $\tilde{M}^{\dbP_\l}$ is a $\dbP_\l$-martingale, and
 \beaa
 \tilde{M}^{\dbP_\l}=M^{\dbP_\l}
 &\mbox{and}&
 \int_0^\cdot\l^T_sH_sds-A^{\dbP_0}=-A^{\dbP_\l},
 ~~\dbP_0-\mbox{a.s.}
 \eeaa
We next introduce the process $\lambda^*$ with $i-$th component proportional to the sign of the $i-$th component of $H$, so that $\lambda^*H=L|H|_1$. Note that $\P^{\l^*}\in\cP$. Then the required decomposition holds with $K:=A^{\dbP_{\l^*}}$.
\\
\no {\bf 2.}\quad
By It\^o's formula, we have
$$A^{\P}_tY_t-\int_0^tY_sdA^{\P}_s=\int_0^t A^{\P}_sdY_s=\int_0^t A^{\P}_sdM^{\P}_s-\frac12(A^{\P}_t)^2,\text{ for all }\P\in\cP.$$
Let $(\t_n)_n$ be a localizing sequence for the $\dbP-$local martingale $\int_0^\cdot A^{\P}_sdM^{\P}_s$. Then,
 $$
 \frac12 \E^{\P}[(A^{\P}_{\t_n})^2]
 \leq 
 2\E^{\P}[\sup_{t\in[0,T]}|Y_t|\cdot A^{\P}_{\t_n}]
 \leq 
 2\Big(\E^\P\big[\sup_{t\in[0,T]}|Y_t|^2\big]\E^\P\big[(A^\P_{\t_n})^2\big]\Big)^\frac12.
 $$
For $\P=\bar\P$ as in Lemma \ref{lem: RCLL version}, we conclude that $\E^{\bar\P}[(A^{\bar\P}_T)^2]<\infty$. Then, one may easily verify that $\E^{\bar\P}\big[\sup_{t\in[0,T]}|M^{\bar\P}|^2\big]<\infty$, and therefore $\E^{\bar\P}[\langle M^{\bar\P}\rangle_T]<\infty$  by the Burkholder-Davis-Gundy inequality. Then, it follows from the Cauchy-Schwartz inequality that 
 \b*
 \E^\P[\langle M^\P\rangle_T^\frac12]
 \;=\;
 \E^\P[\langle M^{\bar\P}\rangle_T^\frac12]
 \;\le\; 
 C(\E^{\bar\P}[\langle M^{\bar\P}\rangle_T])^\frac12<\infty,
 &\text{for all}&
 \P\in\cP,
\e*
and we conclude that $\E^{\P}\big[\sup_{t\in[0,T]}|M^\P_t|\big]<\infty$, by the Burkholder-Davis-Gundy inequality.
\ep

\vspace{5mm}

We next provide some further properties of the previsible nondecreasing process $K$, and we derive an optimal stopping rule.

\begin{Proposition}\label{prop:dK on X=Y}
The processes $Y$ and $K$ are pathwise continuous, $\int_0^T 1_{\{t:X_t<Y_t\}}dK_t=0$, $\dbP_0$-a.s. and the $\F^*$-previsible stopping time $\t^*:=\inf\{t:X_t=Y_t\}$ is an optimal stopping rule. 
\end{Proposition}

In order to prove this result, we introduce the stopping times
 \b*
 D^\eps_t
 \;:=\;
 \inf\{s\geq t:Y_s\leq X_s+\eps\}
 &\mbox{for all}&
 t\in[0,T),~~\eps>0.
 \e*
By the right-continuity of $Y$ and the continuity of $X$, it is clear that $D^\eps_t\in\cT_*^t$. The following two lemmas prepare for the proof of Proposition \ref{prop:dK on X=Y}.

\begin{Lemma}\label{lem:YDesp-Y}
For all $t\in[0,T)$, we have $\ol\cE\big[Y_{D^\eps_t}-Y_t\big]=0$.
\end{Lemma}

\proof Since $Y$ is $\cP$-supermartingale and $D^\eps_t\geq t$, we have $\ol\cE[Y_{D^\eps_t}-Y_t]\le 0$.
On the other hand,  by the dynamic programming principle of Proposition \ref{DPP}, we have
 \beaa
 Y_t
 &=&
 \esup_{\t\in\cT_*^t,\dbP\in\cP}
 \dbE^\dbP\Big[X_\t 1_{\{\t<D^\eps_t\}}+Y_{D^\eps_t}1_{\{\t\geq D^\eps_t\}}|\cF_t\Big],
 ~~\dbP_0-\mbox{a.s.}
 \eeaa
Here, we may prove the lattice property similar to Lemma \ref{lattice}, so that
 \beaa
 \dbE^\dbP[Y_t] 
 \;=\; 
 \sup_{\t\in\cT_*^t,\dbP'\in\cP}
 \dbE^{\dbP\otimes_t\dbP'}\big[X_\t \1_{\{\t<D^\eps_t\}}
                                                      +Y_{D^\eps_t}\1_{\{\t\geq D^\eps_t\}}
                                               \big]
 &\mbox{for all}&
 \dbP\in\cP.
 \eeaa
Then, there exists $(\t_n)_n\subset\cT_*^t$ such that
 \beaa
 \dbE^\dbP[Y_t] 
 &\le&
 \dbE^{\dbP\otimes_t\dbP_n}\big[X_{\t_n}\1_{\{\t_n<D^\eps_t\}}
                                                     +Y_{D^\eps_t}\1_{\{\t_n\geq D^\eps_t\}}
                                              \big]
 +\frac{1}{n}
 \\
 &\leq& 
 \dbE^{\dbP\otimes_t\dbP_n}\big[Y_{\t_n\we D^\eps_t}
                                                      -\eps \1_{\t_n<D^\eps_t}
                                              \big]
 +\frac{1}{n}
 \;\leq\; 
 \dbE^{\dbP\otimes_t\dbP_n}\big[Y_{t}-\eps\1_{\t_n<D^\eps_t}\big] + \frac{1}{n},
 \eeaa
where the last inequality follows from the $\ol\cE-$supermartingale property of $Y$. Note that 
 $$
 \dbE^{\dbP\otimes_t\dbP_n}[Y_t]
 =
 \dbE^{\dbP\otimes\dbP_n}\Big[\dbE^{\dbP\otimes_t\dbP_n}[Y_t|\cF_t]\Big]
 =
 \dbE^{\dbP}\Big[\dbE^{\dbP\otimes_t\dbP_n}[Y_t|\cF_t]\Big]=\dbE^\dbP[Y_t].
 $$
Then $\eps(\dbP\otimes\dbP_n)[\t_n<D^\eps_t]\leq \frac{1}{n}$, and it follows from the previous estimate that:
 \beaa
 \dbE^\dbP[Y_t] 
 &\le& 
 \dbE^{\dbP\otimes\dbP_n}\big[(X_{\t_n}-Y_{D^\eps_t})\1_{\{\t_n<D^\eps_t\}}
                                                  +Y_{D^\eps_t}
                                          \big] + \frac{1}{n}
 \\
 &\le& 
 C(\dbP\otimes\dbP_n)\big[\t_n<D^\eps_t]^\frac12
 + \dbE^{\dbP\otimes\dbP_n}[Y_{D^\eps_t}]
 + \frac{1}{n}
 \;\le\;
 \frac{C}{\sqrt{n\eps}} + \frac{1}{n} + \dbE^{\dbP\otimes\dbP_n}\big[Y_{D^\eps_t}\big],
 \eeaa
by the fact that $\sup_{t\in[0,T]}|X_t|$ and $Y_{D^\epsilon_t}\in[X_{D^\eps_t},X_{D^\eps_t}+\eps]$ are both in $\L^2(\cP)$. Finally, we obtain
 \beaa
 \ol\cE^\cP\big[Y_{D^\eps_t}-Y_t\big]
 \;\ge\; 
 \dbE^{\dbP\otimes\dbP_n}\big[Y_{D^\eps_t}-Y_t\big]
 \;\ge\;
 -\Big(\frac{C}{\sqrt{n\eps}}+\frac{1}{n}\Big)
 \;\longrightarrow\; 0
 &\mbox{as}&
 n\rightarrow \infty. 
 \eeaa
\ep

\begin{Lemma}\label{K indistinguishable}
The processes $\{K_t,t\in[0,T]\}$ and $\{K_{D^\eps_t},t\in[0,T]\}_t$ are indistinguishable.
\end{Lemma}

\proof
By the decomposition of Proposition \ref{prop DM}, we have
 $$
 Y_{D^\eps_t}-Y_t
 =
 +\int_t^{D^\eps_t} H_sdB_s
 -\int_t^{D^\epsilon_t} L|H_s|_1ds
 -K_{D^\epsilon_t}+K_t,\ t\in[0,T],\ \dbP_0\text{-a.s.}
 $$
Since $\ol\cE[Y_{D^\epsilon_t}-Y_t]=0$ by Lemma \ref{lem:YDesp-Y}, we may find a sequence $(\dbP_n)_{n\ge 1}\subset\cP$ such that
 $$
 -\frac{1}{n}
 \;\le\;
 \dbE^{\dbP_n}\big[Y_{D^\eps_t}-Y_t\big]
 \;\le\;
 -\dbE^{\dbP_n}\big[K_{D^\eps_t}-K_t\big]
 \;\le\;
 -\ul\cE\big[K_{D^\eps_t}-K_t\big].
 $$
Then, it follows from the non decrease of $K$ that $\ul\cE\big[K_{D^\eps_t}-K_t\big]=0$, and therefore
 \b*
 K_{D^\eps_t} \;=\; K_t,
 &\dbP_0\text{-a.s. for all}&
 t\in[0,T].
 \e*
Consequently, $\dbP_0[\O']=1$, where $\O':=\big\{K_{D^\epsilon_t}=K_t,\text{ for all }t\in[0,T]\cap \dbQ\big\}$. Further, for any $t\in[0,T)$, let $\{t_n\}_n\subset\dbQ$ and  $t_n\downarrow t$. Since $K$ is nondecreasing, we see that $K_t \le K_{D^\eps_t} \le K_{D^\eps_{t_n}} = K_{t_n}$ on $\O'$. Since $K$ inherits the RCLL property of $Y$, this shows that $K_{D^\eps_t}$ is right continuous on $\O'$, and implies that $\{K_t\}_t$ and $\{K_{D^\epsilon_t}\}_t$ are indistinguishable.
\ep

\vspace{5mm}

\no {\bf Proof of Proposition \ref{prop:dK on X=Y}}
(i) We first prove that $\int(Y-X)dK=0$, $\P_0-$a.s. From Lemma \ref{K indistinguishable}, we have $\dbP_0[\Lambda]=1$, where $\Lambda=\{\o:K_t(\o)=K_{D^\eps_t}(\o)\text{ for all }t\in[0,T]\}$. Next, consider the decomposition of the process $Y$ into a continuous and a purely discontinuous part $Y=Y^c+Y^d$. From the decomposition of Proposition \ref{prop DM} and the fact that $K$ is increasing, we deduce that $\dbP_0[\Lambda']=1$, where $\Lambda':=\{\o:\Delta Y^d_t(\o)\leq 0\text{ for all }t\in[0,T]\}$.

Now fix any $\o\in\Lambda\cap\Lambda'$. For any $t_0\in \{t:X_t(\o)<Y_t(\o)\}$, denote $2c:=Y_{t_0}(\o)-X_{t_0}(\o)>0.$ Since $Y(\o)$ is RCLL with negative jumps, and $X(\o)$ is continuous, there exists $\d$ such that for all $t\in(t_0-\d,t_0]$ we have $Y_t(\o)-X_t(\o)>c$, and
 $$
 t_0\text{ is an interior point of }(t_0-\d,D^c_{t_0-\d}(\o))\subset \{t:X_t(\o)<Y_t(\o)\}.
 $$
Further, it is easy to prove that $\{t:X_t(\o)<Y_t(\o)\}$ can be covered by a countable number of open intervals in the form of $(t_n,D^{\epsilon_n}_{t_n}(\o))$. Finally, we have
 $$
 0\leq\int_0^T 1_{\{t:X_t(\o)<Y_t(\o)\}}dK_t(\o)
 \leq 
 \sum_{n=1}^{\infty}(K_{D^{\eps_n}_{t_n}(\o)}(\o)-K_{t_n}(\o))=0.
 $$ 
(ii) We next prove that $Y$ and $K$ are continuous. Consider the decomposition $K=K^c+K^d$ into a continuous and a purely discontinuous part, and let us show that $\dbP_0\big[K^d_t=0\text{ for all }t\in[0,T]\big]=1$.

Since $K$ is previsible and $\Delta K^d_t=K_t-K_{t-}$, $\Delta K^d$ is also previsible. In the following we set $\inf\emptyset=\infty$. By Theorem 12.3 in Chapter VI of \cite{WilliamsRogers} (p. 333), we know that $\t^\delta=\inf\{t\in(0,T]:\Delta K^d_t>\delta\}$ is a previsible stopping time (defined in Definition 12.1 in Chapter VI of \cite{WilliamsRogers}), for all $\delta>0$. Then, by Theorem 12.6 in Chapter VI of \cite{WilliamsRogers}, $\t^\delta$ can be announced by a sequence of stopping time $\t_n$, i.e. $\t_n<\t^\delta$ and $\t_n\uparrow\t^\delta$, $\dbP_0$-a.s. Then, since $K_t$ and $K_{D_t^\epsilon}$ are indistinguishable by Lemma \ref{K indistinguishable}, it follows from the definition of $\t^\delta$ that $K^d_{D_{\t_n}^\eps}=K^d_{\t_n}<K^d_{\t^\d}$. Then, $\t_n\leq D_{\t_n}^\epsilon<\t^\delta$, $\dbP_0$-a.s. Hence
 \b*
 \dbP_0[\Omega_0]=1,
 &\mbox{where}&
 \Omega_0
 :=\big\{\t_n\uparrow\t^\delta~\mbox{and}~\t_n\le D_{\t_n}^\eps<\t^\delta\big\}.
 \e* 
For all $\o\in \Omega_0$, we can find a sequence $t_n$ such that $D_{\t_n}^\epsilon(\o)\leq t_n<\t^\delta(\o)$ and $Y_{t_n}(\o)\leq X_{t_n}(\o)+\epsilon$. Sending $n\rightarrow \infty$, we get $Y_{\t^\delta(\o)-}(\o)\leq X_{\t^\delta(\o)}(\o)+\epsilon$. So, $Y_{\t^\delta-}\leq X_{\t^\delta}+\epsilon$, $\dbP_0$-a.s. Choosing $\epsilon<\delta$, we see that, whenever $\t^\delta\leq T$, $Y_{\t^\delta}\le Y_{\t^\delta-}-\d<X_{\t^\delta}$, which is therequired contradiction. Hence $\t^\delta=\infty$ for all $\delta>0$, implying that $K^d=0$, $\dbP_0$-a.s.\\
(iii) We now show that $\t^*$ is an optimal stopping rule. The results of (i) and (ii) lead to $K_{\t^*}=0$ $\P_0$-a.s.. Recall the generalized Doob-Meyer decomposition in Proposition \ref{prop DM}. Take 
$\l^*$ such that $\|\l^*\|\leq L$ and $\l^*H=L|H|_1$. Then, by taking expectation under $\P_{\l^*}$, we obtain that
 $$
 Y_0
 \;=\;
 \E^{\P_{\l^*}}\big[Y_{\t^*}\big]
 \;=\;
 \E^{\P_{\l^*}}\big[X_{\t^*}\big].
 $$
The last equality is due to the definition of $\t^*$. Finally, it is clear that $Y_0=\ol\cE^\cP[X_{\t^*}]$. Hence, $\t^*$ is an optimal stopping rule.
\ep

\subsection{Reduction to a standard optimal stopping problem}

As a consequence of the decomposition in Proposition \ref{prop DM} together with Lemma \ref{K indistinguishable}, we obtain the following reduction.

\begin{Proposition}\label{optimal one prob}
There exists a probability $\dbP^*\in\cP$ such that
$$Y_t=\esup_{\t\in\cT_*^t}\dbE^{\dbP^*}[X_\t|\cF_t],\ \dbP_0\text{-a.s.}$$
In particular, there exists a $\dbP^*$-martingale $M^*$ such that $Y=Y_0+M^*-K$, $\dbP_0$-a.s.
\end{Proposition}
\proof
First, for any $\t\in\cT_*^t$ and $\dbP\in\cP$, we have $Y_t\geq \dbE^\dbP[X_\t|\cF_t]$, $\dbP_0$-a.s. Hence, $Y_t\geq \esup_{\t\in\cT_*^t}\dbE^{\dbP}[X_\t|\cF_t],\ \dbP_0\text{-a.s.}$

On the other hand, let $\l^*$ be defined by its $i-$th entry $L\ \text{sgn}(H_t)_i$. From Proposition \ref{prop DM}, we know that $(H\cdot B)-\int_0^\cdot L|H_s|_1ds$ is a $\dbP_{\l^*}$-martingale. Then, it follows from the decomposition of Proposition \ref{prop DM}, together with Lemma \ref{K indistinguishable}, that
 $$
 Y_t
 =
 \dbE^{\dbP_{\l^*}}\big[Y_{D^\eps_t}+K_{D^\eps_t}-K_t|\cF_t\big]
 =
 \dbE^{\dbP_{\l^*}}\big[Y_{D^\epsilon_t}|\cF_t\big]
 \le 
 \dbE^{\dbP_{\l^*}}\big[X_{D^\epsilon_t}|\cF_t\big]+\eps.
 $$
Since $D^\eps\rightarrow D_t:=\inf\{s\geq t:Y_{t-}=X_t\}$, as $\eps\rightarrow 0$, this implies that 
 $$
 Y_t
 \le
 \dbE^{\dbP_{\l^*}}\big[X_{D_t}|\cF_t\big]
 \le
 \esup_{\t\in\cT_*^t}\dbE^{\dbP_{\l^*}}[X_\t|\cF_t],
 ~~\dbP_0\text{-a.s.}
 $$
\ep

\subsection{The $\dbF-$adapted Snell envelop}

Given the continuity of $Y$ in Proposition \ref{prop:dK on X=Y}, we now reduce to an $\dbF-$adapted version.

\begin{Proposition}\label{prop:FSnell}
There is an $\F-$adapted pathwise continuous indistinguishable version of $Y$. 
\end{Proposition}

\proof Define
 \b*
 Y^\dbF_t=\dbE^{\dbP_0}[Y_t|\cF_t]~\mbox{for}~t\in[0,T]\cap\dbQ,
 &\mbox{and}&
 Y^\dbF_t
 :=
 \lim_{s\upuparrows t,s\in\dbQ}\tilde{Y}_s
 ~\mbox{for}~t\in[0,T]\setminus\dbQ.
 \e*
The last limit exists by the pathwise continuity of $Y$, $\dbP_0$-a.s., see Proposition \ref{prop:dK on X=Y}. Clearly, $Y^\dbF_t$ is $\cF_t$-measurable. Since $Y$ is $\dbF^*-$adapted, we have $\dbP_0\big[Y=Y^\dbF\text{ on } [0,T]\cap \dbQ\big]=1$, and by the pathwise continuity of $Y$, we deduce that $\dbP_0\big[Y=Y^\dbF\big]=1$. Hence $Y$ and $Y^\dbF$ are indistinguishable.
\ep

\vspace{5mm}

In this section, we consider the process $Y$ in its version of Proposition \ref{prop:FSnell}, which we call  the $\dbF-$adapted Snell envelop of $X$. We next define:
 \b*
 Z_\t(\o)
 :=
 \sup_{\th\in\cT_*}\ol\cE[X^{\t(\o),\o}_\th],
 &\mbox{for all}&
 \t\in\cT.
 \e*
Clearly, $Y_0=Z_0$. The main result of this subsection is the following.

\begin{Proposition}\label{prop: Y_tau=Z_tau}
Let $Y$ be the $\F-$Snell envelop of $X$. Then, $Y_\t=Z_\t$, $\dbP_0-$\mbox{a.s.} for all $\t\in\cT_*$
\end{Proposition}

In preparation for the proof of this result, we prove two lemmas.

\begin{Lemma}\label{reciprocal}
Let $\hat Y$ be a continuous $\dbF$-adapted process such that, for some $\dbP_0-$martingale $\hat M$ and nondecreasing process $\hat K$:
\\
{\rm (i)} $\hat Y=Y_0+\hat M_t-\max_{|\l|\leq L}\langle \hat M,\int_0^\cdot \l_sdB_s\rangle-\hat K,$ $\dbP_0-$a.s.
\\
{\rm (ii)} $\hat Y\ge X$, $\dbP_0-$-a.s.
\\
{\rm (iii)} $\int_0^T 1_{\{t:X_t<\hat Y_t\}}d\hat K_t=0$, $\dbP_0$-a.s.
\\
Then, $\hat Y_t=\esup_{\t\in\cT_*^t,\dbP\in\cP}\dbE^\dbP[X_\t|\cF_t]$, $\dbP_0$-a.s. 
\end{Lemma}

\proof By martingale representation and the Property (i), there exists $\hat H\in\mathbb{H}_{loc}$ such that $\hat Y=\hat Y_0+(\hat H\cdot B)-L\int_0^. |H_s|_1ds-\hat K$,$\dbP_0-$-a.s. By Girsanov theorem, $\hat M^\l:=\int_0^\cdot \hat H_sdB_s-\int_0^\cdot \l_s^T\hat H_sds$ is $\dbP_\l$-local martingale, and it follows from the previous decomposition that there exists increasing process $\hat K^\l$ such that
 \be\label{DM under Pl}
 \hat Y
 \;=\;
 \hat Y_0+\hat M^\l-\hat K^\l,
 &\dbP_0-\mbox{a.s.}&
 \ee
By the uniqueness of the Doob-Meyer decomposition, we deduce that $\hat M^\l$ is a $\dbP_\l$-martingale, and it follows from \eqref{DM under Pl} and Property (ii) that
 \b*
 \hat Y_t
 \;\ge\; 
 \dbE^{\dbP_\l}\big[\hat Y_\t\big|\cF_t\big]
 \;\ge\;
 \dbE^{\dbP_\l}\big[X_\t\big|\cF_t\big]
 &\text{for all}&
 \t\in\cT_*^t ,
 ~\dbP_\l\in\cP.
 \e*
Hence, $\hat Y_t\geq \esup_{\t\in\cT_*^t,\dbP\in\cP}\dbE^\dbP[X_\t|\cF_t]$. For the reverse inequality, consider the stopping time $D_t:=\inf\{s\geq t:\hat Y_s=X_s\}\in\cT_*^t$. Let $\l^*$ be the process defined by its $i-$th entry $L\;\text{sgn}(\hat H_i)$. Note that $\hat K^{\l^*}=\hat K$ in \eqref{DM under Pl}, and therefore
 \b*
 \hat Y_t=\dbE^{\dbP_{\l^*}}\big[\hat Y_{D_t}+\hat K_{D_t}-\hat K_t\big|\cF_t\big].
 \e*
By property (iii) and the definition of $D_t$, it follows that $\hat K_{D_t}=\hat K_t$, $\dbP_0$-a.s., so that
 $$
 \hat Y_t
 \;=\;
 \dbE^{\dbP_{\l^*}}\big[\hat Y_{D_t}\big|\cF_t\big]
 \;=\;
 \dbE^{\dbP_{\l^*}}\big[X_{D_t}\big|\cF_t\big]
 \;\le\;
 \esup_{\t\in\cT_*^t,\dbP\in\cP}\;\dbE^\dbP[X_\t|\cF_t].
 $$
\ep

\begin{Lemma}\label{lemma full}
Let $M$ be a pathwise continuous $\dbP_0$-martingale with $\dbE^{\dbP_0}\big[\sup_{t\in[0,T]}|M_t|\big]<\infty$. Then, there exists an $\dbF$-adapted indistinguishable version $\tilde{M}$ such that:
 \b*
 \dbP_0\big[\o: \tilde{M}^{\t,\o}~\mbox{is a}~\dbP_0-\mbox{martingale}\big]=1
 &\mbox{for all}&
 \t\in\cT.
 \e*
\end{Lemma}

\proof
{\bf 1.} Let $\tilde{M}_T:=M_T$, and for all $\o\in\O$:
 \b*
 \tilde{M}_s(\o)
 :=
 \dbE^{\dbP_0}\big[M^{s,\o}_T\big]
 ~\mbox{for}~s\in {\dbQ}
 &\mbox{and}&
 \tilde{M}_t(\o)
 :=
 \limsup_{s\uparrow t,s\in {\dbQ}}\tilde{M}_s(\o)
 ~\mbox{for}~t\in[0,T]\setminus {\dbQ}.
 \e*
Clearly, $\tilde{M}$ is $\dbF$-adapted, and $\dbP_0\big[M=\tilde{M}~\mbox{on}~\dbQ\big]=1$. Since $M$ is continuous, it is easy to verify that $\dbP_0[M=\tilde{M}]=1$, i.e. $\tilde{M}$ is an indistinguishable version of $M$.
\\
{\bf 2.} Denote $|\tilde M|^*_t:=\sup_{s\le t}|\tilde M_s|$, and  
 $$
 I_\t
 :=
 \Big\{\o\in \O: \tilde{M}_\t(\o)=\dbE^{\dbP_0}\big[\tilde{M}^{\t,\o}_T\big],
                       ~\dbE^{\dbP_0}\big[|\tilde M|^{*^{\t,\o}}_t\big]<\infty,
                       ~\mbox{and}~\dbP_0\{\tilde M^{\t,\o}\text{ is continuous}\}=1
 \Big\}.
 $$
Since $\tilde{M}$ and $M$ are indistinguishable, $\tilde{M}$ is a $\dbP_0$-martingale and $\dbP_0[I_\t]=1$. For $\eta\in\cT$ with $\eta\ge \tau$, we define a sequence of stopping times $\eta_n:=\frac{[2^n\eta]+1}{2^n}$. Note that $\eta_n$ only take rational values. By the tower property and the definition of $\tilde M_s$ for $s\in\dbQ$, we obtain for $\o\in I_\t$:
 $$
 \tilde{M}_\t(\o)
 =
 \dbE^{\dbP_0}[\tilde{M}_T^{\t,\o}]
 =
 \lim_{n\rightarrow\infty}
 \int_\O \dbE^{\dbP_0}[\tilde{M}_T^{\eta_n,\o\otimes_{\tau} \o'}]\dbP_0(d\o')
 =
 \lim_{n\rightarrow\infty}\dbE^{\dbP_0}[\tilde{M}_{\eta_n}^{\t,\o}].
 $$
Since $\dbE^{\dbP_0}\big[|\tilde M|^{*^{\t,\o}}_t\big]<\infty$, it follows that the family $\{\tilde{M}_{\eta_n}^{\t,\o}\}_{n\in\mathbb{N}}$ is $\dbP_0-$uniformly integrable. Then, it follows from the $\P_0$-a.s. pathwise continuity of $\tilde{M}^{\t,\o}$ that $
 \tilde{M}_\t(\o)=\lim_{n\rightarrow\infty}\dbE^{\dbP_0}[\tilde{M}_{\eta_n}^{\t,\o}]=\dbE^{\dbP_0}[\tilde{M}_{\eta}^{\t,\o}].$ By the arbitrariness of $\eta\in\cT$, this proves that $\tilde M^{\tau,\o}$ is a $\P_0-$martingale.
\ep

\vspace{5mm}

\no {\bf Proof of Proposition \ref{prop: Y_tau=Z_tau}}\quad
Notice that $Y\ge X$, $\dbP_0-$a.s., and by Propositions \ref{prop DM} and \ref{prop:dK on X=Y}, there exists $H\in\mathbb{H}_{loc}$ and nondecreasing previsible process $K$ such that, with $M:=(H\cdot B)$:
 \b*
 Y=Y_0+M-\max_{|\l|\leq L}\langle M,\int_0^\cdot \l_sdB_s\rangle-K
 &\mbox{and}&
 \int_0^T 1_{\{t:X_t<Y_t\}}dK_t=0,
 ~~\dbP_0-\mbox{a.s.}
 \e*
The process $M$ is a pathwise continuous $\dbP_0$-martingale with $\E^{\P_0}\big[\sup_{t\in[0,T]}|M_t|\big]<\infty$, by Proposition \ref{prop DM}. By Lemma \ref{lemma full}, we may consider $M$ as the version for which $M^{\t,\o}$ is a $\dbP_0-$martingale, for $\dbP_0-$a.e.  $\o\in\O.$

Let $T':=T-\t(\o)$, and define $\tilde{M}^{\t,\o}_t(\o'):=M^{\t,\o}_t(\o')-M_\t(\o)$ for $t\in[0,T']$. Then, $\tilde M^{\t,\o}$ is $\dbP_0$-martingale for $\dbP_0-$a.e. $\o\in\O$. We now observe that $(Y^{\t,\o},\tilde M^{\t,\o},K^{\t,\o})$ satisfies the following properties for $\dbP_0$-a.e. $\o\in\O$:
\\
{\rm (i)} $Y^{\t,\o}=Y_{\t(\o)}(\o)+\tilde M^{\t,\o}-\max_{|\l|\leq L}\langle \tilde M^{\t,\o},\int_0^\cdot \l_sdB_s\rangle-K^{\t,\o}$, on $[0,T']$, $\dbP_0-$a.s.
\\
{\rm (ii)} $Y^{\t,\o}\ge X^{\t,\o}$ on $[0,T']$, $\dbP_0-$a.s.
\\
{\rm (iii)} $\int_0^{T'} 1_{\{t:X_t^{\t,\o}<Y_t^{\t,\o}\}}dK^{\t,\o}_t=0$, $\dbP_0-$a.s.

Then, it follows from Lemma \ref{reciprocal} that $Y_\t=Z_\t$, $\dbP_0-$a.s. 
\ep

\section{Appendix: on $\ol\cE$-submartingales}

We say that a process $u$ is an $\ol\cE$-regular submartingale, if
 \b*
 u_t(\o) \;\le\; \ol\cE[u^{t,\o}_\t]
 &\mbox{for all}&
 (t,\o)\in[0,T]\times\O~~\mbox{and}~~\t\in\cT.
 \e*
The main result of this section is the following.

\begin{Proposition}\label{cE submartingale}
Let $u\in C^0_{2,\cP}(\Theta,\R)$ be a $\ol\cE$-regular submartingale. Then, there exists $\dbP^*\in\cP$ such that $u$ is a $\dbP^*$-submartingale.
\end{Proposition}

\proof {\bf 1.} We shall prove in Step 2 below that
 \be\label{vst-step2}
 \ol\cE\big[ u^{t,\o}_{s-t}\big]
 \;=\;
 \esup_{\dbP\in\cP}\dbE^\dbP[u_s|\cF_t]
 &\mbox{for}~~\dbP_0-\mbox{a.e.}~~\o\in\O,&
 \mbox{for all}~~t<s\ge T-t.
 \ee
Let $t^n_k:=kT2^{-n}$, $k\ge 0$, and $\mathbb{I}_n:=\{T\wedge t^n_k:k\ge 0\}$. Since $\cP$ is weakly compact and $u\in C^0_{2,\cP}(\Theta,\dbR)$, we deduce from \eqref{vst-step2} that, for all pair $(n,k)$ with $t^n_k\le T$, there exists $\dbP^{n,k}\in\cP$ such that $u_{t^n_k}\le \dbE^{\dbP^{n,k}}[u_{t^n_{k+1}}|\cF_{t^n_k}]$, $\dbP^0-$a.s.
Defining $\dbP^n:=\dbP^{n,0}\otimes_{t^n_1}\dbP^{n,1}\otimes_{t^n_2}\cdots$, this implies that
 \b*
 u_{t_i^n\wedge T}
 \;\le\; 
 \dbE^{\dbP^n}[u_{t_j^n\wedge T}|\cF_{t_i^n\wedge T}]
 &\dbP_0-\mbox{a.s. for all}&
 0\le i\le j\le n.
 \e*
Since $\cP$ is weakly compact, $\dbP^n$ converges weakly to some $\dbP^*\in\cP$, after possibly passing to a subsequence. Observe that, for $m\ge n$, we have $\dbE^{\dbP^n}\big[u_{t_j^n\wedge T}|\cF_{t_i^n\wedge T}\big]=\dbE^{\dbP^m}\big[u_{t_j^n\wedge T}|\cF_{t_i^n\wedge T}\big]\longrightarrow \dbE^{\dbP^*}\big[u_{t_j^n\wedge T}|\cF_{t_i^n\wedge T}\big]$, as $m\to\infty$, by the fact that $u\in C^0_{2,\cP}(\Theta,\R)$. Hence, $u_t\le \dbE^{\dbP^*}[u_s|\cF_t]$, $\dbP_0-$a.s. for all $t\le s\le T-t$ with $s,t\in\mathbb{I}_n$. By the density of $\mathbb{I}_n$ in $[0,T]$, we further conclude that $u$ is a $\dbP^*$-submartingale. 
\\
{\bf 2.} It remains to prove \eqref{vst-step2}. For $t\leq s$, define a process:
 \b*
 v^s_t:=\esup_{\dbP\in\cP}\dbE^\dbP[u_s|\cF_t]
 &\mbox{for}&
 t\in[0,T],~0\le s\le T-t.
 \e*
Similar to Lemma \ref{lattice}, we may check that the family $\{\dbE^\dbP[u_s|\cF_t];\dbP\in\cP\}$ satisfies the lattice property. Then, for  $t_1\leq t_2\leq s$, we have for all $\dbP\in\cP$:
 $$
 \dbE^\dbP[v^s_{t_2}|\cF_{t_1}]
 =
 \esup_{\dbP'\in\cP}\dbE^{\dbP\otimes_{t_2}\dbP'}[u_s|\cF_{t_1}]
 \leq 
 \esup_{\dbP'\in\cP}\dbE^{\dbP\otimes_{t_1}\dbP'}[u_s|\cF_{t_1}]
 =
 v^s_{t_1}.
 $$
proving that $v^s$ is $\dbP$-supermartingale on $[0,s]$ for all $\dbP\in\cP$. Similar to Lemma \ref{lem: RCLL version}, we may consider $v^s$ in its $\dbF^*-$adapted RCLL version. 

Following the line of argument in the proof of Proposition \ref{prop DM}, there exists $H^s\in\mathbb{H}_{loc}$ and increasing process $K^s$ such that
 \b*
 v^s
 =v^s_0+(H^s\cdot B)-L\int_0^.|H^s_r|_1dr-K^s,
 &\dbP^0-\mbox{a.s.}
 \e*
We next prove that $K^s\equiv 0$, $\dbP^0$-a.s. Indeed, assuming to the contrary that $\dbP^0[K^s_s>0]>0$, it follows that $\ul\cE[K^s_s]>0$. Following the line of argument in Lemma \ref{lem:Y supmartingale}, it can be checked that $\ol\cE[v^s_t]=\ol\cE[u_s]$ for all $t\leq s$. Then, since $v^s_s=u_s$, it follows from the previous decomposition that
 \b*
 \dbE^{\dbP}[v^s_0]
 \;\ge\; 
 \dbE^{\dbP}[u_s+K^s_s]
 \;\ge\; 
 \dbE^{\dbP}[u_s]+\ul\cE[K^s_s]
 &\mbox{for all}&
 \dbP\in\cP,
 \e*
and therefore $\ol\cE[v^s_t]>\ol\cE[u_s]$, which is the required contradiction. This reduces the decomposition of $v^s$ to:
 \b*
 v^s
 =v^s_0+(H^s\cdot B)-L\int_0^.|H^s_r|_1dr,
 &\dbP^0-\mbox{a.s.}
 \e*
so that, with $\l^s$ the process with $i-$th entry $L\;\text{sgn}(H_i)$, we obtain $v^s_t=\dbE^{\dbP^{\l^s}}[u_s|\cF_t]$. We finally prove that \eqref{vst-step2} holds true by following the line of argument in the proof of Proposition \ref{prop: Y_tau=Z_tau}.
\ep

\end{document}